\documentclass{article}

\usepackage{arxiv}

\usepackage[utf8]{inputenc} % allow utf-8 input
\usepackage[T1]{fontenc}    % use 8-bit T1 fonts
\usepackage{hyperref}       % hyperlinks
\usepackage{url}            % simple URL typesetting
\usepackage{booktabs}       % professional-quality tables
\usepackage{amsfonts,amssymb,amsthm,amsmath,amscd,mathrsfs,latexsym}
\usepackage{nicefrac}       % compact symbols for 1/2, etc.
\usepackage{microtype}      % microtypography
\usepackage{lipsum}
\usepackage{graphicx}
\graphicspath{ {./images/} }
\newcommand{\indep}{\perp \!\!\! \perp}

\usepackage{amsmath}

\usepackage[capitalize, noabbrev]{cleveref}

\newtheorem{theorem}{Theorem}[section]
\newtheorem{ex}{Example}
\newtheorem{dfn}{Definition}

\newtheorem{lemma}[theorem]{Lemma}

\usepackage{tikz} % Pour les DAG/SWIG
\usetikzlibrary{positioning, calc, shapes.geometric, shapes, shapes.multipart, arrows.meta, arrows, decorations.markings, external, trees}
\usetikzlibrary{shapes,decorations,arrows,calc,arrows.meta,fit,positioning}
\tikzset{
	-Latex,auto,node distance =1.5 cm and 1.5 cm,semithick,
	state/.style ={minimum width = 1.5 cm},
	point/.style = {inner sep=0.4cm,fill,node contents={}},
	bidirected/.style={Latex-Latex,dashed},
	el/.style = {inner sep=2pt, align=left, sloped}
}
\usepackage{subfig} % Pour les sous-figures
\usepackage{algorithm}
\usepackage{algpseudocode}% http://ctan.org/pkg/algorithmicx

\title{Efficient adjustment sets for time-dependent treatment effect estimation in nonparametric causal graphical models}

\author{
 David Adenyo \\
  Department of Social and Preventive Medicine\\
  Universit\'e Laval\\
   \And
 Mireille E. Schnitzer \\
  Faculty of Pharmacy \& Department of Social and Preventive Medicine\\
  Universit\'e de Montr\'eal\\
  Department of Epidemiology, Biostatistics, and Occupational Health \\
  McGill University \\
  \And
 David Berger \\
  Department of Computer Science\\
  Universit\'e de Montr\'eal \\
  \And
 Jason R. Guertin \\
  Department of Social and Preventive Medicine\\
  Universit\'e Laval\\
  \And
 Denis Talbot \\
  Department of Social and Preventive Medicine\\
  Universit\'e Laval\\
  \texttt{denis.talbot@fmed.ulaval.ca}
}

\begin{document}
\maketitle
\begin{abstract}
Criteria for identifying optimal adjustment sets yielding consistent estimation with minimal asymptotic variance  of average treatment effects in parametric and nonparametric models have recently been established. In a single treatment time point setting, it has been shown that the optimal adjustment set can be identified based on a causal directed acyclic graph alone. In a time-dependent treatment setting, previous work has established graphical rules to compare the asymptotic variance of estimators based on nested time-dependent adjustment sets. However, these rules do not always permit the identification of an optimal time-dependent adjustment set based on a causal graph alone. We extend those results by exploiting conditional independencies that can be read from the graph and demonstrate theoretically and empirically that our results can yield estimators with lower asymptotic variance than those allowed by previous results. We further show how our results allow for the identification of optimal adjustment sets based on a directed acyclic graph alone in the time-dependent treatment setting. 
\end{abstract}

% keywords can be removed

\section{Introduction}
One of the most important challenges for causal inference with observational data lies in the ability to control for confounding. In practice, investigators often identify and measure confounding covariates during study design and estimate the causal effect of treatment or exposure by controlling, or adjusting, for a subset of covariates using an appropriate adjustment method during data analysis. One strategy for selecting which variables to control first entails constructing a causal graph based on expert opinion and subject-matter knowledge. Then, the backdoor criterion can be used to identify sets of variables that, when adequately controlled, allow for the elimination of confounding bias \cite{JudeaPearl2000}. Adjustment sets that allow for consistent estimation of a causal effect are known as ``sufficient adjustment sets''.

For a given causal graph, several sufficient adjustment sets may exist. A popular strategy for choosing among them is to select a minimal sufficient adjustment set, that is, the sufficient adjustment set that includes the fewest variables \cite{Evans2012aa, Greenland:1999aa}. Such a strategy can be appealing to proponents of parsimonious models but offers no guarantee that the resulting treatment effect estimator has a lower variance than an estimator adjusting for more variables. Alternatively, recent work has focused on the identification of the sufficient adjustment set that results in a treatment effect estimator with the lowest possible asymptotic variance. Notably, Henckel et al~\cite{Henckel:2022aa} developed a criterion to identify the optimal sufficient adjustment set for estimating the effect of a single or joint exposure intervention on one or multiple outcomes. This criterion is based on a causal graph and assumes that the relations between the variables in the graph are linear. Rotnitzky and Smucler
\cite{Rotnitzky:2020aa} extended those results in several directions, for instance by showing that the criterion proposed by Henckel et al~\cite{Henckel:2022aa} is valid for a broad class of nonparametric estimators and without assuming linear relations between variables for estimating the average effect of an exposure on an outcome. They further provided a new graphical criterion for comparing adjustment sets in terms of the variance of estimators of the effect of a time-dependent treatment subject to time-dependent confounding \cite{Rotnitzky:2020aa}. An important limitation of the latter result is the apparent impossibility to identify an optimal adjustment set based on a causal graph alone in the time-dependent setting.

In this paper, we extend previous results by proposing an alternative definition of a sufficient time-dependent adjustment set and a novel criterion for comparing the variance of estimators based on different time-dependent adjustment sets. We show theoretically and empirically that employing our criterion can further decrease the variance of estimators of time-dependent treatment effects. We further show that our alternative definition of a sufficient time-dependent adjustment set allows for the identification of optimal time-dependent adjustment sets based on the graph alone. The rest of the paper is organized as follows. In Section 2 we review some background concepts, and present the assumed data structure and the notation. In Section 3 we provide a novel definition of time-dependent adjustment sets. Section 4 presents our theoretical results and numerical illustrations. Section 5 discusses possible applications of our results and potential extensions to our work. All proofs of results stated in the main text can be found in the Appendix.

\section{Background}

We begin by briefly reviewing some concepts of causal graph theory and introducing some notation. The unfamiliar reader is referred to other references for a more detailed introduction (e.g., \cite{JudeaPearl2000} or the appendix of \cite{Vander-Weele:2011aa}). Let $\mathcal{G}$ be a directed acyclic graph (DAG) with nodes $\boldsymbol{V}$ representing the variables, that is, a graph for which all edges are unidirectional arrows and for which there is no cycle. A path is defined as a sequence of variables connected by arrows, independently of their direction. A causal path is a path in which all arrows point in the same direction. Opposingly, non-causal paths are paths in which not all arrows point in the same direction. A variable is said to be a collider on a given path if two arrows point to this variable on this path. Denote by $P$ the unknown distribution of the variables represented in the graph, and by $pa_{\mathcal{G}}(v_k)$ the parents of $v_k$ according to the graph $\mathcal{G}$, that is, $pa_{\mathcal{G}}(v_k)$ is the set of variables in $\mathcal{G}$ that have arrows stemming from them and pointing into $v_k$. Descendants of a node $v_k$ are all the nodes connected to $v_k$ by a causal path that starts with $v_k$. We say that a DAG is a causal graph if all common causes of the variables depicted in the graph are also on the graph \cite{JudeaPearl2000} and if
\begin{align*}
    P = \prod_{v_k \in \boldsymbol{V}} f\left(v_k\mid pa_{\mathcal{G}}(v_k)\right),
\end{align*}
where $f(v_k \mid pa_{\mathcal{G}}(v_k))$ is the distribution of $v_k$ given its parents on the graph. 

Suppose that we have a study with $T+1$ follow-up times and $n$ individuals sampled from a population so that the time-ordered (or otherwise topologically sorted) data take the following form $\mathcal{O}=(\boldsymbol{X}_0,A_0,\boldsymbol{X}_1,A_1,\dots,\boldsymbol{X}_T,A_T,Y)$ where $\boldsymbol{A}=(A_0,...,A_T) \subset \boldsymbol{V}$  is a time-varying exposure or treatment,
$\boldsymbol{X}_t=(X_{t1},X_{t2},...,X_{tn_{xt}}), \ t=0,...,T$ is a vector of time-dependent covariates (with $n_{xt}$ the number of covariates at time $t$) and $Y \subset \boldsymbol{V}$ is the outcome of interest, measured at the end of the study. We assume that $A_t, t=0,...,T$, correspond to finite valued random variables. For time ordered sets of variables, we use an overbar to denote the history of that variable up to a given time point. For example $\bar{\boldsymbol{A}}_T = (A_0, ..., A_T)$. Let $Y^{\boldsymbol{a}}$ be the counterfactual outcome, that is, the outcome that would have been observed if, possibly contrary to fact, exposure had been $\boldsymbol{A}=\boldsymbol{a}$ where $\boldsymbol{a}=(a_0,...,a_T)$.

Based on a causal graph, it is possible to construct a single world intervention graph (SWIG) that represents the causal model (including counterfactual variables) under a hypothetical intervention
\cite{Richardson:2013aa}. As will be seen shortly, we will be using these SWIGs to identify sets of variables that meet the sequential conditional exchangeability assumption necessary to identify time-dependent causal effects. A SWIG can be constructed from a DAG by splitting each treatment node $A_0,...,A_T$ in two parts, one part for the observed value and one for the value set by the (hypothetical) intervention: $A_0 \mid a_0, ..., A_T \mid a_T$. All arrows entering into a node $A_t$ in the DAG also enter $A_t$ in the SWIG, whereas all arrows stemming from a node $A_t$ in the DAG stem from $a_t$ in the SWIG. Further, all descendants of a node $a_t$ in the SWIG are indexed by an $a_t$ superscript to emphasize that it is a counterfactual variable. Figure \ref{fig1art2} presents an example of a DAG and its associated SWIG.

\begin{figure}[ht]
\centering
\subfloat[DAG]{
\begin{tikzpicture}
	\node (Z01) at (-5, 0)  {$X_{01}$};
    \node (Z02) at (-3.5, 0)  {$X_{02}$};
	\node (A0) at (-3,  1)  {$A_0$};
	\node (Z11) at (-2, -0.5)  {$X_{11}$};
    \node (Z12) at (-0.5, -0.5)  {$X_{12}$};
    \node (A1) at (0,  1)  {$A_1$};	
	\node (Y)  at (2,  0)  {$Y$};
 %\path (A1) edge node[midway] {} (A2);
 \path (Z01) edge node[midway] {} (A0);
	%\draw[->] (Z01) to (A0);
	\path (Z01) edge node[midway] {} (Z02);
	\path (Z02) edge node[midway] {} (A1);
         \path (Z02) edge node[midway] {} (Z11);
         \path (A0) edge[bend left = 60] node[midway] {} (Y);
	%\draw[->, bend left = 60] (A0) to (Y);
        \path (A0) edge node[midway] {} (A1);
 	\path (A0) edge node[midway] {} (Z11);
	\path (Z11) edge node[midway] {} (A1);
 	\path (Z11) edge node[midway] {} (Z12);
   	\path (Z12) edge node[midway] {} (Y);
	\path (A1) edge node[midway] {} (Y);
\end{tikzpicture}}
\subfloat[SWIG]{
\begin{tikzpicture}
\node (Z01) at (-5, 0)  {$X_{01}$};
    \node (Z02) at (-3.5, 0)  {$X_{02}$};
	\node (A0) at (-3,  1)  {$A_0 \mid a_0$};
	\node (Z11) at (-2, -0.5)  {$X_{11}^{a_0}$};
    \node (Z12) at (-0.5, -0.5)  {$X_{12}^{a_0}$};
    \node (A1) at (0,  1)  {$A_1^{a_0} \mid a_1$};	
	\node (Y)  at (2,  0)  {$Y^{\boldsymbol{a}}$};
        \path (Z01) edge node[midway] {} (A0);
	\path (Z01) edge node[midway] {} (Z02);
	\path (Z02) edge node[midway] {} (A1);
         \path (Z02) edge node[midway] {} (Z11);
         \path (A0) edge[bend left = 60] node[midway] {} (Y);
        \path (A0) edge node[midway] {} (A1);
 	\path (A0) edge node[midway] {} (Z11);
	\path (Z11) edge node[midway] {} (A1);
 	\path (Z11) edge node[midway] {} (Z12);
   	\path (Z12) edge node[midway] {} (Y);
	\path (A1) edge node[midway] {} (Y);
\end{tikzpicture}}
\caption{A causal directed acyclic graph (DAG; a) and the corresponding single world intervention graph (SWIG; b)}
	     \label{fig1art2}
\end{figure}

Independence between variables can be read from a causal DAG or a SWIG using the concept of d-separation \cite{JudeaPearl2000}. Consider disjoint nodes, $V_j$ and $V_k$, and a set of nodes $\boldsymbol{V}_\ell$. The nodes $V_j$ and $V_k$ are d-separated by $\boldsymbol{V}_\ell$ under a graph $\mathcal{G}$ if all paths between $V_j$ and $V_k$ are blocked given $\boldsymbol{V}_\ell$ and we write $V_j \indep_{\mathcal{G}} V_k \mid \boldsymbol{V}_\ell$. A path is said to be blocked given $\boldsymbol{V}_\ell$ if a) it contains a non-collider that is an element of $\boldsymbol{V}_\ell$, or b) there is a collider on the path where neither it nor its descendants are in $\boldsymbol{V}_\ell$. In the following, we adopt the convention that sets are d-separated from empty sets (i.e., $V_j \indep_{\mathcal{G}} \emptyset$ and $V_j \indep_{\mathcal{G}} \emptyset \mid \boldsymbol{V}_\ell$ are always true). If $V_j$ and $V_k$ are d-separated given $\boldsymbol{V}_\ell$, then $V_j$ and $V_k$ are also independent conditional on $\boldsymbol{V}_\ell$.                

Our focus is on the efficient estimation of contrasts between counterfactual outcome expectations of the form $\Delta(P;\mathcal{G})\equiv \sum_{\boldsymbol{a} \in \mathcal{A}}c_{\boldsymbol{a}} \mathbb{E}_P(Y^{\boldsymbol{a}})$, where $\boldsymbol{a}$ is a vector taking values in the finite set $\mathcal{A}$ of all possible values that $\boldsymbol{A}$ can take and $c_{\boldsymbol{a}} \in \mathbb{R}$. For example, in the context of a two time point study with a binary treatment, $\boldsymbol{A} = (A_0,A_1)$ and $\mathcal{A} = \{(0,0),(1,1),(1,0),(0,1)\}$. If $c_{(0,0)}=-1$, $c_{(1,1)}=1$, and $c_{(1,0)}=c_{(0,1)} = 0$ then $\Delta(P;\mathcal{G})=\mathbb{E}_P(Y^{(1,1)})-\mathbb{E}_P(Y^{(0,0)})$. These causal contrasts can be identified from the observed data by first identifying the counterfactual expectations $\chi_{\boldsymbol{a}}(P;\mathcal{G}) \equiv \mathbb{E}_P(Y^{\boldsymbol{a}})$. \vspace{0.25cm}

\begin{dfn} \label{dfn1}
    Let ${\boldsymbol{K} = (\boldsymbol{K}_0, ..., \boldsymbol{K}_T) \subset \boldsymbol{X}}$ be time ordered, disjoint sets of variables. Similar to Rotnitzky and Smucler~\cite{Rotnitzky:2020aa}, we say that $\boldsymbol{K}$ is a sufficient time-dependent adjustment set if
\begin{align}
    \mathbb{E}_P(Y^{\boldsymbol{a}}) &= \mathbb{E}_P(\mathbb{E}_P[...\mathbb{E}_P\{\mathbb{E}_P(Y \mid \bar{\boldsymbol{A}}_T = \bar{\boldsymbol{a}}_T,\bar{\boldsymbol{K}}_T)\mid \bar{\boldsymbol{A}}_{T-1} = \bar{\boldsymbol{a}}_{T-1},\boldsymbol{\bar{K}}_{T-1}\}...\mid A_0 = a_0,\boldsymbol{K}_0]) \label{eq:id1art2} 
\end{align}
which can equivalently be written as
\begin{align}
    \mathbb{E}_P(Y^{\boldsymbol{a}}) &=
    \mathbb{E}_P\left\{\frac{I_{\boldsymbol{a}}(\boldsymbol{A})Y}{\prod_{t=0}^T {P}(A_t = a_t\mid \bar{\boldsymbol{A}}_{t-1} = \bar{\boldsymbol{a}}_{t-1}, \bar{\boldsymbol{K}}_t)}\right\}, \nonumber
\end{align}
where $I$ is the usual indicator function.
\end{dfn}

One way to identify a sufficient time-dependent adjustment set from a causal graph is to first construct a SWIG from the causal graph, then identify a set $\boldsymbol{K}_0, ..., \boldsymbol{K}_T$ such that the sequential conditional exchangeability assumption, ${Y^{\boldsymbol{a}} \indep A_t \mid  \bar{\boldsymbol{A}}_{t-1}, \bar{\boldsymbol{K}}_t}$ for $t = 0,...,T$, holds, where $\bar{A}_{-1} \equiv \emptyset$. Example \ref{ex:1} illustrates this process.

\begin{ex}[Part 1] \label{ex:1}
Assume that Figure \ref{fig1art2} represents a causal DAG and its associated SWIG. Using d-separation rules, it can be seen that $Y^{\boldsymbol{a}} \indep A_0 \mid  X_{02}$ and that $Y^{\boldsymbol{a}} \indep A_1^{a_0}\mid   X_{12}^{a_0}, X_{02}, A_0$, for example. Under consistency (i.e., if $A_t = a_t$ then $V^{a_t} = V$), the second conditional independence becomes ${Y^{\boldsymbol{a}} \indep A_1\mid  X_{12}, X_{02}, A_0}$. As such, $\mathbb{E}_P(Y^{\boldsymbol{a}}) = \mathbb{E}_P[\mathbb{E}_P\{\mathbb{E}_P(Y\mid A_1 = a_1, X_{12}, A_0 = a_0, X_{02})\mid A_0 = a_0, X_{0,2}\}]$ and $\boldsymbol{K} = (\boldsymbol{K}_0 = X_{02}, \boldsymbol{K}_1 = X_{12})$ is a time-dependent adjustment set.  
\end{ex}

\section{Proposal}

We now propose an alternative definition of a sufficient time-dependent adjustment set. As Rotnitzky and Smucler~\cite{Rotnitzky:2020aa} note, their definition of a sufficient time-dependent adjustment set ignores potential simplifications to the identification formula (\ref{eq:id1art2}) ensuing from conditional independencies. Our alternative definition permits exploiting these conditional independencies. 

\subsection{Extended definition of sufficient time-dependent adjustment sets}
\begin{dfn} \label{dfn2}
Let $\boldsymbol{Z}=(\boldsymbol{Z}_0,...,\boldsymbol{Z}_T)$, with $\boldsymbol{Z}_t \subseteq \{\bar{\boldsymbol{X}}_t, \bar{\boldsymbol{A}}_{t-1}\}$, $t = 0, ..., T$. We say that $\boldsymbol{Z}$ is a sufficient time-dependent adjustment set if 
\begin{align}
\mathbb{E}_P(Y^{\boldsymbol{a}}) &= \mathbb{E}_P(\mathbb{E}_P[...\mathbb{E}_P\{\mathbb{E}_P(Y\mid A_t = a_T,\boldsymbol{Z}_T(\bar{\boldsymbol{a}}_{T-1}))\mid A_{T-1} = a_{T-1},\boldsymbol{Z}_{T-1}(\bar{\boldsymbol{a}}_{T-2})\}...\mid A_0 = a_0,\boldsymbol{Z}_0]) \label{eq:id2art2}
\end{align}
or equivalently
\begin{align*}
    \mathbb{E}_P(Y^{\boldsymbol{a}}) &= \mathbb{E}_P\left\{\frac{I_{\boldsymbol{a}}(\boldsymbol{A})Y}{\prod_{t=0}^T {P}(A_t = a_t\mid \boldsymbol{Z}_t(\bar{\boldsymbol{a}}_{t-1}))}\right\}.
\end{align*}
where we use the notation $\boldsymbol{Z}_t(\bar{\boldsymbol{a}}_{t-1})$ to emphasize that if some or all past treatments are in $\boldsymbol{Z}_t(\bar{\boldsymbol{a}}_{t-1})$ then they are held fixed.
\end{dfn}
Definition \ref{dfn2} differs from Definition \ref{dfn1} in a few subtle but important aspects. First, unlike Definition \ref{dfn1}, Definition~\ref{dfn2} allows, but does not require, for previous treatments to be members of sufficient adjustment sets (i.e., $\bar{\boldsymbol{A}}_{t-1} \in \bar{\boldsymbol{Z}}_t$, $t = 1, ..., T$, is allowed but not required in Definition \ref{dfn2}, whereas $\bar{\boldsymbol{A}}_{t-1}$ is required to be conditioned on in Definition \ref{dfn1}). In addition, the expectations in the identification formula (\ref{eq:id2art2}) are only conditional on $A_t$ and separate $\boldsymbol{Z}_t$ at each time point, not on the histories $\bar{\boldsymbol{A}}_t$ and $\bar{\boldsymbol{K}}_t$, $t = 0, ..., T$, as in Definition \ref{dfn1}. Importantly, all sets that satisfy Definition \ref{dfn1} also satisfy Definition \ref{dfn2}, but the reverse is false. 

\subsection{Construction of a sufficient time-dependent adjustment set according to Definition \ref{dfn2} \label{s:3.2}}

In practice, a sufficient time-dependent adjustment set according to Definition \ref{dfn2} can be constructed by first identifying a time-dependent adjustment set according to Definition \ref{dfn1} and exploiting conditional independencies that can be read from the DAG to exclude unnecessary variables from the identification formula. Algorithm \ref{alg:1} below provides a structured procedure for achieving this. This algorithm employs a recursive procedure, where the problem of estimating the counterfactual expectation under a time-dependent treatment regime $\boldsymbol{A}=\boldsymbol{a}$ is represented as a nested sequence of estimation problems where each problem entails estimating a counterfactual expectation under a given point treatment. Before presenting Algorithm \ref{alg:1}, we present an example of this process with theoretical arguments supporting each step.

To simplify the presentation, we consider a two time point setting, but a similar recursive procedure can be applied for any number of time points. Let $\boldsymbol{K} = (\boldsymbol{K}_0, \boldsymbol{K}_1)$
be a sufficient adjustment set according to Definition \ref{dfn1} and $\mathcal{G}$ be a causal graph for the  distribution $P$ of 
$\{\boldsymbol{X}, \boldsymbol{A}, Y\}$. The counterfactual mean can thus be expressed as
\begin{align}
    \mathbb{E}_P\left(Y^{\boldsymbol{a}}\right) 
    &= \mathbb{E}_P\!\left[ 
        \mathbb{E}_P\!\left\{ 
            \mathbb{E}_P\!\left( 
                Y \,\middle|\, \boldsymbol{A} = \boldsymbol{a}, \bar{\boldsymbol{K}}_1
            \right) 
        \,\middle|\, A_0 = a_0, \boldsymbol{K}_0 
        \right\} 
    \right].
    \label{eq:id1}
\end{align}
This expression highlights two nested estimation problems, one for each of the inner conditional expectations. 

The first estimation problem is
\[
\mathbb{E}_P\!\left(Y \,\middle|\, \boldsymbol{A} = \boldsymbol{a}, \bar{\boldsymbol{K}}_1 \right).
\]
Let $\boldsymbol{B}_1$ and $\boldsymbol{G}_1$ be a partition of $\{A_0, \boldsymbol{K}_0, \boldsymbol{K}_1\}$ such that $Y \indep_{\mathcal{G}}\boldsymbol{B}_1 |\, \boldsymbol{G}_1, A_1$. Then this conditional expectation can be written as
\[
\mathbb{E}_P(Y \mid \boldsymbol{A} = \boldsymbol{a}, \bar{\boldsymbol{K}}_1) 
= \mathbb{E}_P(Y \mid A_1 = a_1, \boldsymbol{G}_1).
\]
Remark that if $A_0$ is part of $\boldsymbol{G}_1$, it is fixed at $A_0 = a_0$. As such, the problem of estimating $\mathbb{E}_P(Y \mid \boldsymbol{A} = \boldsymbol{a}, \bar{\boldsymbol{K}}_1)$ entails estimating a counterfactual expectation under the point-treatment $A_1 = a_1$ for each possible value that the variables in $\boldsymbol{G}_1$ can take. 

Defining $Q_0(a_1, \boldsymbol{G}_1) \equiv \mathbb{E}_P(Y \mid A_1 = a_1, \boldsymbol{G}_1)$, which we simply denote as $Q_0$ when no ambiguity arises, the second estimation problem can be written as
\begin{align*}
    \mathbb{E}(Q_0 \mid A_0 = a_0, \boldsymbol{K}_0).
\end{align*}
Let $\boldsymbol{B}_0$ and $\boldsymbol{G}_0$ be disjoint subsets of $\boldsymbol{K}_0$ such that $Q_0 \indep_{\mathcal{G}} \boldsymbol{B}_0 \mid \boldsymbol{G}_0$. The second estimation problem can therefore be written as
\[
\mathbb{E}_P\!\left( Q_0 \,\middle|\, A_0 = a_0, \boldsymbol{G}_0 \right).
\]
This again corresponds to the estimation of the counterfactual expectation under a point-treatment.  

A key step in the above procedure is the correct identification of $\boldsymbol{B}_0$ such that $Q_0 \indep_{\mathcal{G}} \boldsymbol{B}_0 \mid \boldsymbol{G}_0$. To simplify this task, we can construct a subgraph $\mathcal{G}_0$ based on the initial DAG $\mathcal{G}$. The idea is to simplify the structure of the initial DAG to focus on the problem of evaluating $Q_0 \indep_{\mathcal{G}} \boldsymbol{B}_0\mid \boldsymbol{G}_0$ without losing any information in the initial graph that is relevant to this problem. These subgraphs will also play an important role in identifying optimal time-dependent adjustment sets. The first step of the construction of $\mathcal{G}_0$ is to remove the node $Y$ (since it is no longer relevant) and to add a node $Q_0$. For a given $a_1$, remark that $Q_0$ is only a function of $\boldsymbol{G}_1$. As such, add arrows from all variables in $\boldsymbol{G}_1$ to $Q_0$. Next, because $A_1$ is fixed at $a_1$ for all individuals when calculating $Q_0$, no association between other nodes and $Q_0$ can pass through $A_1$. Consequently, the node $A_1$ and all arrows pointing into it or from it are removed from $\mathcal{G}_0$. Nodes $\boldsymbol{X}_1$ are also superfluous since our interest lies in $\{\boldsymbol{B}_0, \boldsymbol{G}_0\} = \boldsymbol{K}_0 \subseteq \boldsymbol{X}_0$. However, before removing these nodes, we need to verify if some of the association between $Q_0$ and $\boldsymbol{X}_0$ passes only through $\boldsymbol{X}_1$. To do so, we replace any causal path from variables in $\boldsymbol{X}_0$ to $Q_0$ passing only through variables in $\boldsymbol{X}_1$ as a direct arc from these variables to $Q_0$. The same process is applied for all causal paths from $A_0$ to $Q_0$ passing only through variables in $\boldsymbol{X}_1$. Variables $\boldsymbol{X}_1$ can now be removed from $\mathcal{G}_0$. The resulting subgraph can now be used to more easily evaluate relevant d-separation conditions and identify $\boldsymbol{B}_0$. We now present Algorithm \ref{alg:1}.

\begin{algorithm}[H]
\caption{Construction of a sufficient time-dependent adjustment set according to Definition~\ref{dfn2}. \label{alg:1}}
\begin{algorithmic}
    \State 1. Identify a sufficient time-dependent adjustment set $\boldsymbol{K}$ according to Definition \ref{dfn1} using the SWIG associated with the initial DAG. Express this sufficient time-dependent adjustment set as $\boldsymbol{Z}$ using the notation of Definition \ref{dfn2}. In the following steps, the DAG (and not the SWIG) will be used to read d-separation conditions.
    \State 2. Read from the initial DAG d-separation between $Y$ and variables in $\{A_T, \boldsymbol{Z}_T\}$ and exclude variables from $\boldsymbol{Z}_T$ that are d-separated from $Y$ by $A_T$ and the other variables in $\boldsymbol{Z}_T$. Denote as $\boldsymbol{G}_T$ the set of variables that are retained. 
    \State 3. Let $Q_{T} = Y$.
    \For{$t = T-1$ to $0$}
    \State 3.1 Construct a subgraph $\mathcal{G}_t$ following these steps in order:
    \begin{itemize}
        \item[i)] remove node $Q_{t+1}$ and add a novel node $Q_{t}$.
        \item[ii)] add arrows from all nodes in $\boldsymbol{G}_{t+1}$ to $Q_{t}$. 
        \item[iii)] remove the node $A_{t+1}$ and all arrows pointing into it or from it.
        \item[iv)] Change all causal paths from any variable in $\{\bar{A}_t, \bar{\boldsymbol{X}}_t\}$ to $Q_t$ that only pass through variables in $\boldsymbol{X}_{t+1}$ into an arrow from the starting variable to $Q_t$, that is, for all $k$, $0 \leq k \leq t$, change all paths of the form ${X}_{k \ell} \rightarrow X_{(t+1)j_1} \rightarrow \cdots \rightarrow X_{(t+1)j_m}\rightarrow Q_{t}$ to ${X}_{k \ell}\rightarrow Q_{t}$, and ${A}_{k} \rightarrow X_{(t+1)j_1} \rightarrow \cdots \rightarrow X_{(t+1)j_m}\rightarrow Q_{t}$ to ${A}_{k}\rightarrow Q_{t}$ with $ \ell \in \{1, \dots,n_{xk}\}$, $j_1,\dots,j_m\in \{1,2,\dots,n_{x(t+1)}\}, m \geq 1$.
       \item[v)] remove all nodes $\boldsymbol{X}_{t+1}$ and all arrows pointing into them or from them.
    \end{itemize}
    \State 3.2 Read from the subgraph $\mathcal{G}_t$ d-separation between $Q_{t}$ and variables in $\{A_t, \boldsymbol{X}_t\}$ and exclude variables from $\boldsymbol{X}_t$ that are d-separated from $Q_{t}$ by $A_t$ and the other variables in $\boldsymbol{X}_t$. Denote as $\boldsymbol{G}_t$ the set of variables that are retained.
\EndFor
\State 4. Return $\boldsymbol{G}=(\boldsymbol{G}_0,\dots,\boldsymbol{G}_T)$ a sufficient time-dependent adjustment set according to Definition~\ref{dfn2}.
\end{algorithmic}
\end{algorithm}

We return to Example \ref{ex:1} to illustrate the application of Algorithm \ref{alg:1} concretely.

\addtocounter{ex}{-1}
\begin{ex}[Part 2] \label{ex:1c}
Step 1 of the algorithm was achieved in the first part of Example \ref{ex:1}, where the sufficient time-dependent adjustment set $\boldsymbol{K} = (\boldsymbol{K}_0 = X_{02}, \boldsymbol{K}_1 = X_{12})$ was found. Using the notation introduced in Definition \ref{dfn2}, this sufficient time-dependent adjustment set can be expressed as ${\boldsymbol{Z} = (\boldsymbol{Z}_0 = X_{02}, \boldsymbol{Z}_1 = \{X_{12}, A_0, X_{02}\})}$, which leads to the same identification formula as previously, that is 
$$\mathbb{E}_P(Y^{\boldsymbol{a}}) = \mathbb{E}_P[\mathbb{E}_P\{\mathbb{E}_P(Y\mid A_1 = a_1, X_{12}, A_0 = a_0, X_{02})\mid A_0 = a_0, X_{02}\}].$$ 

In Step 2, by inspecting the DAG, it can be seen that $Y \indep_\mathcal{G} X_{02} \mid  A_1, X_{12}, A_0$. As such, 
\begin{align*}
    \mathbb{E}_p(Y^{\boldsymbol{a}}) &= \mathbb{E}_P[\mathbb{E}_P\{\mathbb{E}_P(Y\mid A_1 = a_1, X_{12}, A_0 = a_0, X_{02})\mid A_0 = a_0, X_{02}\}]\\ &
    = \mathbb{E}_P[\mathbb{E}_P\{\mathbb{E}_P(Y\mid A_1 = a_1, X_{12}, A_0 = a_0)\mid A_0 = a_0, X_{02}\}]
\end{align*}
and we deduce that  $(\boldsymbol{Z}_0 = X_{02}, \boldsymbol{G}_1 = \{X_{12}, A_0\})$ is a sufficient time-dependent adjustment set according to Definition~\ref{dfn2}. 

For $t = 0$, in Step 3.1, we produce the subgraph $\mathcal{G}_0$ in Figure \ref{fig2art2}. To obtain this graph, the following steps were followed: i) $Y$ has been removed and $Q_0=\mathbb{E}_P(Y\mid A_1 = a_1, X_{12}, A_0 = a_0)$ has been added, ii) arrows $X_{12} \to Q_0$ and $A_0 \to Q_0$ have been added, iii) $A_1$ and all arrows pointing into it or from it were removed, iv) an arrow $X_{02} \to Q_0$ was added because there was a causal path from it to $Q_0$ only passing through variables in $\boldsymbol{X}_1$ ($X_{02} \to X_{11} \to X_{12} \to Q_0$), and v) the nodes $X_{11}$ and $X_{12}$ were removed.
\begin{figure}[H]
\centering
\begin{tikzpicture}
	\node (Z01) at (-5, 0)  {$X_{01}$};
    \node (Z02) at (-1, -2)  {$X_{02}$};
	\node (A0) at (-2,  1)  {$A_0$};
	\node (Y)  at (2,  0)  {$Q_0$};	
	\path (Z01) edge node[midway] {} (A0);
	\path (Z01) edge node[midway] {} (Z02);
	\path (A0) edge node[midway] {} (Y);
        \path (Z02) edge node[midway] {} (Y);
\end{tikzpicture}
\caption{The subgraph $\mathcal{G}_0$ obtained in Step 3.1 of Example \ref{ex:1c}}
	     \label{fig2art2}
\end{figure}
In Step 3.2, we can see from Figure \ref{fig2art2} that $X_{02}$ is associated with $Q_0$ conditional on $A_0$. No additional simplification is thus possible. The resulting simplified time-dependent adjustment set according to Definition \ref{dfn2} is $\boldsymbol{G} = (\boldsymbol{G}_0 = X_{02}, \boldsymbol{G}_1 = \{X_{12}, A_0\})$. Remark that this time-dependent adjustment set is not admissible according to Definition \ref{dfn1}. 
\end{ex}

In the next section, we turn our attention to comparing the asymptotic variance of nonparametric efficient estimators of causal contrasts $\Delta(P;\mathcal{G})\equiv \sum_{\boldsymbol{a} \in \mathcal{A}}c_{\boldsymbol{a}} \mathbb{E}_P(Y^{\boldsymbol{a}})$ based on different sufficient time-dependent adjustment sets according to our proposed definition.

\section{Comparison of sufficient time-dependent adjustment sets}

We begin by defining some additional notation:
\begin{equation*}
b_{{a}_t}(\boldsymbol{Z}_t;P)\equiv
\mathbb{E}_P[...\mathbb{E}_P\{\mathbb{E}_P(Y\mid {a}_T,\boldsymbol{Z}_T(\bar{\boldsymbol{a}}_{T-1}))\mid {a}_{T-1},\boldsymbol{Z}_{T-1}(\bar{\boldsymbol{a}}_{T-2})\}...\mid {a}_t,\boldsymbol{Z}_t(\bar{\boldsymbol{a}}_{t-1})],
\end{equation*}
\begin{equation*}
\pi_{a_t}(\boldsymbol{Z}_t;P)\equiv P(A_t=a_t\mid \boldsymbol{Z}_t(\bar{\boldsymbol{a}}_{t-1})),
\end{equation*}
\begin{equation*}
\psi_{P,\boldsymbol{a}}(\boldsymbol{Z};\mathcal{G})=
\dfrac{I_{\bar{\boldsymbol{a}}_{T}}(\bar{\boldsymbol{A}}_{T})}{\lambda_{\bar{\boldsymbol{a}}_T}(\bar{\boldsymbol{Z}}_T;P)}\{Y-b_{{a}_T}(\boldsymbol{Z}_T;P)\}+\displaystyle{\sum_{t=0}^{T}}
\dfrac{I_{\bar{{\boldsymbol{a}}}_{t-1}}(\bar{\boldsymbol{A}}_{t-1})\{b_{{a}_t}(\boldsymbol{Z}_t;P)-b_{{a}_{t-1}}(\boldsymbol{Z}_{t-1};P)\}}{\lambda_{\bar{\boldsymbol{a}}_{t-1}}(\bar{\boldsymbol{Z}}_{t-1};P)}, 
\end{equation*}
where 
\begin{align*}
    {b_{a_{-1}}(\boldsymbol{Z}_{-1};P) \equiv \chi_{\boldsymbol{a}}(P;\mathcal{G})}, \\
    \lambda_{\bar{\boldsymbol{a}}_{t-1}}(\bar{\boldsymbol{Z}}_{t-1};P) \equiv \displaystyle{\prod_{k=0}^{t-1}} 	\pi_{a_k}(\boldsymbol{Z}_k;P) \text{ and} \\
    I_{\bar{\boldsymbol{a}}_{-1}}(\bar{\boldsymbol{A}}_{-1})\{\lambda_{\bar{\boldsymbol{a}}_{-1}}(\bar{\boldsymbol{Z}}_{-1};P)\}^{-1} \equiv 1.
\end{align*}
The quantity $\psi_{P,\boldsymbol{a}}(\boldsymbol{Z};\mathcal{G})$ is the efficient influence function of the counterfactual expectation $\mathbb{E}_P(Y^{\boldsymbol{a}})$ under the nonparametric model using the time-dependent adjustment set $\boldsymbol{Z}$ \cite{robins1995}. Similar to the Cram\'er-Rao bound, the variance of this efficient influence function scaled by a factor $1/n$ defines a lower bound for the asymptotic variance of nonparametric estimators of $\mathbb{E}_P(Y^{\boldsymbol{a}})$ for a given time-dependent adjustment set \cite{Hines:2022aa, AnastasiosTsiatis}. We further define $\sigma^{2}_{\boldsymbol{a},\boldsymbol{Z}}(P)$ and $\sigma^{2}_{\Delta,\boldsymbol{Z}}(P)$ as the asymptotic variances of nonparametric efficient estimators of the counterfactual mean $\chi_{\boldsymbol{a}}(P;\mathcal{G})$ and of causal contrast $\Delta(P;\mathcal{G})$ under the time-dependent adjustment set $\boldsymbol{Z}$. Note that the latter can be obtained from the former using the functional delta method \cite{Zepeda-Tello:2022aa} or tricks delineated in \cite{Kennedy:2022aa}. Using this notation, the following two lemmas compare the asymptotic variance of nonparametric efficient estimators of causal contrasts between nested sufficient time-dependent adjustment sets. 

\begin{lemma}(Inclusion of additional variables $(\boldsymbol{G},\boldsymbol{D})$ in a sufficient time-dependent adjustment set $\boldsymbol{Z} = (\boldsymbol{B},\boldsymbol{C})$)\label{lemma:1}
Consider a DAG \(\mathcal{G}\) with vertex set \(\boldsymbol{V}\). 
Let \(\boldsymbol{A}\subset \boldsymbol{V}\) and \(Y \in \boldsymbol{V} \setminus \boldsymbol{A} \). For \( t = T, \dots, 0 \), define \(\mathcal{G}_t \) as the subgraph of \(\mathcal{G}\) constructed by applying Algorithm \ref{alg:1}. Let $\boldsymbol{K}=(\boldsymbol{K}_0,\boldsymbol{K}_1,...,\boldsymbol{K}_T)$ be a sufficient time-dependent adjustment set according to Definition \ref{dfn1}
and let $\boldsymbol{B}_t, \boldsymbol{G}_t,\boldsymbol{C}_t,\boldsymbol{D}_t\subseteq \{\bar{\boldsymbol{A}}_{t-1}, \bar{\boldsymbol{K}}_t\}$, 
%$\boldsymbol{G}_t \subseteq \{\bar{\boldsymbol{A}}_{t-1},\bar{\boldsymbol{K}}_t\}$, $\boldsymbol{C}_t \subseteq \{\bar{\boldsymbol{A}}_{t-1},\bar{\boldsymbol{K}}_t\}$ and $\boldsymbol{D}_t \subseteq \{\bar{\boldsymbol{A}}_{t-1},\bar{\boldsymbol{K}}_t\}
$t=0,...,T$, be a partition of $\{\bar{\boldsymbol{A}}_{t-1},\bar{\boldsymbol{K}}_t\}$.
%such that $(\boldsymbol{G}_t,\boldsymbol{B}_t,  \boldsymbol{C}_t,\boldsymbol{D}_t) = \{\bar{\boldsymbol{A}}_{t-1},\bar{\boldsymbol{K}}_t\}$. 
Suppose $(\boldsymbol{B},\boldsymbol{C})=\{(\boldsymbol{B}_0,\boldsymbol{C}_0),(\boldsymbol{B}_1,\boldsymbol{C}_1),\dots,(\boldsymbol{B}_T,\boldsymbol{C}_T)\}$ is a sufficient time-dependent adjustment set according to Definition \ref{dfn2} that satisfies
\begin{equation}
	A_t\indep_{\mathcal{G}_t} \boldsymbol{G}_t, \boldsymbol{D}_t \ \mid   \ \boldsymbol{B}_t,\boldsymbol{C}_t, \ for \ t=0,...,T. \label{eq1.lemma1}
\end{equation}
Then, for any distribution $P$ compatible with the graph $\mathcal{G}$, $(\boldsymbol{G},\boldsymbol{B},  \boldsymbol{C},\boldsymbol{D})=\{(\boldsymbol{G}_0,\boldsymbol{B}_0,  \boldsymbol{C}_0,\boldsymbol{D}_0),(\boldsymbol{G}_1,\boldsymbol{B}_1,  \boldsymbol{C}_1,\boldsymbol{D}_1),\dots,(\boldsymbol{G}_T,\boldsymbol{B}_T,  \boldsymbol{C}_T,\boldsymbol{D}_T)\}$ is also a sufficient time-dependent adjustment set  according to Definition \ref{dfn2} and  
\begin{enumerate}
\item 
$
\sigma^{2}_{\boldsymbol{a},\boldsymbol{B},\boldsymbol{C}}(P)-\sigma^{2}_{\boldsymbol{a},\boldsymbol{G},\boldsymbol{B},  \boldsymbol{C},\boldsymbol{D}}(P)\geq 0 \ \ and  \ \  \sigma^{2}_{\Delta,\boldsymbol{B},\boldsymbol{C}}(P)-\sigma^{2}_{\Delta,\boldsymbol{G},\boldsymbol{B},  \boldsymbol{C},\boldsymbol{D}}(P)\geq 0.
$
\item For any $\boldsymbol{G}' \subseteq \boldsymbol{G}$ and $\boldsymbol{D}' \subseteq \boldsymbol{D}$, ($\boldsymbol{G}',\boldsymbol{B},\boldsymbol{C},\boldsymbol{D}'$) is also sufficient time-dependent adjustment set according to Definition \ref{dfn2} and 
\begin{equation*}
    \sigma^2_{\boldsymbol{a},\boldsymbol{B},\boldsymbol{C}}(P) \;\;\ge\;\; 
\sigma^2_{\boldsymbol{a},\boldsymbol{G}',\boldsymbol{B},\boldsymbol{C},\boldsymbol{D}'}(P) \;\;\ge\;\; 
\sigma^2_{\boldsymbol{a},\boldsymbol{G},\boldsymbol{B},\boldsymbol{C},\boldsymbol{D}}(P).
\end{equation*}
\end{enumerate}
\end{lemma}

Lemma 11 of Rotnitzky and Smucler~\cite{Rotnitzky:2020aa} and Lemma \ref{lemma:1} rest on the same basic mechanism: supplementing a valid adjustment set with so-called \emph{precision variables} that help predict the outcome but, conditional on the current adjustment set, are not associated with the treatment, cannot introduce bias and will not increase the asymptotic variance of the causal estimator. 

\begin{lemma}(Exclusion of variables $(\boldsymbol{B},\boldsymbol{D})$ from a sufficient time-dependent adjustment set $\boldsymbol{Z} = (\boldsymbol{G}, \boldsymbol{B},\boldsymbol{C},\boldsymbol{D})$) \label{lemma:2}
Consider a DAG \(\mathcal{G}\) with vertex set \(\boldsymbol{V}\). 
Let \(\boldsymbol{A}\subset \boldsymbol{V}\) and \(Y \in \boldsymbol{V} \setminus \boldsymbol{A} \). For \( t = T, \dots, 0 \), define \(\mathcal{G}_t \) as the subgraph of \(\mathcal{G}\) constructed by applying Algorithm \ref{alg:1}. Let $\boldsymbol{K}=(\boldsymbol{K}_0,\boldsymbol{K}_1,...,\boldsymbol{K}_T)$ be a sufficient time-dependent adjustment set according to Definition \ref{dfn1}
and let $\boldsymbol{B}_t, \boldsymbol{G}_t,\boldsymbol{C}_t,\boldsymbol{D}_t\subseteq \{\bar{\boldsymbol{A}}_{t-1}, \bar{\boldsymbol{K}}_t\}$, 
$t=0,...,T$, be a partition of $\{\bar{\boldsymbol{A}}_{t-1},\bar{\boldsymbol{K}}_t\}$. 
%such that $(\boldsymbol{G}_t,\boldsymbol{B}_t,  \boldsymbol{C}_t,\boldsymbol{D}_t) = \{\bar{\boldsymbol{A}}_{t-1},\bar{\boldsymbol{K}}_t\}$. 
Suppose $(\boldsymbol{G},\boldsymbol{B},  \boldsymbol{C},\boldsymbol{D})=\{(\boldsymbol{G}_0,\boldsymbol{B}_0,  \boldsymbol{C}_0,\boldsymbol{D}_0),(\boldsymbol{G}_1,\boldsymbol{B}_1,  \boldsymbol{C}_1,\boldsymbol{D}_1),\dots,(\boldsymbol{G}_T,\boldsymbol{B}_T,  \boldsymbol{C}_T,\boldsymbol{D}_T)\}$ is a sufficient time-dependent adjustment set according to Definition \ref{dfn2} that satisfies
\begin{equation}
Y\indep_{\mathcal{G}_T} \boldsymbol{B}_T,\boldsymbol{D}_T \ \mid  \ \boldsymbol{G}_T,\boldsymbol{C}_T,{A}_T \label{eq:lemma2.1}
\end{equation}
and
\begin{equation}
Q_{t-1}\indep_{\mathcal{G}_{t-1}} \boldsymbol{B}_{t-1},\boldsymbol{D}_{t-1} \ \mid  \ \boldsymbol{G}_{t-1},\boldsymbol{C}_{t-1}, {A}_{t-1}\ for \  t=T,...,1. \label{eq:lemma2.2}
\end{equation}
Then, for any distribution $P$ compatible with the graph $\mathcal{G}$, $(\boldsymbol{G},\boldsymbol{C})=\{(\boldsymbol{G}_0,\boldsymbol{C}_0),(\boldsymbol{G}_1,\boldsymbol{C}_1),\dots,(\boldsymbol{G}_T,\boldsymbol{C}_T)\}
$ is also a sufficient time-dependent adjustment set according to Definition \ref{dfn2} and
\begin{enumerate}
\item 
$\sigma^{2}_{\boldsymbol{a},\boldsymbol{G},\boldsymbol{B},\boldsymbol{C},\boldsymbol{D}}(P)-\sigma^{2}_{\boldsymbol{a},\boldsymbol{G},\boldsymbol{C}}(P)\geq 0 \ \ and  \ \  \sigma^{2}_{\Delta,\boldsymbol{G},\boldsymbol{B},\boldsymbol{C},\boldsymbol{D}}(P)-\sigma^{2}_{\Delta,\boldsymbol{G},\boldsymbol{C}}(P)\geq 0.
$
\item For any $\boldsymbol{B}' \subseteq \boldsymbol{B}$ and $\boldsymbol{D}' \subseteq \boldsymbol{D}$, ($\boldsymbol{G},\boldsymbol{B}',\boldsymbol{C},\boldsymbol{D}'$) is also sufficient time-dependent adjustment set according to Definition \ref{dfn2} and 
\begin{equation*}
\sigma^2_{\boldsymbol{a},\boldsymbol{G},\boldsymbol{B},\boldsymbol{C},\boldsymbol{D}}(P) 
\;\;\ge\;\; \sigma^2_{\boldsymbol{a},\boldsymbol{G},\boldsymbol{B}',\boldsymbol{C},\boldsymbol{D}'}(P)
\;\;\ge\;\; \sigma^2_{\boldsymbol{a},\boldsymbol{G},\boldsymbol{C}}(P).
\end{equation*}
\end{enumerate}
\end{lemma}
As in Lemma 12 of Rotnitzky and Smucler~\cite{Rotnitzky:2020aa}, Lemma \ref{lemma:2} states that retaining so-called superfluous adjustment variables that are associated with treatment but provide no additional predictive value for the outcome after conditioning on the existing adjustment set leads to an increase, or at least no reduction, in the asymptotic variance of the estimators.

The following theorem combines both of these lemmas together.

\begin{theorem} \label{theorem:1}
Consider a DAG \(\mathcal{G}\) with vertex set \(\boldsymbol{V}\). 
Let \(\boldsymbol{A}\subset \boldsymbol{V}\) and \(Y \in \boldsymbol{V} \setminus \boldsymbol{A} \). For \( t = T, \dots, 0 \), define \(\mathcal{G}_t \) as the subgraph of \(\mathcal{G}\) constructed by applying Algorithm \ref{alg:1}. Let $\boldsymbol{K}=(\boldsymbol{K}_0,\boldsymbol{K}_1,...,\boldsymbol{K}_T)$ be a sufficient time-dependent adjustment set according to Definition \ref{dfn1}
and let $\boldsymbol{B}_t, \boldsymbol{G}_t,\boldsymbol{C}_t,\boldsymbol{D}_t\subseteq \{\bar{\boldsymbol{A}}_{t-1}, \bar{\boldsymbol{K}}_t\}$, 
$t=0,...,T$, be a partition of $\{\bar{\boldsymbol{A}}_{t-1},\bar{\boldsymbol{K}}_t\}$. 
%such that $(\boldsymbol{G}_t,\boldsymbol{B}_t,  \boldsymbol{C}_t,\boldsymbol{D}_t) = \{\bar{\boldsymbol{A}}_{t-1},\bar{\boldsymbol{K}}_t\}$.
Suppose  
$(\boldsymbol{B},\boldsymbol{C})=\{(\boldsymbol{B}_0,\boldsymbol{C}_0),(\boldsymbol{B}_1,\boldsymbol{C}_1),\dots,(\boldsymbol{B}_T,\boldsymbol{C}_T)\}
$
and
$(\boldsymbol{G},\boldsymbol{C})=\{(\boldsymbol{G}_0,\boldsymbol{C}_0),(\boldsymbol{G}_1,\boldsymbol{C}_1),\dots,(\boldsymbol{G}_T,\boldsymbol{C}_T)\}
$
are two sufficient time-dependent adjustment sets according to Definition \ref{dfn2} such that 
\begin{equation}\label{eqv1}
    A_T\indep_{\mathcal{G}_T} \boldsymbol{G}_T,\boldsymbol{D}_T  \ \mid  \ \boldsymbol{B}_T,\boldsymbol{C}_T,
\end{equation}
\begin{equation}\label{eqv2}
    Y\indep_{\mathcal{G}_T} \boldsymbol{B}_T,\boldsymbol{D}_T \ \mid  \ \boldsymbol{G}_T,\boldsymbol{C}_T,{A}_T,
\end{equation}
and
\begin{equation}\label{eqv3}
    A_{t-1}\indep_{\mathcal{G}_{t-1}} \boldsymbol{G}_{t-1},\boldsymbol{D}_{t-1}  \ \mid  \ \boldsymbol{B}_{t-1},\boldsymbol{C}_{t-1},
\end{equation}
\begin{equation}\label{eqv4}
    Q_{t-1}\indep_{\mathcal{G}_{t-1}} \boldsymbol{B}_{t-1},\boldsymbol{D}_{t-1} \ \mid  \ \boldsymbol{G}_{t-1}, \boldsymbol{C}_{t-1},{A}_{t-1}, \  for \  t = T,\dots,1.
\end{equation}
Then for all distributions $P$ compatible with the causal graph $\mathcal{G}$,
\begin{enumerate}
\item 
$\sigma^{2}_{\boldsymbol{a},\boldsymbol{B},\boldsymbol{C}}(P)\geq\sigma^{2}_{\boldsymbol{a},\boldsymbol{G},\boldsymbol{B},\boldsymbol{C},\boldsymbol{D}}(P)\geq \sigma^{2}_{\boldsymbol{a},\boldsymbol{G},\boldsymbol{C}} \ \ and  \ \  \sigma^{2}_{\Delta,\boldsymbol{B},\boldsymbol{C}}(P)\geq\sigma^{2}_{\Delta,\boldsymbol{G},\boldsymbol{B},\boldsymbol{C},\boldsymbol{D}}(P)\geq\sigma^{2}_{\Delta,\boldsymbol{G},\boldsymbol{C}}(P).
$
\item For any $\boldsymbol{G}' \subseteq \boldsymbol{G}$, $\boldsymbol{B}' \subseteq \boldsymbol{B}$ and $\boldsymbol{D}' \subseteq \boldsymbol{D}$, ($\boldsymbol{G}',\boldsymbol{B},\boldsymbol{C},\boldsymbol{D}'$) and ($\boldsymbol{G},\boldsymbol{B}',\boldsymbol{C},\boldsymbol{D}'$) are also sufficient time-dependent adjustment set according to Definition \ref{dfn2}
and
\begin{equation*}
    \sigma^2_{\boldsymbol{a},\boldsymbol{B},\boldsymbol{C}}(P) \;\;\ge\;\; 
\sigma^2_{\boldsymbol{a},\boldsymbol{G}',\boldsymbol{B},\boldsymbol{C},\boldsymbol{D}'}(P) \;\;\ge\;\; 
\sigma^2_{\boldsymbol{a},\boldsymbol{G},\boldsymbol{B},\boldsymbol{C},\boldsymbol{D}}(P) \;\;\ge\;\; \sigma^2_{\boldsymbol{a},\boldsymbol{G},\boldsymbol{B}',\boldsymbol{C},\boldsymbol{D}'}(P)
\;\;\ge\;\; \sigma^2_{\boldsymbol{a},\boldsymbol{G},\boldsymbol{C}}(P), \ \ and
\end{equation*}
\begin{equation*}
    \sigma^2_{\Delta,\boldsymbol{B},\boldsymbol{C}}(P) \;\;\ge\;\; 
\sigma^2_{\Delta,\boldsymbol{G}',\boldsymbol{B},\boldsymbol{C},\boldsymbol{D}'}(P) \;\;\ge\;\; 
\sigma^2_{\Delta,\boldsymbol{G},\boldsymbol{B},\boldsymbol{C},\boldsymbol{D}}(P) \;\;\ge\;\; \sigma^2_{\Delta,\boldsymbol{G},\boldsymbol{B}',\boldsymbol{C},\boldsymbol{D}'}(P)
\;\;\ge\;\; \sigma^2_{\Delta,\boldsymbol{G},\boldsymbol{C}}(P).
\end{equation*}
\end{enumerate}
\end{theorem}

\subsection{Optimal sufficient time-dependent adjustment set}

We now extend a reasoning similar to Rotnitzky and Smucler\cite{Rotnitzky:2020aa} and Henckel et al~\cite{Henckel:2022aa} to the longitudinal case to demonstrate how to obtain an optimal sufficient time-dependent adjustment set. Our process is based on the subgraphs $\mathcal{G}_t$ built as described in Algorithm \ref{alg:1}. For each subgraph $\mathcal{G}_t, t=T,\dots,0$, we define the causal nodes $\mbox{cn}(A_t,Q_{t},\mathcal{G}_t)$ as the set of all nodes in $\mathcal{G}_t$ (excluding $A_t$) that lie on a causal path from $A_t$ to $Q_{t}$ with $Q_{T}=Y$ and $\mathcal{G}_T=\mathcal{G}$ the initial DAG. The corresponding forbidden set is given by $\mbox{forb}(A_t,Q_{t},\mathcal{G}_t)=\mbox{de}(\mbox{cn}(A_t,Q_{t},\mathcal{G}_t))\cup \{A_t\}$; the forbidden set thus includes $A_t$ as well as the descendants of the causal nodes. An optimal adjustment set for time $t=0,\dots,T$, denoted $\boldsymbol{O}_{t}(A_{t}, Q_{t}, \mathcal{G}_t)$, or simply as $\boldsymbol{O}_t$ whenever no confusion can arise, is then defined as: 
$$
\boldsymbol{O}_{t} \equiv \mbox{pa}_{\mathcal{G}_{t}}(\mbox{cn}(A_{t}, Q_{t}, \mathcal{G}_{t}))\backslash \mbox{forb}(A_{t}, Q_{t}, \mathcal{G}_t).
$$
An optimal sufficient time-dependent adjustment set is obtained by combining these sets: $\boldsymbol{\bar{O}} = (\boldsymbol{O}_0,...,\boldsymbol{O}_T)$.

The set $\boldsymbol{O}_{t}$ is an optimal adjustment set at time $t$ in the sense that it is a sufficient adjustment set for estimating the effect of $A_t$ on $Q_t$ and yielding estimators with the lowest possible asymptotic variance among all sufficient adjustment sets. This is established in Lemmas E.4 and E.5 of Henckel et al \cite{Henckel:2022aa}, as well as Corollary 8 in Rotnitzky and Smucler \cite{Rotnitzky:2020aa}. Because, as we argued in Section \ref{s:3.2}, the problem of estimating a time-dependent counterfactual mean can be seen as a sequence of estimation of point-treatment counterfactual means, it follows that $\boldsymbol{\bar{O}}$ is also a sufficient time-dependent adjustment set according to Definition \ref{dfn2}. By composition, $\boldsymbol{\bar{O}}$ further yields estimators with the smallest possible variance among all time-dependent adjustment sets. Indeed, for any $\boldsymbol{Z}$ that is a sufficient time-dependent adjustment set according to Definition \ref{dfn2}, we can use Theorem \ref{theorem:1} to show that the asymptotic variance of estimators based on $\bar{\boldsymbol{O}}$ is smaller or equal than the variance of estimators based on $\boldsymbol{Z}$. To show this, for all $t = 0, ..., T$, let $\boldsymbol{C}_t = \boldsymbol{O}_t \cap \boldsymbol{Z}_t$, $\boldsymbol{G}_t = \boldsymbol{O}_t \backslash \boldsymbol{C}_t$, $\boldsymbol{B}'_t = \boldsymbol{Z}_t \backslash \boldsymbol{C}_t$, $\boldsymbol{D}_t = \emptyset$, and choose any $\boldsymbol{B}_t \supseteq \boldsymbol{B}'_t$ disjoint of $(\boldsymbol{G}_t$, $\boldsymbol{C}_t)$ and such that $(\boldsymbol{G}, \boldsymbol{B}, \boldsymbol{C})$ is a sufficient adjustment set according to Definition \ref{dfn1}. It can then be verified that conditions (\ref{eqv1}), (\ref{eqv2}) (\ref{eqv3}) and (\ref{eqv4}) are met, which implies the desired result. Consequently, we have the following Theorem:

\begin{theorem}\label{th5.1}
Consider a DAG \(\mathcal{G}\) with vertex set \(\boldsymbol{V}\). 
Let \(\boldsymbol{A}\subset \boldsymbol{V}\) and \(Y \in \boldsymbol{V} \setminus \boldsymbol{A} \). For \( t = T, \dots, 0 \), define \(\mathcal{G}_t \) as the subgraph of \(\mathcal{G}\) constructed by applying Algorithm \ref{alg:1}. 
Then
\[
\boldsymbol{\bar{O}} = (\boldsymbol{O}_0, ..., \boldsymbol{O}_T)
\]
is a sufficient time-dependent adjustment set according to Definition \ref{dfn2} and, for all $P$ and all sets  $\boldsymbol{Z}$ that are sufficient time-dependent adjustment sets according to Definition \ref{dfn2},
\[
\sigma^2_{\boldsymbol{a},\boldsymbol{Z}}(P) - \sigma^2_{\boldsymbol{a},\boldsymbol{\bar{O}}}(P) \ge 0.
\]
$\boldsymbol{\bar{O}}$ is thus an optimal sufficient time-dependent adjustment set.
\end{theorem}

\subsection{Numerical illustrations}

We now provide two numerical illustrations of our theoretical results. We first return one last time to Example \ref{ex:1} to illustrate how our theorems can be used in practice to identify a time-dependent adjustment set yielding an estimator with reduced variance or an optimal adjustment set. Next, we revisit Example 4 of Rotnitzky and Smucler (2020) \cite{Rotnitzky:2020aa}. In this example, an optimal adjustment set cannot be identified based on the causal graph alone when considering Definition \ref{dfn1} and associated lemmas and theorem. Opposingly, our extended Definition \ref{dfn2} and associated theorems resolve this problem. 

\addtocounter{ex}{-1}
\begin{ex}[Part 3] \label{ex:1d}
Table \ref{tab:1} lists all time-dependent adjustment sets according to Definition \ref{dfn2}. Note that only adjustment sets 1-9 satisfy Definition \ref{dfn1}. Based on Theorem \ref{theorem:1}, the asymptotic variance of estimators based on sets 5 and 14 can be compared, for example. Let $\boldsymbol{K} = (\boldsymbol{K_0} = X_{02}, \boldsymbol{K_1} = X_{12}, A_0, X_{02})$. 
Because 
\begin{align*}
  A_1 &\indep_\mathcal{G} \emptyset \mid X_{02}, X_{12}, A_0 \\
  A_0 &\indep_{\mathcal{G}_0} \emptyset \mid X_{02} \\
  Y &\indep_\mathcal{G} X_{02} \mid  X_{12}, A_0, A_1 \\
  X_{02} &\indep_{\mathcal{G}_0} \emptyset \mid X_{02}, A_0
\end{align*}
we can apply Theorem \ref{theorem:1} with $\boldsymbol{B} = (\boldsymbol{B_0} = \emptyset, \boldsymbol{B_1} = X_{02})$, $\boldsymbol{C} = (\boldsymbol{C_0} = X_{02}, \boldsymbol{C_1} = X_{12}, A_0)$, $\boldsymbol{D} = \emptyset$, $\boldsymbol{G} = \emptyset$. 
%we know by theorem \ref{theorem:1} that $\sigma^{2}_{\boldsymbol{a},\boldsymbol{G},\boldsymbol{B},\boldsymbol{C},\boldsymbol{D}}(P)\geq \sigma^{2}_{\boldsymbol{a},\boldsymbol{G},\boldsymbol{C}}$, which is equivalent to $\sigma^{2}_{a,(\{X_{01},X_{02}\},\{X_{02},X_{11},X_{12},X_{01},A_0\})}(P)\geq\sigma^{2}_{a, (\{X_{02}\},\{X_{12},A_0\})}(P)$. 
Consequently, the asymptotic variance of nonparametric efficient estimators based on set 14 is lower or equal to the asymptotic variance of nonparametric efficient estimators based on set 5. We can further show that set 14 is an optimal sufficient time-dependent adjustment set. Indeed, we have $\text{cn}(A_1, Y, \mathcal{G}) = Y$, $\text{forb}(A_1, Y, \mathcal{G}) = A_1$, $\boldsymbol{O}_{1} = pa_{\mathcal{G}}(\text{cn}(A_1, Y, \mathcal{G})) \backslash \text{forb}(A_1, Y, \mathcal{G}) = \{A_0, X_{12}\}$ at $t = 1$, and $\text{cn}(A_0, Q_0, \mathcal{G}_0) = Q_0$, $\text{forb}(A_0, Q_0, \mathcal{G}_0) = A_0$, $\boldsymbol{O}_{0} = pa_{\mathcal{G}_0}(\text{cn}(A_0, Q_0, \mathcal{G}_0)) \backslash \text{forb}(A_0, Q_0, \mathcal{G}_0) = X_{02}$ at $t = 0$. We provide a numerical illustration based on simulated data to support these theoretical results. We simulated $10,000$ datasets of size $n =  1,000$ using the following data generating equations:
\begin{align*}
    X_{01} &\sim TN(\mu = 0, \sigma^2 = 1, \min = -2, \max = 2) \\
    X_{02} &\sim X_{01} + TN(\mu = 0, \sigma^2 = 1, \min = -2, \max = 2) \\
    A_0 &\sim Bernouilli(p = \text{expit}(X_{01})) \\
    X_{11} &\sim A_0 + X_{02} + TN(\mu = 0, \sigma^2 = 1, \min = -2, \max = 2) \\
    X_{12} &\sim X_{11} + TN(\mu = 0, \sigma^2 = 1, \min = -2, \max = 2) \\
    A_1 &\sim Bernouilli(p = \text{expit}(A_0+X_{02}+X_{11})) \\
    Y &\sim A_1 + A_0 + X_{12} + TN(\mu = 0, \sigma^2 = 1, \min = -2, \max = 2), 
\end{align*}
where $TN$ refers to a truncated normal distribution (in order to avoid practical positivity violations) and $\text{expit}(\cdot) = \exp(\cdot)/(1 + \exp(\cdot))$ is the inverse of the $\text{logit}$ function. The number of replications was chosen to ensure that the Monte Carlo standard error for the standard deviation was $\approx$ 0.001, which is sufficiently small to accurately identify the adjustment set with the lowest variance. We estimated $\mathbb{E}(Y^{11})$ using a pooled longitudinal targeted maximum likelihood estimator, which is based on the nonparametric efficient influence function $\psi_{P,\boldsymbol{a}}(\boldsymbol{Z};\mathcal{G})$. This was achieved in R v4.3.2 using the \texttt{ltmle} package v2.0.0 and generalized linear models to model both the outcome and exposure. 

\begin{table}[ht]
\caption{Sufficient time-dependent adjustment sets in Example \ref{ex:1} according to Definition \ref{dfn2} \label{tab:1}}
\centering
\begin{tabular}{cccc}
Adjustment set & $\boldsymbol{Z_0}$         & $\boldsymbol{Z_1}$                                & Monte Carlo\\
               &                        &                                               & standard deviation\\  \hline
1              & $X_{01}$              & $\{X_{11}, A_0, X_{01}\}$                   & 0.125 \\
2              & $X_{01}$              & $\{X_{12}, A_0, X_{01}\}$                   & 0.106 \\
3              & $X_{01}$              & $\{X_{11}, X_{12}, A_0, X_{01}\}$          & 0.110 \\
4              & $X_{02}$              & $\{X_{11}, A_0, X_{02}\}$                   & 0.118 \\
5              & $X_{02}$              & $\{X_{12}, A_0, X_{02}\}$                   & 0.099 \\
6              & $X_{02}$              & $\{X_{11}, X_{12}, A_0, X_{02}\}$          & 0.101 \\
7              & $\{X_{01}, X_{02}\}$ & $\{X_{11}, A_0, X_{01}, X_{02}\}$          & 0.123 \\
8              & $\{X_{01}, X_{02}\}$ & $\{X_{12}, A_0, X_{01}, X_{02}\}$          & 0.104 \\
9              & $\{X_{01}, X_{02}\}$ & $\{X_{11}, X_{12}, A_0, X_{01}, X_{02}\}$ & 0.106 \\ \hline
10             & $X_{01}$              & $\{X_{11}, A_0\}$                            & 0.117 \\
11             & $X_{01}$              & $\{X_{12}, A_0\}$                            & 0.099 \\
12             & $X_{01}$              & $\{X_{11}, X_{12}, A_0\}$                   & 0.106 \\
13             & $X_{02}$              & $\{X_{1,1}, A_0\}$                            & 0.107 \\
14             & $X_{02}$              & $\{X_{12}, A_0\}$                            & 0.088 \\
15             & $X_{02}$              & $\{X_{11}, X_{12}, A_0\}$                   & 0.095 \\
16             & $\{X_{01}, X_{02}\}$ & $\{X_{11}, A_0, X_{02}\}$                   & 0.123 \\
17             & $\{X_{01}, X_{02}\}$ & $\{X_{12}, A_0, X_{02}\}$                   & 0.104 \\
18             & $\{X_{0,1}, X_{02}\}$ & $\{X_{11}, X_{12}, A_0, X_{02}\}$          & 0.106 \\
19             & $\{X_{01}, X_{02}\}$ & $\{X_{11}, A_0, X_{01}\}$                   & 0.120 \\
20             & $\{X_{01}, X_{02}\}$ & $\{X_{12}, A_0, X_{01}\}$                   & 0.100 \\
21             & $\{X_{01}, X_{02}\}$ & $\{X_{11}, X_{12}, A_0, X_{01}\}$          & 0.105 \\
22             & $\{X_{01}, X_{02}\}$ & $\{X_{11}, A_0\}$                            & 0.112 \\
23             & $\{X_{01}, X_{02}\}$ & $\{X_{12}, A_0\}$                            & 0.093 \\
24             & $\{X_{01}, X_{02}\}$ & $\{X_{11}, X_{12}, A_0\}$                   & 0.100
\end{tabular}
\end{table}
\end{ex}
\clearpage

\begin{ex} \label{ex:2}
As a second example, we revisit Example 4 from \cite{Rotnitzky:2020aa}. The associated causal graph is represented in Figure \ref{fig4art2}. We assume the following time (or topological) ordering of variables $\{H, A_0, R, Q, A_1, Y\}$. In this example, the optimal sufficient adjustment set cannot be uniquely identified from the causal graph alone when considering Definition~\ref{dfn1}. Indeed, as illustrated in  \cite{Rotnitzky:2020aa}, the optimal sufficient adjustment set is either 1 or 8 depending on the data-generating equations, more specifically on the relative strength of the $H \rightarrow R \rightarrow Q$ and $H \rightarrow A_1$ pathways.

\begin{figure}[ht]
\centering
\begin{tikzpicture}
	\node (A0) at (-3,  0)  {$A_0$};
	\node (R)  at (-2, 0)  {$R$};
    \node (H)  at (-4, 1)  {$H$};
	\node (Q)  at (-1,  -1)  {$Q$};	
 	\node (A1) at (0,  0)  {$A_1$};
   	\node (Y) at (1,  0)  {$Y$};
	\path (A0) edge node[midway] {} (R);
 	\path (H)  edge node[midway] {} (R);
   	\path (R)  edge node[midway] {} (Q);
    \path (H)  edge node[midway] {} (A1);
    \path (Q)  edge node[midway] {} (Y);
    \path (A1)  edge node[midway] {} (Y);
\end{tikzpicture}
\caption{Causal directed acyclic graph for Example \ref{ex:2} (Example 4 in Rotnitzky and Smucler, 2020)}
	     \label{fig4art2}
\end{figure}

Table \ref{tab:2} lists all sufficient adjustment sets according to Definition \ref{dfn2}. Only sets 1-11 satisfy Definition \ref{dfn1}. Using Theorem \ref{theorem:1}, it can be found that estimators based on set 24 have a lower variance than those based on either sets 1 or 8. We first compare sets 8 and 24 by applying Theorem \ref{theorem:1},  letting $\boldsymbol{K}= (\boldsymbol{K_0} = H, \boldsymbol{K_1} = \{A_0, Q, H\})$, noticing that
\begin{align*}
    A_1 &\indep_\mathcal{G} \emptyset \mid A_0, H, Q \\
    A_0 &\indep_{\mathcal{G}_0} H \\
    Y &\indep_\mathcal{G} A_0, H \mid Q, A_1 \\
    Q_0 &\indep_{\mathcal{G}_0} \emptyset \mid H, A_0
\end{align*}
and thus setting $\boldsymbol{B} = (\boldsymbol{B_0} = \emptyset, \boldsymbol{B_1} = A_0, H)$, $\boldsymbol{C} = (\boldsymbol{C_0} = \emptyset, \boldsymbol{C_1} = Q)$, $\boldsymbol{D} = \emptyset$, $\boldsymbol{G} = (\boldsymbol{G}_0 = H, \boldsymbol{G}_1 = \emptyset)$. The comparison between sets 1 and 24 then follows from this previous result by setting $\boldsymbol{B}'=(\boldsymbol{B}_0 = \emptyset, \boldsymbol{B}_1 = A_0)$. Set 24 is moreover an optimal time-dependent adjustment set. At $t = 1$, we have $\text{cn}(A_1, Y, \mathcal{G}) = Y$, $\text{forb}(A_1, Y, \mathcal{G}) = A_1$, $\boldsymbol{O}_{1} = pa_{\mathcal{G}}(\text{cn}(A_1, Y, \mathcal{G})) \backslash \text{forb}(A_1, Y, \mathcal{G}) = Q$. The subgraph $\mathcal{G}_0$ can thus be shown to be $A_0 \rightarrow Q_0 \leftarrow H$. Therefore, at $t = 0$ we have $\text{cn}(A_0, Q_0, \mathcal{G}_0) = Q_0$, $\text{forb}(A_0, Q_0, \mathcal{G}_0) = A_0$, $\boldsymbol{O}_{0} = pa_{\mathcal{G}_0}(\text{cn}(A_0, Q_0, \mathcal{G}_0)) \backslash \text{forb}(A_0, Q_0, \mathcal{G}_0) = H$. 

To illustrate empirically these results, we generated 20,000 datasets of size $n$ = 1,000 according to the following equations
\begin{align*}
    A_0 &\sim Bernoulli(p = \text{expit}(0.5)) \\
    H &\sim TN(\mu = 0, \sigma^2 = 1, \min = -2, \max = 2) \\ 
    R &\sim A_0 + 2.5H + TN(\mu = 0, \sigma^2 = 1, \min = -2, \max = 2) \\ 
    Q &\sim R + TN(\mu = 0, \sigma^2 = 1, \min = -2, \max = 2) \\ 
    A_1 &\sim Bernoulli(p = \text{expit}(3H)) \\
    Y &\sim A_1 + Q + TN(\mu = 0, \sigma^2 = 1, \min = -2, \max = 2), 
\end{align*}
\noindent to illustrate the case with a stronger $H \rightarrow A_1$ pathway, and another 20,000 datasets of size $n$ = 1,000 where the data-generating equation of $R$ was modified to
\begin{align*}
    R &\sim A_0 + 4H + TN(\mu = 0, \sigma^2 = 1, \min = -2, \max = 2) 
\end{align*}
\noindent to illustrate the case with a stronger $H \rightarrow R \rightarrow Q$ pathway. The Monte Carlo standard error is $\leq 0.002$. It can be seen that the set with the lowest variance indeed depends on the scenario when only considering sets 1-11, but not when additionally considering sets 12-26. 

\begin{table}[ht]
\caption{Sufficient time-dependent adjustment sets in Example \ref{ex:2} according to Definition \ref{dfn2} \label{tab:2}}
\centering
\begin{tabular}{ccccc}
Adjustment set & $\boldsymbol{Z_0}$         & $\boldsymbol{Z_1}$       & \multicolumn{2}{c}{Monte Carlo SD} \\
               &                        &                      & Stronger $H \rightarrow A_1$ & Stronger $H \rightarrow R \rightarrow Q$ \\
               &                        &                      & pathway & pathway \\ \hline
1              & $\emptyset$          & $\{A_0, Q\}$         & 0.146 & 0.199 \\
2              & $\emptyset$          & $\{A_0, R\}$         & 0.175 & 0.222 \\
3              & $\emptyset$          & $\{A_0, H\}$         & 0.214 & 0.248 \\
4              & $\emptyset$          & $\{A_0, Q, R\}$      & 0.154 & 0.203 \\
5              & $\emptyset$          & $\{A_0, Q, H\}$      & 0.165 & 0.207 \\
6              & $\emptyset$          & $\{A_0, R, H\}$      & 0.191 & 0.229 \\
7              & $\emptyset$          & $\{A_0, Q, R, H\}$   & 0.165 & 0.207 \\
8              & $H$                    & $\{A_0, Q, H\}$      & 0.148 & 0.173 \\
9              & $H$                    & $\{A_0, R, H\}$      & 0.177 & 0.199 \\
10             & $H$                    & $\{A_0, R, Q, H\}$   & 0.148 & 0.173 \\
11             & $H$                    & $\{A_0, H\}$         & 0.201 & 0.220 \\ \hline
12             & $\emptyset$          & $Q$                  & 0.139 & 0.194 \\
13             & $\emptyset$          & $R$                  & 0.158 & 0.213 \\
14             & $\emptyset$          & $\{Q, R\}$           & 0.143 & 0.197 \\
15             & $\emptyset$          & $\{Q, H\}$           & 0.165 & 0.208 \\
16             & $\emptyset$         & $\{R, H\}$           & 0.191 & 0.229 \\
17             & $\emptyset$          & $\{Q, R, H\}$        & 0.165 & 0.208 \\
18             & $H$                    & $\{Q, H\}$           & 0.148 & 0.174 \\
19             & $H$                    & $\{R, H\}$           & 0.177 & 0.199 \\
20             & $H$                    & $\{R, Q, H\}$        & 0.148 & 0.174 \\
21             & $H$                    & $\{A_0, Q\}$         & 0.128 & 0.164 \\
22             & $H$                    & $\{A_0, R\}$         & 0.159 & 0.191 \\
23             & $H$                    & $\{A_0, R, Q\}$      & 0.136 & 0.168 \\
24             & $H$                    & $Q$                  & 0.119 & 0.157 \\
25             & $H$                    & $R$                  & 0.140 & 0.180 \\
26             & $H$                    & $\{R, Q\}$           & 0.123 & 0.161 \\
\end{tabular}\\
SD = Standard deviation
\end{table}
\end{ex}

\clearpage

\section{Conclusion}
In this article, we built on the work of Rotnitzky and Smucler \cite{Rotnitzky:2020aa} and proposed an alternative definition of sufficient time-dependent adjustment sets that takes into account potential simplifications to the identification formula using conditional independencies that can be read from the causal graph. We proposed two lemmas and a theorem that allow comparing the asymptotic variance of efficient estimators based on our definition of sufficient time-dependent adjustment sets and showed that further variance reduction can be obtained as compared to estimators based on previous results. We have further demonstrated that it is possible to identify an optimal sufficient time-dependent adjustment set based on the graph alone when considering our definition. We provided two numerical illustrations of our results. In one of those examples, the optimal adjustment set varies according to the data-generating equations for the same DAG when employing Rotnitzky and Smucler's definition \cite{Rotnitzky:2020aa}. In contrast, among the sufficient time-dependent adjustment sets that meet our definition, the optimal set remains the same regardless of the data-generating equations for the same DAG. The reduction in the variability of estimates allowed by our results was observed to be appreciable: the standard deviation of the estimates obtained using our optimal adjustment sets was between 10\% and 20\% lower than the standard deviation using the optimal sets according to previous results.

Our current contribution has several important implications. For example, our results can be used directly by data-analysts when estimating the effect of a time-dependent treatment in order to improve the precision of their estimates. Moreover, our results can be used for developing data-driven variable selection procedures. Indeed, several confounder selection approaches for a time-fixed covariate have been developed for targeting the optimal adjustment set \cite{shortreed, koch}. Extensions to the time-varying treatment setting based on the results of Rotnitzky and Smucler \cite{Rotnitzky:2020aa} and ours are already underway \cite{schnitzer2024adaptive}. 

Our current work could be extended in several directions. First, we have only considered DAGs where all variables are observed. In situations where some DAG vertices are unobservable, but some adjustment sets are observable, Henckel et al \cite{Henckel:2022aa} and Rotnitzky and Smucler \cite{Rotnitzky:2020aa} have provided theoretical results indicating that it is not always possible to find a uniformly optimal set among observable adjustment sets. It would interesting to similarly extend our results to situations where not all variables are observed. In studies with time-dependent confounding, it is also very common for data from some individuals to be missing or incomplete due to loss to follow-up or other factors (i.e., censored data). Ignoring these losses to follow-up when identifying the adjustment set may induce selection bias. Correa et al \cite{Correa} proposed a necessary and sufficient graphical criterion to estimate causal effects under confounding and selection bias for the cross-sectional case. Our results directly extend to the case where the interest lies specifically in estimating the effect of a treatment regime under no censoring, by simply treating censoring nodes as treatment nodes. Finally, to facilitate the use of our proposed rules, it would be greatly valuable to develop efficient algorithms and implement them in open-source software, similar to how DAGitty can be used to identify a minimal sufficient adjustment set based on a causal graph \cite{textor}, or how the \texttt{optAdjSet} function of the \texttt{pcalg} package in R can be used to identify an optimal adjustment set based on a causal graph in the single treatment time point setting. 

\section*{Acknowledgements}
The authors would like to thank Geneviève Lefebvre for her comments on earlier version of this manuscript. This work was supported by a Pierre Jacob Durand -- Citoyens du monde scholarship from Université Laval to DA; a Discovery Grant from the Natural Sciences and Engineering Research Council of Canada
to DT (RGPIN-2023-04911); Fonds de recherche du Québec -- Santé research career awards to DT (\url{https://doi.org/10.69777/312198}) and JRG (\url{https://doi.org/10.69777/330678}); and a Canada Research Chair to MES.  

\newpage
\bibliographystyle{unsrt}  
\bibliography{references}  %%% Remove comment to use the external .bib file (using bibtex).

\newpage

\appendix

To simplify notation, and whenever no confusion can arise, we will denote $\boldsymbol{L}_t = (\boldsymbol{D}_t,\boldsymbol{G}_t)$ and 
$\boldsymbol{N}_t = (\boldsymbol{B}_t,\boldsymbol{C}_t), t=0,\dots,T$
in Lemma \ref{lemma:1} and $\boldsymbol{M}_t = (\boldsymbol{D}_t,\boldsymbol{B}_t)$ and 
$\boldsymbol{U}_t = (\boldsymbol{G}_t,\boldsymbol{C}_t), t=0,\dots,T$ in Lemma \ref{lemma:2}.
\section*{Proof of Lemma 4.1} 
%\footnotesize
First, because $A_t\indep \boldsymbol{L}_t \ \mid   \ \boldsymbol{N}_t, \ \text{for} \ t=0,...,T$,
	\begin{eqnarray} 
\pi_{a_t}(\bar{\boldsymbol{a}}_{t-1},\bar{\boldsymbol{K}}_t;P)
 &=& \pi_{a_t}(\boldsymbol{L}_t,\boldsymbol{N}_t;P) \nonumber \\
		&\equiv& P(A_t=a_t\mid \boldsymbol{L}_t,\boldsymbol{N}_t) \nonumber \\
  &=& P(A_t=a_t\mid \boldsymbol{N}_t) \nonumber \\
		&=& 	\pi_{a_t}(\boldsymbol{N}_t;P). \label{eq:proof4.1.1}
	\end{eqnarray}
Because $\boldsymbol{N}$ is a sufficient time-dependent adjustment set and as a consequence of the previous result
	\begin{eqnarray} 
\mathbb{E}_p(Y^{\boldsymbol{a}})
         &=&
	\mathbb{E}_P\left[ \left\{ \displaystyle{\prod_{t=0}^{T}}\pi_{a_t}(\bar{\boldsymbol{a}}_{t-1},\bar{\boldsymbol{K}}_t;P)\right\}^{-1} I_{\boldsymbol{a}}(\boldsymbol{A})Y\right] \nonumber \\
		&=&
		\mathbb{E}_P\left[ \left\{ \displaystyle{\prod_{t=0}^{T}}\pi_{a_t}(\boldsymbol{N}_t;P)\right\}^{-1} I_{\boldsymbol{a}}(\boldsymbol{A})Y\right] \nonumber \\
		&=&\mathbb{E}_P\left[  \left\{ \displaystyle{\prod_{t=0}^{T}}\pi_{a_t}(\boldsymbol{L}_t,\boldsymbol{N}_t;P)\right\}^{-1} I_{\boldsymbol{a}}(\boldsymbol{A})Y\right]\nonumber.
	\end{eqnarray}
This shows that $(\boldsymbol{G,B,C,D})$ is also a sufficient time-dependent adjustment set. 

We next consider the efficient influence function $\psi_{P,\boldsymbol{a}}(\boldsymbol{N};\mathcal{G})$ in order to compare the asymptotic variance of nonparametric estimators based on based $\boldsymbol{N}$ versus $(\boldsymbol{N}, \boldsymbol{L})$. First, using some algebraic manipulations, we have 
{\footnotesize
\begin{eqnarray} 
    \psi_{P,\boldsymbol{a}}(\boldsymbol{N};\mathcal{G})
	&=& 
	\dfrac{I_{\bar{\boldsymbol{a}}_{T}}(\bar{\boldsymbol{A}}_{T})}{\lambda_{\bar{\boldsymbol{a}}_T}(\bar{\boldsymbol{N}}_T;P)}\{Y-b_{\boldsymbol{a}_T}(\boldsymbol{N}_T;P)\}+\displaystyle{\sum_{t=0}^{T}}
	\dfrac{I_{\bar{\boldsymbol{a}}_{t-1}}(\bar{\boldsymbol{A}}_{t-1})\{b_{\boldsymbol{a}_t}(\boldsymbol{N}_t;P)-b_{\boldsymbol{a}_{t-1}}(\boldsymbol{N}_{t-1};P)\}}{\lambda_{\bar{\boldsymbol{a}}_{t-1}}(\bar{\boldsymbol{N}}_{t-1};P)} \nonumber \\ %%%% Par dfn 
	&=& 
	\dfrac{I_{\bar{\boldsymbol{a}}_{T}}(\bar{\boldsymbol{A}}_{T})Y}{\lambda_{\bar{\boldsymbol{a}}_T}(\bar{\boldsymbol{N}}_T;P)}-\dfrac{I_{\bar{\boldsymbol{a}}_{T}}(\bar{\boldsymbol{A}}_{T})}{\lambda_{\bar{\boldsymbol{a}}_T}(\bar{\boldsymbol{N}}_T;P)}b_{\boldsymbol{a}_T}(\boldsymbol{N}_T;P) \nonumber \\
    &+&
    \displaystyle{\sum_{t=0}^{T}}\dfrac{I_{\bar{\boldsymbol{a}}_{t-1}}(\bar{\boldsymbol{A}}_{t-1})\{b_{\boldsymbol{a}_t}(\boldsymbol{N}_t;P)-b_{\boldsymbol{a}_{t-1}}(\boldsymbol{N}_{t-1};P)\}}{\lambda_{\bar{\boldsymbol{a}}_{t-1}}(\bar{\boldsymbol{N}}_{t-1};P)} \nonumber \\ %%%% Distribution du premier terme
	&=& 
	\dfrac{I_{\bar{\boldsymbol{a}}_{T}}(\bar{\boldsymbol{A}}_{T})Y}{\lambda_{\bar{\boldsymbol{a}}_T}(\bar{\boldsymbol{N}}_T;P)}-\displaystyle{\sum_{t=0}^{T}}\dfrac{I_{\bar{\boldsymbol{a}}_{t-1}}(\bar{\boldsymbol{A}}_{t-1})}{\lambda_{\bar{\boldsymbol{a}}_{t-1}}(\bar{\boldsymbol{N}}_{t-1};P)}\left\{\dfrac{I_{{\boldsymbol{a}}_{t}}(\boldsymbol{A}_{t})}{\pi_{a_t}(\boldsymbol{N}_t;P)}-1\right\}b_{\boldsymbol{a}_t}(\boldsymbol{N}_t;P) -b_{\boldsymbol{a}_{-1}}(\boldsymbol{N}_{-1};P)  \nonumber. %%%% Réécriture
\end{eqnarray}}
Next, using (\ref{eq:proof4.1.1}) and recalling that $b_{\boldsymbol{a}_{-1}}(\boldsymbol{N}_{-1};P) = \chi_{\boldsymbol{a}}(P;\mathcal{G})$, we get
\begin{eqnarray}
\psi_{P,\boldsymbol{a}}(\boldsymbol{N};\mathcal{G})
	&=& 
    \dfrac{I_{\bar{\boldsymbol{a}}_{T}}(\bar{\boldsymbol{A}}_{T})Y}{\lambda_{\bar{\boldsymbol{a}}_T}(\bar{\boldsymbol{L}}_T,\bar{\boldsymbol{N}}_T;P)}-\displaystyle{\sum_{t=0}^{T}}\dfrac{I_{\bar{\boldsymbol{a}}_{t-1}}(\bar{\boldsymbol{A}}_{t-1})}{\lambda_{\bar{\boldsymbol{a}}_{t-1}}(\bar{\boldsymbol{L}}_{t-1},\bar{\boldsymbol{N}}_{t-1};P)}\left\{\dfrac{I_{{\boldsymbol{a}}_{t}}(\boldsymbol{A}_{t})}{\pi_{a_t}(\boldsymbol{L}_t,\boldsymbol{N}_t;P)}-1\right\}b_{\boldsymbol{a}_t}(\boldsymbol{N}_t;P) \nonumber \\
    &-& \chi_{\boldsymbol{a}}(P;\mathcal{G}). \nonumber
\end{eqnarray} %%%% On remplace tous les pi(B) par pi(G,B)
Adding and subtracting 
$\displaystyle{\sum_{t=0}^{T}}\dfrac{I_{\bar{\boldsymbol{a}}_{t-1}}(\bar{\boldsymbol{A}}_{t-1})}{\lambda_{\bar{\boldsymbol{a}}_{t-1}}(\bar{\boldsymbol{L}}_{t-1},\bar{\boldsymbol{N}}_{t-1};P)}\left\{\dfrac{I_{{\boldsymbol{a}}_{t}}(\boldsymbol{A}_{t})}{\pi_{a_t}(\boldsymbol{L}_t,\boldsymbol{N}_t;P)}-1\right\}b_{\boldsymbol{a}_t}(\boldsymbol{L}_t,\boldsymbol{N}_t;P)$, we have:
%\footnotesize
{\footnotesize
\begin{eqnarray}
\psi_{P,\boldsymbol{a}}(\boldsymbol{N};\mathcal{G})
	&=& 
    \dfrac{I_{\bar{\boldsymbol{a}}_{T}}(\bar{\boldsymbol{A}}_{T})Y}{\lambda_{\bar{\boldsymbol{a}}_T}(\bar{\boldsymbol{L}}_T,\bar{\boldsymbol{N}}_T;P)}-\displaystyle{\sum_{t=0}^{T}}\dfrac{I_{\bar{\boldsymbol{a}}_{t-1}}(\bar{\boldsymbol{A}}_{t-1})}{\lambda_{\bar{\boldsymbol{a}}_{t-1}}(\bar{\boldsymbol{L}}_{t-1},\bar{\boldsymbol{N}}_{t-1};P)}\left\{\dfrac{I_{{\boldsymbol{a}}_{t}}(\boldsymbol{A}_{t})}{\pi_{a_t}(\boldsymbol{L}_t,\boldsymbol{N}_t;P)}-1\right\}b_{\boldsymbol{a}_t}(\boldsymbol{N}_t;P) \nonumber \\
    &-& \chi_{\boldsymbol{a}}(P;\mathcal{G}) + \displaystyle{\sum_{t=0}^{T}}\dfrac{I_{\bar{\boldsymbol{a}}_{t-1}}(\bar{\boldsymbol{A}}_{t-1})}{\lambda_{\bar{\boldsymbol{a}}_{t-1}}(\bar{\boldsymbol{L}}_{t-1},\bar{\boldsymbol{N}}_{t-1};P)}\left\{\dfrac{I_{{\boldsymbol{a}}_{t}}(\boldsymbol{A}_{t})}{\pi_{a_t}(\boldsymbol{L}_t,\boldsymbol{N}_t;P)}-1\right\}b_{\boldsymbol{a}_t}(\boldsymbol{L}_t,\boldsymbol{N}_t;P)  \nonumber \\
	&-& 	\displaystyle{\sum_{t=0}^{T}}\dfrac{I_{\bar{\boldsymbol{a}}_{t-1}}(\bar{\boldsymbol{A}}_{t-1})}{\lambda_{\bar{\boldsymbol{a}}_{t-1}}(\bar{\boldsymbol{L}}_{t-1},\bar{\boldsymbol{N}}_{t-1};P)}\left\{\dfrac{I_{{\boldsymbol{a}}_{t}}(\boldsymbol{A}_{t})}{\pi_{a_t}(\boldsymbol{L}_t,\boldsymbol{N}_t;P)}-1\right\}b_{\boldsymbol{a}_t}(\boldsymbol{L}_t,\boldsymbol{N}_t;P)  \nonumber \\		
	&=& \dfrac{I_{\bar{\boldsymbol{a}}_{T}}(\bar{\boldsymbol{A}}_{T})Y}{\lambda_{\bar{\boldsymbol{a}}_T}(\bar{\boldsymbol{L}}_T,\bar{\boldsymbol{N}}_T;P)}-\displaystyle{\sum_{t=0}^{T}}\dfrac{I_{\bar{\boldsymbol{a}}_{t-1}}(\bar{\boldsymbol{A}}_{t-1})}{\lambda_{\bar{\boldsymbol{a}}_{t-1}}(\bar{\boldsymbol{L}}_{t-1},\bar{\boldsymbol{N}}_{t-1};P)}\left\{\dfrac{I_{{\boldsymbol{a}}_{t}}(\boldsymbol{A}_{t})}{\pi_{a_t}(\boldsymbol{L}_{t},\boldsymbol{N}_t;P)}-1\right\}b_{\boldsymbol{a}_t}(\boldsymbol{L}_{t},\boldsymbol{N}_t;P)\nonumber \\ 
    &-& \chi_{\boldsymbol{a}}(P;\mathcal{G})+\displaystyle{\sum_{t=0}^{T}}\dfrac{I_{\bar{\boldsymbol{a}}_{t-1}}(\bar{\boldsymbol{A}}_{t-1})}{\lambda_{\bar{\boldsymbol{a}}_{t-1}}(\bar{\boldsymbol{N}}_{t-1};P)}\left\{\dfrac{I_{{\boldsymbol{a}}_{t}}(\boldsymbol{A}_{t})}{\pi_{a_t}(\boldsymbol{N}_t;P)}-1\right\}\{b_{\boldsymbol{a}_t}(\boldsymbol{L}_{t},\boldsymbol{N}_t;P)-b_{\boldsymbol{a}_t}(\boldsymbol{N}_{t};P)\} \nonumber \\
	&=& \psi_{P,\boldsymbol{a}}(\boldsymbol{L},\boldsymbol{N};\mathcal{G})+\sum_{t=0}^{T}\dfrac{I_{\bar{\boldsymbol{a}}_{t-1}}(\bar{\boldsymbol{A}}_{t-1})}{\lambda_{\bar{\boldsymbol{a}}_{t-1}}(\bar{\boldsymbol{N}}_{t-1};P)}\left\{\dfrac{I_{{\boldsymbol{a}}_{t}}(\boldsymbol{A}_{t})}{\pi_{a_t}(\boldsymbol{N}_t;P)}-1\right\}\{b_{\boldsymbol{a}_t}(\boldsymbol{L}_{t},\boldsymbol{N}_t;P) \nonumber \\
    &-& b_{\boldsymbol{a}_t}  (\boldsymbol{N}_{t};P)\}. \label{eq:proof4.1.2} %%%% Remplacement des lambda et pi(G, B) par lambda et pi(B)
\end{eqnarray}}
Defining 
	\[
	r_{_t}(\bar{\boldsymbol{A}}_t,\boldsymbol{L}_t,\bar{\boldsymbol{N}}_t;s_{\boldsymbol{a},t},P)\equiv \dfrac{I_{\bar{\boldsymbol{a}}_{t-1}}(\bar{\boldsymbol{A}}_{t-1})}{\lambda_{\bar{\boldsymbol{a}}_{t-1}}(\bar{\boldsymbol{N}}_{t-1};P)}\left\{\dfrac{I_{{\boldsymbol{a}}_{t}}(\boldsymbol{A}_{t})}{\pi_{a_t}(\boldsymbol{N}_t;P)}-1\right\}s_{\boldsymbol{a},t}(\boldsymbol{L}_{t},\boldsymbol{N}_{t})
	\]
	with
	\[
	s_{\boldsymbol{a},t}(\boldsymbol{L}_{t},\boldsymbol{N}_{t})\equiv b_{\boldsymbol{a}_t}(\boldsymbol{L}_{t},\boldsymbol{N}_t;P)-b_{\boldsymbol{a}_t}(\boldsymbol{N}_{t};P),
	\]
relation (\ref{eq:proof4.1.2}) becomes
\begin{equation}
   \psi_{P,\boldsymbol{a}}(\boldsymbol{N};\mathcal{G})	= \psi_{P,\boldsymbol{a}}(\boldsymbol{L},\boldsymbol{N};\mathcal{G})+\sum_{t=0}^{T}r_{_t}(\bar{\boldsymbol{A}}_t,\boldsymbol{L}_t,\bar{\boldsymbol{N}}_t;s_{\boldsymbol{a},t},P). \label{eq:proof4.1.3} 
\end{equation}
Using relation (\ref{eq:proof4.1.3})
	\begin{eqnarray}
	var_P\{\psi_{P,a}(\boldsymbol{N};\mathcal{G})\}
 &=& var_P\{\psi_{P,a}(\boldsymbol{L},\boldsymbol{N};\mathcal{G})\}+\sum_{t=0}^{T}var_P\{r_{_t}(\bar{\boldsymbol{A}}_t,\boldsymbol{L}_t,\bar{\boldsymbol{N}}_t;s_{\boldsymbol{a},t},P)\} \nonumber \\
 &+& 2\sum_{t=0}^{T} cov_P\left\{\psi_{P,a}(\boldsymbol{L},\boldsymbol{N};\mathcal{G}),r_{_t}(\bar{\boldsymbol{A}}_t,\boldsymbol{L}_t,\bar{\boldsymbol{N}}_t;s_{\boldsymbol{a},t},P)\right\} \nonumber \\
 &+& 2\sum_{0\leq i<t\leq T}cov_P\left\{r_{_t}(\bar{\boldsymbol{A}}_t,\boldsymbol{L}_t,\bar{\boldsymbol{N}}_t;s_{\boldsymbol{a},t},P),r_{i}(\bar{\boldsymbol{A}}_{i},\boldsymbol{L}_{i},\bar{\boldsymbol{N}}_{i};s_{\boldsymbol{a},i},P)\right\}. 
	\end{eqnarray}
For any $0\leq i<t \leq T$, we have
    
\resizebox{1\textwidth}{!}{$
\begin{aligned}
cov_P\left\{r_{_t}(\bar{\boldsymbol{A}}_t,\boldsymbol{L}_t,\bar{\boldsymbol{N}}_t;s_{\boldsymbol{a},t},P),r_{i}(\bar{\boldsymbol{A}}_{i},\boldsymbol{L}_{i},\bar{\boldsymbol{N}}_{i};s_{\boldsymbol{a},i},P)\right\}
&= \mathbb{E}_P\{r_{_t}(\bar{\boldsymbol{A}}_t,\boldsymbol{L}_t,\bar{\boldsymbol{N}}_t;s_{\boldsymbol{a},t},P)r_{i}(\bar{\boldsymbol{A}}_{i},\boldsymbol{L}_{i},\bar{\boldsymbol{N}}_{i};s_{\boldsymbol{a},i},P)\}  \\
&- \mathbb{E}_P\{r_{_t}(\bar{\boldsymbol{A}}_t,\boldsymbol{L}_t,\bar{\boldsymbol{N}}_t;s_{\boldsymbol{a},t},P)\}E_P\{r_{i}(\bar{\boldsymbol{A}}_{i},\boldsymbol{L}_{i},\bar{\boldsymbol{N}}_{i};s_{\boldsymbol{a},i},P)\}
\end{aligned} 
$
}

We now show that $\mathbb{E}_P\{r_{_t}(\bar{\boldsymbol{A}}_t,\boldsymbol{L}_t,\bar{\boldsymbol{N}}_t;s_{\boldsymbol{a},t},P)\} = 0$.

{\footnotesize
\begin{eqnarray} 
		\mathbb{E}_P\{r_{_t}(\bar{\boldsymbol{A}}_t,\boldsymbol{L}_t,\bar{\boldsymbol{N}}_t;s_{\boldsymbol{a},t},P)\}
		&=&\mathbb{E}_P[\mathbb{E}_P\{r_{_t}(\bar{\boldsymbol{A}}_t,\boldsymbol{L}_t,\bar{\boldsymbol{N}}_t;s_{\boldsymbol{a},t},P)\big |   \boldsymbol{L}_t,\bar{\boldsymbol{N}}_t,\bar{\boldsymbol{A}}_{t-1}=\bar{\boldsymbol{a}}_{t-1}\}] \nonumber \\ %%%% Total expectation
		&=& 
		\mathbb{E}_P\left( \mathbb{E}_P\left[\dfrac{I_{\bar{\boldsymbol{a}}_{t-1}}(\bar{\boldsymbol{A}}_{t-1})s_{\boldsymbol{a},t}(\boldsymbol{L}_{t},\boldsymbol{N}_{t})}{\lambda_{\bar{\boldsymbol{a}}_{t-1}}(\bar{\boldsymbol{N}}_{t-1};P)}\left\{\dfrac{I_{{\boldsymbol{a}}_{t}}(\boldsymbol{A}_{t})}{\pi_{a_t}(\boldsymbol{N}_t;P)}-1\right\}\Bigg |  \boldsymbol{L}_t,\bar{\boldsymbol{N}}_t,\bar{\boldsymbol{A}}_{t-1}=\bar{\boldsymbol{a}}_{t-1}\right]\right)  \nonumber \\ %%%% Definition + Juste le cas A_k = a_k n'est pas = 0
		&=& 
        \mathbb{E}_P\left(\dfrac{I_{\bar{\boldsymbol{a}}_{t-1}}(\bar{\boldsymbol{A}}_{t-1})s_{\boldsymbol{a},t}(\boldsymbol{L}_{t},\boldsymbol{N}_{t})}{\lambda_{\bar{\boldsymbol{a}}_{t-1}}(\bar{\boldsymbol{N}}_{t-1};P)}\left[\mathbb{E}_P\left\{\dfrac{I_{{\boldsymbol{a}}_{t}}(\boldsymbol{A}_{t})}{\pi_{a_t}(\boldsymbol{N}_t;P)}-1\right\} \Bigg |   \boldsymbol{L}_t,\bar{\boldsymbol{N}}_t,\bar{\boldsymbol{A}}_{t-1}=\bar{\boldsymbol{a}}_{t-1}\right]\right) \nonumber \\ %%%% Le terme sorti de l'espérance est une constante sachant la condition
		&=& 
		\mathbb{E}_P\left( \dfrac{I_{\bar{\boldsymbol{a}}_{t-1}}(\bar{\boldsymbol{A}}_{t-1})s_{\boldsymbol{a},t}(\boldsymbol{L}_{t},\boldsymbol{N}_{t})}{\lambda_{\bar{\boldsymbol{a}}_{t-1}}(\bar{\boldsymbol{N}}_{t-1};P)}\left[\mathbb{E}_P\left\{\dfrac{I_{{\boldsymbol{a}}_{t}}(\boldsymbol{A}_{t})}{\pi_{a_t}(\boldsymbol{N}_t;P)}\Bigg |   \boldsymbol{L}_t,\boldsymbol{N}_t\right\}-1\right]\right)  \label{eq10} \\ %%%% Sortie du -1 et Gk, Bk incluent Gbar, Bbar, Abar
  &=& \mathbb{E}_P\left( \dfrac{I_{\bar{\boldsymbol{a}}_{t-1}}(\bar{\boldsymbol{A}}_{t-1})s_{\boldsymbol{a},t}(\boldsymbol{L}_{t},\boldsymbol{N}_{t})}{\lambda_{\bar{\boldsymbol{a}}_{t-1}}(\bar{\boldsymbol{N}}_{t-1};P)}\left[\mathbb{E}_P\left\{\dfrac{I_{{\boldsymbol{a}}_{t}}(\boldsymbol{A}_{t})}{\pi_{a_t}(\boldsymbol{N}_t;P)}\Bigg |   \boldsymbol{N}_t\right\}-1\right]\right)  \label{eq1} \\ %%%% A_k indep G_k sachant B_k, ce qui permet de retirer G_k 
		&=&
		\mathbb{E}_P\left( \dfrac{I_{\bar{\boldsymbol{a}}_{t-1}}(\bar{\boldsymbol{A}}_{t-1})s_{\boldsymbol{a},t}(\boldsymbol{L}_{t},\boldsymbol{N}_{t})}{\lambda_{\bar{\boldsymbol{a}}_{t-1}}(\bar{\boldsymbol{N}}_{t-1};P)}\left[\dfrac{\mathbb{E}_P\left\{I_{\boldsymbol{a}_{t}}(\boldsymbol{A}_{t})\big |   \boldsymbol{N}_t\right\}}{\pi_{a_t}(\boldsymbol{N}_t;P)}-1\right]\right)   \nonumber \\
		&=&
		\mathbb{E}_P\left[ \dfrac{I_{\bar{\boldsymbol{a}}_{t-1}}(\bar{\boldsymbol{A}}_{t-1})s_{\boldsymbol{a},t}(\boldsymbol{L}_{t},\boldsymbol{N}_{t})}{\lambda_{\bar{\boldsymbol{a}}_{t-1}}(\bar{\boldsymbol{N}}_{t-1};P)}\left\{\dfrac{P(\boldsymbol{A}_{t}=a_t\big |\boldsymbol{N}_{t})}{\pi_{a_t}(\boldsymbol{N}_t;P)}-1\right\}\right]  \nonumber \\
		&=&
		\mathbb{E}_P\left[ \dfrac{I_{\bar{\boldsymbol{a}}_{t-1}}(\bar{\boldsymbol{A}}_{t-1})s_{\boldsymbol{a},t}(\boldsymbol{L}_{t},\boldsymbol{N}_{t})}{\lambda_{\bar{\boldsymbol{a}}_{t-1}}(\bar{\boldsymbol{N}}_{t-1};P)}\left\lbrace\dfrac{\pi_{a_t}(\boldsymbol{N}_t;P)}{\pi_{a_t}(\boldsymbol{N}_t;P)}-1\right\rbrace\right]  = 0 \nonumber 
	\end{eqnarray}}
where equality (\ref{eq10}) follows from the fact that 
 $\{\bar{\boldsymbol{A}}_{t-1},\bar{\boldsymbol{N}}_{t-1}\} \subset (\boldsymbol{L}_t,\boldsymbol{N}_t)$ and equality (\ref{eq1}) follows from (\ref{eq1.lemma1}). 
 
 We now turn our attention to showing that $\mathbb{E}_P\{r_{_t}(\bar{\boldsymbol{A}}_t,\boldsymbol{L}_t,\bar{\boldsymbol{N}}_t;P)r_{i}(\bar{\boldsymbol{A}}_{i},\boldsymbol{L}_{i},\bar{\boldsymbol{N}}_{i};P)\} = 0$. Together with the previous result, this will show that $cov_P\left\{r_{_t}(\bar{\boldsymbol{A}}_t,\boldsymbol{L}_t,\bar{\boldsymbol{N}}_t;s_{\boldsymbol{a},t},P),r_{i}(\bar{\boldsymbol{A}}_{i},\boldsymbol{L}_{i},\bar{\boldsymbol{N}}_{i};s_{\boldsymbol{a},i},P)\right\} = 0$. 
 \begin{footnotesize}
 \begin{align*}
 &\mathbb{E}_P\{r_{_t}(\bar{\boldsymbol{A}}_t,\boldsymbol{L}_t,\bar{\boldsymbol{N}}_t;P)r_{i}(\bar{\boldsymbol{A}}_{i},\boldsymbol{L}_{i},\bar{\boldsymbol{N}}_{i};P)\} \\ %%%% Par dfn
  & = \mathbb{E}_P\left[\dfrac{I_{\bar{\boldsymbol{a}}_{t-1}}(\bar{\boldsymbol{A}}_{t-1})I_{\bar{\boldsymbol{a}}_{i-1}}(\bar{\boldsymbol{A}}_{i-1})s_{\boldsymbol{a},t}(\boldsymbol{L}_{t},\boldsymbol{N}_{t})s_{\boldsymbol{a},i}(\boldsymbol{L}_{i},\boldsymbol{N}_{i})}{\lambda_{\bar{\boldsymbol{a}}_{t-1}}(\bar{\boldsymbol{N}}_{t-1};P)\lambda_{\bar{\boldsymbol{a}}_{i-1}}(\bar{\boldsymbol{N}}_{i-1};P)}\left\{\dfrac{I_{{\boldsymbol{a}}_{t}}(\boldsymbol{A}_{t})}{\pi_{a_t}(\boldsymbol{N}_t;P)}-1\right\}
		\left\{\dfrac{I_{{\boldsymbol{a}}_{i}}(\boldsymbol{A}_{i})}{\pi_{a_{i}}(\boldsymbol{N}_{i};P)}-1\right\}\right]  \\ %%%% Esp totales
  & = \mathbb{E}_P\left(\mathbb{E}_P\left[\dfrac{I_{\bar{\boldsymbol{a}}_{t-1}}(\bar{\boldsymbol{A}}_{t-1})I_{\bar{\boldsymbol{a}}_{i-1}}(\bar{\boldsymbol{A}}_{i-1})s_{\boldsymbol{a},t}(\boldsymbol{L}_{t},\boldsymbol{N}_{t})s_{\boldsymbol{a},i}(\boldsymbol{L}_{i},\boldsymbol{N}_{i})}{\lambda_{\bar{\boldsymbol{a}}_{t-1}}(\bar{\boldsymbol{N}}_{t-1};P)\lambda_{\bar{\boldsymbol{a}}_{i-1}}(\bar{\boldsymbol{N}}_{i-1};P)} \times \right. \right. \\ 
  & \qquad \qquad \left. \left.  \left\{\dfrac{I_{{\boldsymbol{a}}_{t}}(\boldsymbol{A}_{t})}{\pi_{a_t}(\boldsymbol{N}_t;P)}-1\right\}
		\left\{\dfrac{I_{{\boldsymbol{a}}_{i}}(\boldsymbol{A}_{i})}{\pi_{a_{i}}(\boldsymbol{N}_{i};P)}-1\right\}\right]\Bigg |   \boldsymbol{L}_t, \bar{\boldsymbol{N}}_{t}, \boldsymbol{L}_i, \bar{\boldsymbol{N}}_{i},\bar{\boldsymbol{A}}_{t-1}=\bar{\boldsymbol{a}}_{t-1},\bar{\boldsymbol{A}}_{i}=\bar{\boldsymbol{a}}_{i}\right) \\ %%%% On sort tout ce qui devient constant, on ne garde que le conditionnement Gk Bk qui inclut tout le reste (k > i), le reste suit le meme principe que la partie precedente
  & = \mathbb{E}_P\Biggl(\dfrac{I_{\bar{\boldsymbol{a}}_{t-1}}(\bar{\boldsymbol{A}}_{t-1})I_{\bar{\boldsymbol{a}}_{i-1}}(\bar{\boldsymbol{A}}_{i-1})s_{\boldsymbol{a},t}(\boldsymbol{L}_{t},\boldsymbol{N}_{t})s_{\boldsymbol{a},i}(\boldsymbol{L}_{i},\boldsymbol{N}_{i})}{\lambda_{\bar{\boldsymbol{a}}_{t-1}}(\bar{\boldsymbol{N}}_{t-1};P)\lambda_{\bar{\boldsymbol{a}}_{i-1}}(\bar{\boldsymbol{N}}_{i-1};P)}\left\{\dfrac{I_{{\boldsymbol{a}}_{i}}(\boldsymbol{A}_{i})}{\pi_{a_{i}}(\boldsymbol{N}_{i};P)}-1\right\}
  \times  
  \left[\mathbb{E}_P\Biggl\{\dfrac{I_{{\boldsymbol{a}}_{t}}(\boldsymbol{A}_{t})}{\pi_{a_t}(\boldsymbol{N}_t;P)}-1\Bigg |   \boldsymbol{L}_t, \boldsymbol{N}_{t}\Biggr\}\right] \Biggr)\\
  & = \mathbb{E}_P\Biggl[\dfrac{I_{\bar{\boldsymbol{a}}_{t-1}}(\bar{\boldsymbol{A}}_{t-1})I_{\bar{\boldsymbol{a}}_{i-1}}(\bar{\boldsymbol{A}}_{i-1})s_{\boldsymbol{a},t}(\boldsymbol{L}_{t},\boldsymbol{N}_{t})s_{\boldsymbol{a},i}(\boldsymbol{L}_{i},\boldsymbol{N}_{i})}{\lambda_{\bar{\boldsymbol{a}}_{t-1}}(\bar{\boldsymbol{N}}_{t-1};P)\lambda_{\bar{\boldsymbol{a}}_{i-1}}(\bar{\boldsymbol{N}}_{i-1};P)}\left\{\dfrac{I_{{\boldsymbol{a}}_{i}}(\boldsymbol{A}_{i})}{\pi_{a_{i}}(\boldsymbol{N}_{i};P)}-1\right\} \times 
  \mathbb{E}_P\Biggl\{\dfrac{I_{{\boldsymbol{a}}_{t}}(\boldsymbol{A}_{t})}{\pi_{a_t}(\boldsymbol{N}_t;P)}-1\Bigg |   \boldsymbol{N}_t\Biggr\}\Biggr]\\
  & = \mathbb{E}_P\Biggl(\dfrac{I_{\bar{\boldsymbol{a}}_{t-1}}(\bar{\boldsymbol{A}}_{t-1})I_{\bar{\boldsymbol{a}}_{i-1}}(\bar{\boldsymbol{A}}_{i-1})s_{\boldsymbol{a},t}(\boldsymbol{L}_{t},\boldsymbol{N}_{t})s_{\boldsymbol{a},i}(\boldsymbol{L}_{i},\boldsymbol{N}_{i})}{\lambda_{\bar{\boldsymbol{a}}_{t-1}}(\bar{\boldsymbol{N}}_{t-1};P)\lambda_{\bar{\boldsymbol{a}}_{i-1}}(\bar{\boldsymbol{N}}_{i-1};P)}\left\{\dfrac{I_{{\boldsymbol{a}}_{i}}(\boldsymbol{A}_{i})}{\pi_{a_{i}}(\boldsymbol{N}_{i};P)}-1\right\} 
  \times 
  \left[\mathbb{E}_P\Biggl\{\dfrac{I_{{\boldsymbol{a}}_{t}}(\boldsymbol{A}_{t})}{\pi_{a_t}(\boldsymbol{N}_t;P)}\Bigg |   \boldsymbol{N}_t\Biggr\}-1\right]\Biggr) =0.     
 \end{align*}
 \end{footnotesize}
We next show that $cov_P\left\{\psi_{P,a}(\boldsymbol{L},\boldsymbol{N};\mathcal{G}),r_t(\bar{\boldsymbol{A}}_t,\boldsymbol{L}_t,\bar{\boldsymbol{N}}_t;s_{\boldsymbol{a},t},P)\right\} = 0$: 
\begin{align}
& cov_P\left\{\psi_{P,a}(\boldsymbol{L},\boldsymbol{N};\mathcal{G}),r_t(\bar{\boldsymbol{A}}_t,\boldsymbol{L}_t,\bar{\boldsymbol{N}}_t;s_{\boldsymbol{a},t},P)\right\} \nonumber \\
&= cov_P\Biggl[\dfrac{I_{\bar{\boldsymbol{a}}_{T}}(\bar{\boldsymbol{A}}_{T})\{Y-b_{\boldsymbol{a}_T}(\boldsymbol{L}_T,\boldsymbol{N}_T;P)\}}{\lambda_{\bar{\boldsymbol{a}}_T}(\bar{\boldsymbol{L}}_T,\bar{\boldsymbol{N}}_T;P)},r_{_t}(\bar{\boldsymbol{A}}_t,\boldsymbol{L}_t,\bar{\boldsymbol{N}}_t;s_{\boldsymbol{a},t},P)\Biggr] \nonumber \\ 
& \qquad + cov_P\Biggl[\displaystyle{\sum_{i=0}^{T}}
 \dfrac{I_{\bar{\boldsymbol{a}}_{i-1}}(\bar{\boldsymbol{A}}_{i-1})\{b_{\boldsymbol{a}_i}(\boldsymbol{L}_i,\boldsymbol{N}_i;P)-b_{\boldsymbol{a}_{i-1}}(\boldsymbol{L}_{i-1},\boldsymbol{N}_{i-1};P)\}}{\lambda_{\bar{\boldsymbol{a}}_{i-1}}(\bar{\boldsymbol{L}}_{i-1},\bar{\boldsymbol{N}}_{i-1};P)}
 ,r_{_t}(\bar{\boldsymbol{A}}_t,\boldsymbol{L}_t,\bar{\boldsymbol{N}}_t;s_{\boldsymbol{a},t},P)\Biggr] \nonumber \\ %%%% Dfn de psi et division en deux cov
  & = cov_P\Biggl[\dfrac{I_{\bar{\boldsymbol{a}}_{T}}(\bar{\boldsymbol{A}}_{T})\{Y-b_{\boldsymbol{a}_T}(\boldsymbol{L}_T,\boldsymbol{N}_T;P)\}}{\lambda_{\bar{\boldsymbol{a}}_T}(\bar{\boldsymbol{L}}_T,\bar{\boldsymbol{N}}_T;P)},r_{_t}(\bar{\boldsymbol{A}}_t,\boldsymbol{L}_t,\bar{\boldsymbol{N}}_t;s_{\boldsymbol{a},t},P)\Biggr]\nonumber \\ 
  & +\displaystyle{\sum_{i=0}^{T}}cov_P\Biggl[
 \dfrac{I_{\bar{\boldsymbol{a}}_{i-1}}(\bar{\boldsymbol{A}}_{i-1})\{b_{\boldsymbol{a}_i}(\boldsymbol{L}_i,\boldsymbol{N}_i;P)-b_{\boldsymbol{a}_{i-1}}(\boldsymbol{L}_{i-1},\boldsymbol{N}_{i-1};P)\}}{\lambda_{\bar{\boldsymbol{a}}_{i-1}}(\bar{\boldsymbol{L}}_{i-1},\bar{\boldsymbol{N}}_{i-1};P)}
 ,r_{_t}(\bar{\boldsymbol{A}}_t,\boldsymbol{L}_t,\bar{\boldsymbol{N}}_t;s_{\boldsymbol{a},t},P)\Biggr]. \label{eq12} %%%% sortie de la somme de la 2e cov
\end{align}

We have
\begin{align*}
 &cov_P\Biggl[\dfrac{I_{\bar{\boldsymbol{a}}_{T}}(\bar{\boldsymbol{A}}_{T})\{Y-b_{\boldsymbol{a}_T}(\boldsymbol{L}_T,\boldsymbol{N}_T;P)\}}{\lambda_{\bar{\boldsymbol{a}}_T}(\bar{\boldsymbol{L}}_T,\bar{\boldsymbol{N}}_T;P)} ,r_{_t}(\bar{\boldsymbol{A}}_t,\boldsymbol{L}_t,\bar{\boldsymbol{N}}_t;s_{\boldsymbol{a},t},P)\Biggr] \\ 
 & = cov_P\Biggl[\dfrac{I_{\bar{\boldsymbol{a}}_{T}}(\bar{\boldsymbol{A}}_{T})\{Y-b_{\boldsymbol{a}_T}(\boldsymbol{L}_T,\boldsymbol{N}_T;P)\}}{\lambda_{\bar{\boldsymbol{a}}_T}(\bar{\boldsymbol{L}}_T,\bar{\boldsymbol{N}}_T;P)},\dfrac{I_{\bar{\boldsymbol{a}}_{t-1}}(\bar{\boldsymbol{A}}_{t-1})}{\lambda_{\bar{\boldsymbol{a}}_{t-1}}(\bar{\boldsymbol{N}}_{t-1};P)}\left\{\dfrac{I_{{\boldsymbol{a}}_{t}}(\boldsymbol{A}_{t})}{\pi_{a_t}(\boldsymbol{N}_t;P)}-1\right\}s_{\boldsymbol{a},t}(\boldsymbol{L}_{t},\boldsymbol{N}_{t})\Biggr]\\ %%%% Dfn de rk
 &= \mathbb{E}_P\Biggl[\dfrac{I_{\bar{\boldsymbol{a}}_{T}}(\bar{\boldsymbol{A}}_{T})s_{\boldsymbol{a},t}(\boldsymbol{L}_{t},\boldsymbol{N}_{t})\{Y-b_{\boldsymbol{a}_T}(\boldsymbol{L}_T,\boldsymbol{N}_T;P)\}}{\lambda_{\bar{\boldsymbol{a}}_T}(\bar{\boldsymbol{L}}_T,\bar{\boldsymbol{N}}_T;P)}\dfrac{I_{\bar{\boldsymbol{a}}_{t-1}}(\bar{\boldsymbol{A}}_{t-1})}{\lambda_{\bar{\boldsymbol{a}}_{t-1}}(\bar{\boldsymbol{N}}_{t-1};P)}\left\{\dfrac{I_{{\boldsymbol{a}}_{t}}(\boldsymbol{A}_{t})}{\pi_{a_t}(\boldsymbol{N}_t;P)}-1\right\}\Biggr]\\  %%%% Cov = E[produit] - produit des E, mais E[rk] = 0
&=\mathbb{E}_P\Biggl(\mathbb{E}_P\Biggl[\dfrac{I_{\bar{\boldsymbol{a}}_{T}}(\bar{\boldsymbol{A}}_{T})I_{\bar{\boldsymbol{a}}_{t-1}}(\bar{\boldsymbol{A}}_{t-1})\{Y-b_{\boldsymbol{a}_T}(\boldsymbol{L}_T,\boldsymbol{N}_T;P)s_{\boldsymbol{a},t}(\boldsymbol{L}_{t},\boldsymbol{N}_{t})\}}{\lambda_{\bar{\boldsymbol{a}}_T}(\bar{\boldsymbol{L}}_T,\bar{\boldsymbol{N}}_T;P)\lambda_{\bar{\boldsymbol{a}}_{t-1}}(\bar{\boldsymbol{N}}_{t-1};P)} \times  \\
& \qquad \qquad  \left\{\dfrac{I_{{\boldsymbol{a}}_{t}}(\boldsymbol{A}_{t})}{\pi_{a_t}(\boldsymbol{N}_t;P)}-1\right\}\Bigg |   \boldsymbol{\bar{L}}_{T},\boldsymbol{\bar{N}}_{T}, \boldsymbol{\bar{N}}_{t}, \boldsymbol{{L}}_{t},\boldsymbol{\bar{A}}_{T}=\boldsymbol{\bar{a}}_{T},\boldsymbol{\bar{A}}_{t}=\boldsymbol{\bar{a}}_{t} \Biggr]\Biggr)\\ %%%% Meme principe que d'habitude
&=\mathbb{E}_P\Biggl(\dfrac{I_{\bar{\boldsymbol{a}}_{T}}(\bar{\boldsymbol{A}}_{T})I_{\bar{\boldsymbol{a}}_{t-1}}(\bar{\boldsymbol{A}}_{t-1})s_{\boldsymbol{a},t}(\boldsymbol{L}_{t},\boldsymbol{N}_{t})}{\lambda_{\bar{\boldsymbol{a}}_T}(\bar{\boldsymbol{L}}_T,\bar{\boldsymbol{N}}_T;P)\lambda_{\bar{\boldsymbol{a}}_{t-1}}(\bar{\boldsymbol{N}}_{t-1};P)}\left\{\dfrac{I_{{\boldsymbol{a}}_{t}}(\boldsymbol{A}_{t})}{\pi_{a_t}(\boldsymbol{N}_t;P)}-1\right\} \times \\
& \qquad \qquad \mathbb{E}_P[\{Y-b_{\boldsymbol{a}_T}(\boldsymbol{L}_T,\boldsymbol{N}_T;P)\}\big |  \boldsymbol{\bar{L}}_{T},\boldsymbol{\bar{N}}_{T}, \boldsymbol{\bar{N}}_{t}, \boldsymbol{{L}}_{t},\boldsymbol{\bar{A}}_{T}=\boldsymbol{\bar{a}}_{T},\boldsymbol{\bar{A}}_{t}=\boldsymbol{\bar{a}}_{t}]\Biggr) =0,
\end{align*}
because
\begin{align*}
   & \mathbb{E}_P[\{Y-b_{\boldsymbol{a}_T}(\boldsymbol{L}_T,\boldsymbol{N}_T;P)\}\big |   \boldsymbol{\bar{L}}_{T},\boldsymbol{\bar{N}}_{T}, \boldsymbol{\bar{N}}_{t}, \boldsymbol{{L}}_{t},\boldsymbol{\bar{A}}_{T}=\boldsymbol{\bar{a}}_{T},\boldsymbol{\bar{A}}_{t}=\boldsymbol{\bar{a}}_{t}]\\
   & = \mathbb{E}_P(Y\big |   \boldsymbol{{L}}_{T},\boldsymbol{{N}}_{T},\boldsymbol{\bar{L}}_{T-1},\boldsymbol{\bar{N}}_{T-1}, \boldsymbol{\bar{N}}_{t}, \boldsymbol{{L}}_{t},\boldsymbol{\bar{A}}_{T}=\boldsymbol{\bar{a}}_{T},\boldsymbol{\bar{A}}_{t}=\boldsymbol{\bar{a}}_{t})\\
   & \qquad \qquad - \mathbb{E}_P\{b_{\boldsymbol{a}_T}(\boldsymbol{L}_T,\boldsymbol{N}_T;P)\big |   \boldsymbol{{L}}_{T},\boldsymbol{{N}}_{T},\boldsymbol{\bar{L}}_{T-1},\boldsymbol{\bar{N}}_{T-1}, \boldsymbol{\bar{N}}_{t}, \boldsymbol{{L}}_{t},\boldsymbol{\bar{A}}_{T}=\boldsymbol{\bar{a}}_{T},\boldsymbol{\bar{A}}_{t}=\boldsymbol{\bar{a}}_{t}\} \\
   & = \scalebox{0.9}{$
   \mathbb{E}_P(Y\big |  \boldsymbol{{L}}_{T},\boldsymbol{{N}}_{T},\boldsymbol{a}_T)-\mathbb{E}_P\{b_{\boldsymbol{a}_T}(\boldsymbol{L}_T,\boldsymbol{N}_T;P)\big |  \boldsymbol{{L}}_{T},\boldsymbol{{N}}_{T},\boldsymbol{a}_T\}=b_{\boldsymbol{a}_T}(\boldsymbol{L}_T,\boldsymbol{N}_T;P)-b_{\boldsymbol{a}_T}(\boldsymbol{L}_T,\boldsymbol{N}_T;P)=0.$}
\end{align*}
It remains to show that 
$$cov_P\Biggl[
 \dfrac{I_{\bar{\boldsymbol{a}}_{i-1}}(\bar{\boldsymbol{A}}_{i-1})\{b_{\boldsymbol{a}_i}(\boldsymbol{L}_i,\boldsymbol{N}_i;P)-b_{\boldsymbol{a}_{i-1}}(\boldsymbol{L}_{i-1},\boldsymbol{N}_{i-1};P)\}}{\lambda_{\bar{\boldsymbol{a}}_{i-1}}(\bar{\boldsymbol{L}}_{i-1},\bar{\boldsymbol{N}}_{i-1};P)}
 ,r_{_t}(\bar{\boldsymbol{A}}_t,\boldsymbol{L}_t,\bar{\boldsymbol{N}}_t;s_{\boldsymbol{a},t},P)\Biggr] = 0.$$
We consider separately the case where $t\geq i$ and the case where $t < i$. For $t\geq i$, we have 
\begin{align*}
& cov_P\Biggl[
 \dfrac{I_{\bar{\boldsymbol{a}}_{i-1}}(\bar{\boldsymbol{A}}_{i-1})\{b_{\boldsymbol{a}_i}(\boldsymbol{L}_i,\boldsymbol{N}_i;P)-b_{\boldsymbol{a}_{i-1}}(\boldsymbol{L}_{i-1},\boldsymbol{N}_{i-1};P)\}}{\lambda_{\bar{\boldsymbol{a}}_{i-1}}(\bar{\boldsymbol{L}}_{i-1},\bar{\boldsymbol{N}}_{i-1};P)}
 ,r_{_t}(\bar{\boldsymbol{A}}_t,\boldsymbol{L}_t,\bar{\boldsymbol{N}}_t;s_{\boldsymbol{a},t},P)\Biggr]\\
 &= \scalebox{0.9}{$
 {cov_P\Biggl[
 \dfrac{I_{\bar{\boldsymbol{a}}_{i-1}}(\bar{\boldsymbol{A}}_{i-1})\{b_{\boldsymbol{a}_i}(\boldsymbol{L}_i,\boldsymbol{N}_i;P)-b_{\boldsymbol{a}_{i-1}}(\boldsymbol{L}_{i-1},\boldsymbol{N}_{i-1};P)\}}{\lambda_{\bar{\boldsymbol{a}}_{i-1}}(\bar{\boldsymbol{L}}_{i-1},\bar{\boldsymbol{N}}_{i-1};P)} , 
 \dfrac{I_{\bar{\boldsymbol{a}}_{t-1}}(\bar{\boldsymbol{A}}_{t-1})s_{\boldsymbol{a},t}(\boldsymbol{L}_{t},\boldsymbol{N}_{t})}{\lambda_{\bar{\boldsymbol{a}}_{t-1}}(\bar{\boldsymbol{N}}_{t-1};P)}\left\{\dfrac{I_{{\boldsymbol{a}}_{t}}(\boldsymbol{A}_{t})}{\pi_{a_{t}}(\boldsymbol{N}_{t};P)}-1\right\}\Biggr]}$}\\ %%%% Dfn
 &= \scalebox{0.9}{$
 \mathbb{E}_P\Biggl[
\dfrac{I_{\bar{\boldsymbol{a}}_{i-1}}(\bar{\boldsymbol{A}}_{i-1})I_{\bar{\boldsymbol{a}}_{t-1}}(\bar{\boldsymbol{A}}_{t-1})s_{\boldsymbol{a},t}(\boldsymbol{L}_{t},\boldsymbol{N}_{t})\{b_{\boldsymbol{a}_i}(\boldsymbol{L}_i,\boldsymbol{N}_i;P)-b_{\boldsymbol{a}_{i-1}}(\boldsymbol{L}_{i-1},\boldsymbol{N}_{i-1};P)\}}{\lambda_{\bar{\boldsymbol{a}}_{i-1}}(\bar{\boldsymbol{L}}_{i-1},\bar{\boldsymbol{N}}_{i-1};P)\lambda_{\bar{\boldsymbol{a}}_{t-1}}(\bar{\boldsymbol{N}}_{t-1};P)}
\left\{\dfrac{I_{{\boldsymbol{a}}_{t}}(\boldsymbol{A}_{t})}{\pi_{a_{t}}(\boldsymbol{N}_{t};P)}-1\right\}\Biggr]$}\\ %%% E[produit] - produit des E, mais E[rk] = 0
&=\mathbb{E}_P\Biggl(\mathbb{E}_P\Biggl[
\dfrac{I_{\bar{\boldsymbol{a}}_{i-1}}(\bar{\boldsymbol{A}}_{i-1})I_{\bar{\boldsymbol{a}}_{t-1}}(\bar{\boldsymbol{A}}_{t-1})s_{\boldsymbol{a},t}(\boldsymbol{L}_{t},\boldsymbol{N}_{t})\{b_{\boldsymbol{a}_i}(\boldsymbol{L}_i,\boldsymbol{N}_i;P)-b_{\boldsymbol{a}_{i-1}}(\boldsymbol{L}_{i-1},\boldsymbol{N}_{i-1};P)\}}{\lambda_{\bar{\boldsymbol{a}}_{i-1}}(\bar{\boldsymbol{L}}_{i-1},\bar{\boldsymbol{N}}_{i-1};P)\lambda_{\bar{\boldsymbol{a}}_{t-1}}(\bar{\boldsymbol{N}}_{t-1};P)} \times  \\
& \qquad \qquad  \left\{\dfrac{I_{{\boldsymbol{a}}_{t}}(\boldsymbol{A}_{t})}{\pi_{a_{t}}(\boldsymbol{N}_{t};P)}-1\right\}\Bigg |  \bar{\boldsymbol{L}}_{i},\bar{\boldsymbol{N}}_{i},\bar{\boldsymbol{N}}_{t},\boldsymbol{L}_{t},\bar{\boldsymbol{A}}_{i-1}=\bar{\boldsymbol{a}}_{i-1},\bar{\boldsymbol{A}}_{t-1}=\bar{\boldsymbol{a}}_{t-1}\Biggr]\Biggr) \\ %%%% Esperances totales
& = \mathbb{E}_P\Biggl(
\dfrac{I_{\bar{\boldsymbol{a}}_{i-1}}(\bar{\boldsymbol{A}}_{i-1})I_{\bar{\boldsymbol{a}}_{t-1}}(\bar{\boldsymbol{A}}_{t-1})s_{\boldsymbol{a},t}(\boldsymbol{L}_{t},\boldsymbol{N}_{t})\{b_{\boldsymbol{a}_i}(\boldsymbol{L}_i,\boldsymbol{N}_i;P)-b_{\boldsymbol{a}_{i-1}}(\boldsymbol{L}_{i-1},\boldsymbol{N}_{i-1};P)\}}{\lambda_{\bar{\boldsymbol{a}}_{i-1}}(\bar{\boldsymbol{L}}_{i-1},\bar{\boldsymbol{N}}_{i-1};P)\lambda_{\bar{\boldsymbol{a}}_{t-1}}(\bar{\boldsymbol{N}}_{t-1};P)} \times \\
& \qquad \qquad  \mathbb{E}_P\Biggl[\left\lbrace\dfrac{I_{{\boldsymbol{a}}_{t}}(\boldsymbol{A}_{t})}{\pi_{a_{t}}(\boldsymbol{N}_{t};P)}-1\right\rbrace \Bigg | \bar{\boldsymbol{L}}_{i},\bar{\boldsymbol{N}}_{i},\bar{\boldsymbol{N}}_{t},\boldsymbol{L}_{t},\bar{\boldsymbol{A}}_{i-1}=\bar{\boldsymbol{a}}_{i-1},\bar{\boldsymbol{A}}_{t-1}=\bar{\boldsymbol{a}}_{t-1}\Biggr]\Biggr)\\ %%% On sort les constantes
& = \mathbb{E}_P\Biggl(
\dfrac{I_{\bar{\boldsymbol{a}}_{i-1}}(\bar{\boldsymbol{A}}_{i-1})I_{\bar{\boldsymbol{a}}_{t-1}}(\bar{\boldsymbol{A}}_{t-1})s_{\boldsymbol{a},t}(\boldsymbol{L}_{t},\boldsymbol{N}_{t})\{b_{\boldsymbol{a}_i}(\boldsymbol{L}_i,\boldsymbol{N}_i;P)-b_{\boldsymbol{a}_{i-1}}(\boldsymbol{L}_{i-1},\boldsymbol{N}_{i-1};P)\}}{\lambda_{\bar{\boldsymbol{a}}_{i-1}}(\bar{\boldsymbol{L}}_{i-1},\bar{\boldsymbol{N}}_{i-1};P)\lambda_{\bar{\boldsymbol{a}}_{t-1}}(\bar{\boldsymbol{N}}_{t-1};P)} \times \\
& \qquad \qquad \mathbb{E}_P\Biggl[\left\lbrace\dfrac{I_{{\boldsymbol{a}}_{t}}(\boldsymbol{A}_{t})}{\pi_{a_{t}}(\boldsymbol{N}_{t};P)}-1\right\rbrace \Bigg | \boldsymbol{N}_{t}\Biggr]\Biggr)=0. %%%% Gk et Bk = tout le reste + independant de Gk sachant Bk
\end{align*}
Similarly, for $t < i$,
\begin{align*}
 &cov_P\Biggl[
 \dfrac{I_{\bar{\boldsymbol{a}}_{i-1}}(\bar{\boldsymbol{A}}_{i-1})\{b_{\boldsymbol{a}_i}(\boldsymbol{L}_i,\boldsymbol{N}_i;P)-b_{\boldsymbol{a}_{i-1}}(\boldsymbol{L}_{i-1},\boldsymbol{N}_{i-1};P)\}}{\lambda_{\bar{\boldsymbol{a}}_{i-1}}(\bar{\boldsymbol{L}}_{i-1},\bar{\boldsymbol{N}}_{i-1};P)}
 ,r_{_t}(\bar{\boldsymbol{A}}_t,\boldsymbol{L}_t,\bar{\boldsymbol{N}}_t;s_{\boldsymbol{a},t},P)\Biggr]\\ 
 &= \mathbb{E}_P\Biggl(
\dfrac{I_{\bar{\boldsymbol{a}}_{i-1}}(\bar{\boldsymbol{A}}_{i-1})I_{\bar{\boldsymbol{a}}_{t-1}}(\bar{\boldsymbol{A}}_{t-1})s_{\boldsymbol{a},t}(\boldsymbol{L}_{t},\boldsymbol{N}_{t})}{\lambda_{\bar{\boldsymbol{a}}_{i-1}}(\bar{\boldsymbol{L}}_{i-1},\bar{\boldsymbol{N}}_{i-1};P)\lambda_{\bar{\boldsymbol{a}}_{t-1}}(\bar{\boldsymbol{N}}_{t-1};P)}
\left\{\dfrac{I_{{\boldsymbol{a}}_{t}}(\boldsymbol{A}_{t})}{\pi_{a_{t}}(\boldsymbol{N}_{t};P)}-1\right\} \times \\
& \mathbb{E}_P[\{b_{\boldsymbol{a}_i}(\boldsymbol{L}_i,\boldsymbol{N}_i;P)-b_{\boldsymbol{a}_{i-1}}(\boldsymbol{L}_{i-1},\boldsymbol{N}_{i-1};P)\}\mid\bar{\boldsymbol{L}}_{i-1},\bar{\boldsymbol{N}}_{i-1},\bar{\boldsymbol{N}}_{t},\boldsymbol{L}_{t},\bar{\boldsymbol{A}}_{i-1}=\bar{\boldsymbol{a}}_{i-1},\bar{\boldsymbol{A}}_{t}=\bar{\boldsymbol{a}}_{t}]\Biggr)
  \\ %%%% Dfn + (cov = E) + esperances totales + sortie des constantes
& = \mathbb{E}_P\Biggl(
\dfrac{I_{\bar{\boldsymbol{a}}_{i-1}}(\bar{\boldsymbol{A}}_{i-1})I_{\bar{\boldsymbol{a}}_{t-1}}(\bar{\boldsymbol{A}}_{t-1})s_{\boldsymbol{a},t}(\boldsymbol{L}_{t},\boldsymbol{N}_{t})}{\lambda_{\bar{\boldsymbol{a}}_{i-1}}(\bar{\boldsymbol{L}}_{i-1},\bar{\boldsymbol{N}}_{i-1};P)\lambda_{\bar{\boldsymbol{a}}_{t-1}}(\bar{\boldsymbol{N}}_{t-1};P)}
\left\{\dfrac{I_{{\boldsymbol{a}}_{t}}(\boldsymbol{A}_{t})}{\pi_{a_{t}}(\boldsymbol{N}_{t};P)}-1\right\} \times \\
& \qquad \qquad \mathbb{E}_P[\{b_{\boldsymbol{a}_i}(\boldsymbol{L}_i,\boldsymbol{N}_i;P)-b_{\boldsymbol{a}_{i-1}}(\boldsymbol{L}_{i-1},\boldsymbol{N}_{i-1};P)\}\mid\boldsymbol{a}_{i-1},\boldsymbol{L}_{i-1},\boldsymbol{N}_{i-1}]\Biggr) \\ %%%% On ne garde que le conditionnement p/r a i qui inclut le reste
& = \mathbb{E}_P\Biggl(
\dfrac{I_{\bar{\boldsymbol{a}}_{i-1}}(\bar{\boldsymbol{A}}_{i-1})I_{\bar{\boldsymbol{a}}_{t-1}}(\bar{\boldsymbol{A}}_{t-1})s_{\boldsymbol{a},t}(\boldsymbol{L}_{t},\boldsymbol{N}_{t})}{\lambda_{\bar{\boldsymbol{a}}_{i-1}}(\bar{\boldsymbol{L}}_{i-1},\bar{\boldsymbol{N}}_{i-1};P)\lambda_{\bar{\boldsymbol{a}}_{t-1}}(\bar{\boldsymbol{N}}_{t-1};P)}
\left\{\dfrac{I_{{\boldsymbol{a}}_{t}}(\boldsymbol{A}_{t})}{\pi_{a_{t}}(\boldsymbol{N}_{t};P)}-1\right\} \times \\
& \qquad \qquad [\mathbb{E}_P\{b_{\boldsymbol{a}_i}(\boldsymbol{L}_i,\boldsymbol{N}_i;P)\mid a_{i-1},\boldsymbol{L}_{i-1},\boldsymbol{N}_{i-1}\}-b_{\boldsymbol{a}_{i-1}}(\boldsymbol{L}_{i-1},\boldsymbol{N}_{i-1};P)]\Biggr) \\
& = \mathbb{E}_P\Biggl[
\dfrac{I_{\bar{\boldsymbol{a}}_{i-1}}(\bar{\boldsymbol{A}}_{i-1})I_{\bar{\boldsymbol{a}}_{t-1}}(\bar{\boldsymbol{A}}_{t-1})s_{\boldsymbol{a},t}(\boldsymbol{L}_{t},\boldsymbol{N}_{t})}{\lambda_{\bar{\boldsymbol{a}}_{i-1}}(\bar{\boldsymbol{L}}_{i-1},\bar{\boldsymbol{N}}_{i-1};P)\lambda_{\bar{\boldsymbol{a}}_{t-1}}(\bar{\boldsymbol{N}}_{t-1};P)}
\left\{\dfrac{I_{{\boldsymbol{a}}_{t}}(\boldsymbol{A}_{t})}{\pi_{a_{t}}(\boldsymbol{N}_{t};P)}-1\right\} \times \\
& \qquad \qquad \{b_{\boldsymbol{a}_{i-1}}(\boldsymbol{L}_{i-1},\boldsymbol{N}_{i-1};P)-b_{\boldsymbol{a}_{i-1}}(\boldsymbol{L}_{i-1},\boldsymbol{N}_{i-1};P)\}\Biggr]=0.
\end{align*}

We conclude that $var_P\{\psi_{P,\boldsymbol{a}}(\boldsymbol{N};\mathcal{G})\}=var_P\{\psi_{P,a}(\boldsymbol{L},\boldsymbol{N};\mathcal{G})\}+\sum_{t=0}^{T}var_P\{r_{_t}(\bar{\boldsymbol{A}}_t,\boldsymbol{L}_t,\bar{\boldsymbol{N}}_t;s_{\boldsymbol{a},t},P)\}$, which shows that $\sigma^{2}_{\boldsymbol{a},\boldsymbol{N}}(P)-\sigma^{2}_{\boldsymbol{a},\boldsymbol{L},\boldsymbol{N}}(P)\geq 0$. 

It only remains to show that $\sigma^{2}_{\Delta,\boldsymbol{N}}(P)-\sigma^{2}_{\Delta,\boldsymbol{L},\boldsymbol{N}}(P)\geq 0$. Let $\boldsymbol{c}\equiv (c_a)_{a \in \mathcal{\boldsymbol{A}}}$ and $\psi_{P, \Delta}(\boldsymbol{Z};\mathcal{G})\equiv \{\psi_{P,a}(\boldsymbol{Z};\mathcal{G})\}_{a \in \mathcal{\boldsymbol{A}}}$ for all $\boldsymbol{Z}$. 
	We can then write
	\begin{eqnarray} 
		\sum_{a \in \mathcal{\boldsymbol{A}}} c_a\psi_{P,a}(\boldsymbol{Z};\mathcal{G})
		&=&\boldsymbol{c}^{T}\psi_{P,\Delta}(\boldsymbol{Z};\mathcal{G}). \nonumber 
	\end{eqnarray}
	Following steps similar to previously, a result analogous to relation (\ref{eq:proof4.1.2}) can be obtained: 
	\begin{eqnarray*} 
		\boldsymbol{c}^{T}\psi_{P,\Delta}(\boldsymbol{N};\mathcal{G})
		&=& \boldsymbol{c}^{T}\psi_{P,\Delta}(\boldsymbol{L},\boldsymbol{N};\mathcal{G})+\sum_{t=0}^{T} h_t(\boldsymbol{\bar{A}}_t,\boldsymbol{L}_t,\boldsymbol{\bar{N}}_t,P)
	\end{eqnarray*}
	where
	
	$$	h_t(\boldsymbol{\bar{A}}_t,\boldsymbol{L}_t,\boldsymbol{\bar{N}}_t,P)\equiv\displaystyle{\sum_{\boldsymbol{a} \in \mathcal{A}}}c_{\boldsymbol{a}}\dfrac{I_{\bar{\boldsymbol{a}}_{t-1}}(\bar{\boldsymbol{A}}_{t-1})}{\lambda_{\bar{\boldsymbol{a}}_{t-1}}(\bar{\boldsymbol{N}}_{t-1};P)}\left\{\dfrac{I_{{\boldsymbol{a}}_{t}}(\boldsymbol{A}_{t})}{\pi_{a_t}(\boldsymbol{N}_t;P)}-1\right\}\{b_{\boldsymbol{a}_t}(\boldsymbol{L}_{t},\boldsymbol{N}_t;P)-b_{\boldsymbol{a}_t}(\boldsymbol{N}_{t};P)\}.
	$$
	Also following similar steps to previously, it can be shown that for any $0\leq i < t\leq T$, $h_i(\boldsymbol{\bar{A}}_i,\boldsymbol{L}_i,\boldsymbol{\bar{N}}_i,P)$ is independent of $h_t(\boldsymbol{\bar{A}}_t,\boldsymbol{L}_t,\boldsymbol{\bar{N}}_t,P)$ under $P$, and that $\boldsymbol{c}^{T}\psi_{P,\Delta}(\boldsymbol{L},\boldsymbol{N};\mathcal{G})$ is independent of $h_i(\boldsymbol{\bar{A}}_i,\boldsymbol{L}_i,\boldsymbol{\bar{N}}_i,P)$ under $P$. We thus obtain
	\begin{eqnarray} 
		\sigma_{\Delta,\boldsymbol{N}}^{2}(P)
		&=& var_P\{\boldsymbol{c}^{T}\psi_{P,\Delta}(\boldsymbol{N};\mathcal{G})\}\nonumber \\
		&=& var_P\{\boldsymbol{c}^{T}\psi_{P,\Delta}(\boldsymbol{L},\boldsymbol{N};\mathcal{G})\}+\sum_{t=0}^{T} var_P\{h_t(\boldsymbol{\bar{A}}_t,\boldsymbol{L}_t,\boldsymbol{\bar{N}}_t,P)\} \nonumber \\
		&=& \sigma_{\Delta,\boldsymbol{L},\boldsymbol{N}}^{2}(P)+\sum_{t=0}^{T} var_P\{h_t(\boldsymbol{\bar{A}}_t,\boldsymbol{L}_t,\boldsymbol{\bar{N}}_t,P)\} \nonumber \\
		\Leftrightarrow \sigma_{\Delta,\boldsymbol{N}}^{2}(P)-	\sigma_{\Delta,\boldsymbol{L},\boldsymbol{N}}^{2}(P)
		&=& \sum_{t=0}^{T} var_P\{h_t(\boldsymbol{\bar{A}}_t,\boldsymbol{L}_t,\boldsymbol{\bar{N}}_t,P)\} \geq 0. \nonumber
	\end{eqnarray} 
Since $A_t\indep \boldsymbol{L}_t\mid \boldsymbol{N}_t$,
it follows that for every subset \(\boldsymbol{L}'_{t} \subseteq \boldsymbol{L}_t\), one also has 
$A_t\indep \boldsymbol{L}'_{t} \mid \boldsymbol{N}_t$,
that is, conditional independence is stable under taking subsets.  
Applying exactly the same reasoning with \(\boldsymbol{L}'_{t}\) in place of \(\boldsymbol{L}_t\), we obtain $
\sigma^2_{\boldsymbol{a},\boldsymbol{N}}(P) \;\;\ge\;\; 
\sigma^2_{\boldsymbol{a},\boldsymbol{L'},\boldsymbol{N}}(P)\;\;\ge\;\;\sigma^2_{\boldsymbol{a},\boldsymbol{L},\boldsymbol{N}}(P)
$ and $\sigma^2_{\Delta,\boldsymbol{N}}(P) \;\;\ge\;\; 
\sigma^2_{\Delta,\boldsymbol{L'},\boldsymbol{N}}(P)\;\;\ge\;\; 
\sigma^2_{\Delta,\boldsymbol{L},\boldsymbol{N}}(P).$

This concludes the proof of Lemma \ref{lemma:1}.
\section*{Proof of Lemma 4.2} 
Set $Q_{T}({a}_{T+1}, \boldsymbol{U}_{T+1})=Y$. Note that $Q_{t-1}=Q_{t-1}({a}_{t}, \boldsymbol{U}_{t}) \equiv \mathbb{E}_P(Q_{t}({a}_{t+1}, \boldsymbol{U}_{t+1})\mid {A}_{t} = {a}_{t}, \boldsymbol{U}_{t}), t=T,\dots,1$.

Because $Q_{t-1}\indep_{\mathcal{G}_{t-1}}\boldsymbol{M}_{t-1} \ \mid  \ \boldsymbol{U}_{t-1},{A}_{t-1}$ for $  t=T,...,1$, we have by construction \begin{equation}
\boldsymbol{U}_{t}\backslash \boldsymbol{\bar{A}}_{t-1}\indep_{\mathcal{G}_{t-1}} \boldsymbol{M}_{t-1} \ \mid  \ \boldsymbol{U}_{t-1}, {A}_{t-1}\ \text{for} \  t=T,...,1.\label{eq18}
\end{equation}

We first show by reverse induction that for all $t \in \{0,...,T\}$
	\begin{equation}
			b_{{a}_t}(\boldsymbol{U}_t,\boldsymbol{M}_t;P) = b_{{a}_t}(\boldsymbol{U}_t;P). \label{eq:lemma2.3}
	\end{equation}
Note that
	\begin{eqnarray} 
b_{{a}_T}(\boldsymbol{U}_T,\boldsymbol{M}_T;P)
	&\equiv& \mathbb{E}_P(Y\mid \boldsymbol{U}_T,\boldsymbol{M}_T,{A}_T={a}_T) \nonumber \\
	&=& \mathbb{E}_P(Y\mid \boldsymbol{U}_T,{A}_T={a}_T)	 \nonumber \\ %%%% hypothèse 1 du lemme
	&=& 	b_{{a}_T}(\boldsymbol{U}_T;P)\nonumber %%%% dfn
\end{eqnarray}
where the second equality follows by (\ref{eq:lemma2.1}). This shows that (\ref{eq:lemma2.3}) holds for $t = T$. Next, assume that (\ref{eq:lemma2.3}) holds for $t \in \{t^{*}+1,...,T\}$ for some $t^{*}\geq0$. We now show that this implies that (\ref{eq:lemma2.3}) also holds for $t=t^{*}$. We have
	\begin{eqnarray} 
		b_{{a}_{t^{*}}}(\boldsymbol{U}_{t^{*}},\boldsymbol{M}_{t^{*}};P)
		&\equiv& \mathbb{E}_P\{b_{{a}_{t^{*}+1}}(\boldsymbol{U}_{t^{*}+1},\boldsymbol{M}_{t^{*}+1};P)\mid \boldsymbol{M}_{t^{*}},\boldsymbol{U}_{t^{*}},{A}_{t^{*}}={a}_{t^{*}}\} \nonumber \\ %%%% dfn
		&=& 	\mathbb{E}_P\{b_{{a}_{t^{*}+1}}(\boldsymbol{U}_{t^{*}+1};P)\mid \boldsymbol{M}_{t^{*}},\boldsymbol{U}_{t^{*}},{A}_{t^{*}}={a}_{t^{*}}\} \nonumber \\ %%%% par hypothese
		&=& 	\mathbb{E}_P\{b_{{a}_{t^{*}+1}}(\boldsymbol{U}_{t^{*}+1};P)\mid \boldsymbol{U}_{t^{*}},{A}_{t^{*}}={a}_{t^{*}}\} \equiv b_{{a}_t^{*}}(\boldsymbol{U}_t^{*};P),\nonumber %%%% hypothese 2 du lemme
	\end{eqnarray}
	where the second equality is obtained by the inductive hypothesis and the third by (\ref{eq18}) applied to  $t=t^{*}+1$. We thus conclude that (\ref{eq:lemma2.3}) holds for all $t \in \{0,...,T\}$. Consequently, $\mathbb{E}_T(Y^{\boldsymbol{a}})=\mathbb{E}_P\{b_{{a}_0}(\boldsymbol{U}_0,\boldsymbol{M}_0;P)\}= \mathbb{E}_P\{b_{{a}_0}(\boldsymbol{U}_0;P)\}$, where the first equality follows from the assumption that $(\boldsymbol{U}, \boldsymbol{M})$ is a sufficient time-dependent adjustment set and the second follows from (\ref{eq:lemma2.3}) applied to $t = 0$. Thus, $\boldsymbol{U} = (\boldsymbol{U}_0,\boldsymbol{U}_1,\cdots,\boldsymbol{U}_T)$ is a sufficient time-dependent adjustment set.
 
We now turn our attention to showing $\sigma^{2}_{\boldsymbol{a},\boldsymbol{U},\boldsymbol{M}}(P)-\sigma^{2}_{\boldsymbol{a},\boldsymbol{U}}(P)\geq 0 \ \ \text{and}  \ \  \sigma^{2}_{\Delta,\boldsymbol{U},\boldsymbol{M}}(P)-\sigma^{2}_{\Delta,\boldsymbol{U}}(P)\geq 0$. To do this, we first consider the efficient influence function. We have 
\begin{eqnarray} 
\psi_{P,\boldsymbol{a}}(\boldsymbol{U},\boldsymbol{M};\mathcal{G})
        &=&\dfrac{I_{\bar{\boldsymbol{a}}_{T}}(\boldsymbol{\bar{A}}_{T})}{\lambda_{\bar{\boldsymbol{a}}_T}(\bar{\boldsymbol{U}}_T,\bar{\boldsymbol{M}}_T;P)}\{Y-b_{\boldsymbol{a}_T}(\boldsymbol{U}_T,\boldsymbol{M}_T;P)\}\nonumber \\
        &+&\displaystyle{\sum_{t=0}^{T}}
	\dfrac{I_{\bar{\boldsymbol{a}}_{t-1}}(\bar{\boldsymbol{A}}_{t-1})\{b_{\boldsymbol{a}_t}(\boldsymbol{U}_t,\boldsymbol{M}_{t};P)-b_{\boldsymbol{a}_{t-1}}(\boldsymbol{U}_{t-1},\boldsymbol{M}_{t-1};P)\}}{\lambda_{\bar{\boldsymbol{a}}_{t-1}}(\bar{\boldsymbol{U}}_{t-1},\bar{\boldsymbol{M}}_{t-1};P)} \nonumber
\end{eqnarray} 
Denoting
\begin{align*}
    r_{_t}(\boldsymbol{\bar{A}}_{t-1},\bar{\boldsymbol{U}}_{t},\bar{\boldsymbol{M}}_{t-1};s_{\boldsymbol{a},t},P) & \equiv  \dfrac{I_{\bar{\boldsymbol{a}}_{t-1}}(\bar{\boldsymbol{A}}_{t-1})\{b_{\boldsymbol{a}_t}(\boldsymbol{U}_t,\boldsymbol{M}_{t};P)-b_{\boldsymbol{a}_{t-1}}(\boldsymbol{U}_{t-1},\boldsymbol{M}_{t-1};P)\}}{\lambda_{\bar{\boldsymbol{a}}_{t-1}}(\boldsymbol{\bar{U}}_{t-1},\bar{\boldsymbol{M}}_{t-1};P)}
\end{align*}
 we have
\begin{eqnarray}
 var_P\{\psi_{P,a}(\boldsymbol{U},\boldsymbol{M};\mathcal{G})\}
 &=&var_P\left[\dfrac{I_{\bar{\boldsymbol{a}}_{T}}(\boldsymbol{\bar{A}}_{T})}{\lambda_{\bar{\boldsymbol{a}}_T}(\bar{\boldsymbol{U}}_T,\bar{\boldsymbol{M}}_T;P)}\{Y-b_{\boldsymbol{a}_T}(\boldsymbol{U}_T,\boldsymbol{M}_T;P)\}\right] \nonumber \\
 &+&\sum_{t=0}^{T}var_P\{r_{_t}(\boldsymbol{\bar{A}}_{t-1},\bar{\boldsymbol{U}}_{t},\bar{\boldsymbol{M}}_{t-1};s_{\boldsymbol{a},t},P)\} \nonumber \\
 &+& 2 \sum_{t=0}^{T}\scalebox{0.75}{$ cov_P\left[\dfrac{I_{\bar{\boldsymbol{a}}_{T}}(\boldsymbol{\bar{A}}_{T})}{\lambda_{\bar{\boldsymbol{a}}_T}(\bar{\boldsymbol{U}}_T,\bar{\boldsymbol{M}}_T;P)}\{Y-b_{\boldsymbol{a}_T}(\boldsymbol{U}_T,\boldsymbol{M}_T;P)\},r_{_t}(\bar{\boldsymbol{A}}_{t-1},\bar{\boldsymbol{U}}_{t},\bar{\boldsymbol{M}}_{t-1};s_{\boldsymbol{a},t},P)\right]$} \nonumber \\
  &+&2\sum_{0\leq i<t\leq T}cov_P\left\{r_{_t}(\bar{\boldsymbol{A}}_t,\bar{\boldsymbol{U}}_{t},\bar{\boldsymbol{M}}_t;s_{\boldsymbol{a},t},P),r_{i}(\bar{\boldsymbol{A}}_{i},\bar{\boldsymbol{U}}_{i},\bar{\boldsymbol{M}}_{i};s_{\boldsymbol{a},i},P)\right\}. \nonumber
\end{eqnarray}
To work on this expression, we start by showing that some components involved are equal to 0. First,
{\footnotesize
 \begin{align*}
     &\mathbb{E}_P\{r_{_t}(\bar{\boldsymbol{A}}_{t-1},\bar{\boldsymbol{U}}_{t},\bar{\boldsymbol{M}}_{t-1};s_{\boldsymbol{a},t},P)\} \\
     & \qquad = \mathbb{E}_P[\mathbb{E}_P\{r_{_t}(\bar{\boldsymbol{A}}_{t-1},\bar{\boldsymbol{U}}_{t},\bar{\boldsymbol{M}}_{t-1};s_{\boldsymbol{a},t},P)\mid \bar{\boldsymbol{U}}_{t-1},\bar{\boldsymbol{M}}_{t-1},\bar{\boldsymbol{A}}_{t-1}\}]\\ %%%% Total expectation
  & \qquad =\mathbb{E}_P\left(\mathbb{E}_P\left[\dfrac{I_{\bar{\boldsymbol{a}}_{t-1}}(\bar{\boldsymbol{A}}_{t-1})\{b_{\boldsymbol{a}_t}(\boldsymbol{U}_t,\boldsymbol{M}_{t};P)-b_{\boldsymbol{a}_{t-1}}(\boldsymbol{U}_{t-1},\boldsymbol{M}_{t-1};P)\}}{\lambda_{\bar{\boldsymbol{a}}_{t-1}}(\bar{\boldsymbol{U}}_{t-1},\bar{\boldsymbol{M}}_{t-1};P)}\bigg | \bar{\boldsymbol{U}}_{t-1},\bar{\boldsymbol{M}}_{t-1},\bar{\boldsymbol{A}}_{t-1} =\bar{\boldsymbol{a}}_{t-1}\right]\right)\\ %%%% Par définition; l'espérance externe est une somme/intégrale sur les valeurs des variables de conditionnement, seulement le terme avec A_k-1 = a_k-1 \neq 0
  & \qquad= 
  \scalebox{0.9}{$
  \mathbb{E}_P\left(\dfrac{I_{\bar{\boldsymbol{a}}_{t-1}}(\bar{\boldsymbol{A}}_{t-1})}{\lambda_{\bar{\boldsymbol{a}}_{t-1}}(\bar{\boldsymbol{U}}_{t-1},\bar{\boldsymbol{M}}_{t-1};P)}\mathbb{E}_P\left[\{b_{\boldsymbol{a}_t}(\boldsymbol{U}_t, \boldsymbol{M}_t;P)-b_{\boldsymbol{a}_{t-1}}(\boldsymbol{U}_{t-1},\boldsymbol{M}_{t-1};P)\}\bigg | \bar{\boldsymbol{U}}_{t-1},\bar{\boldsymbol{M}}_{t-1},\bar{\boldsymbol{A}}_{t-1}=\bar{\boldsymbol{a}}_{t-1}\right]\right)$}\\ %%%% On sort de l'espérance une constante sachant les conditions; 
& \qquad=\mathbb{E}_P\left(\dfrac{I_{\bar{\boldsymbol{a}}_{t-1}}(\bar{\boldsymbol{A}}_{t-1})}{\lambda_{\bar{\boldsymbol{a}}_{t-1}}(\bar{\boldsymbol{U}}_{t-1},\bar{\boldsymbol{M}}_{t-1};P)}\mathbb{E}_P\left[\{b_{\boldsymbol{a}_t}(\boldsymbol{U}_t;P)-b_{\boldsymbol{a}_{t-1}}(\boldsymbol{U}_{t-1};P)\}\bigg | \bar{\boldsymbol{U}}_{t-1},\bar{\boldsymbol{M}}_{t-1},\bar{\boldsymbol{A}}_{t-1}=\bar{\boldsymbol{a}}_{t-1}\right]\right)\\ %%%% Resultat (13)
& \qquad=\mathbb{E}_P\left(\dfrac{I_{\bar{\boldsymbol{a}}_{t-1}}(\bar{\boldsymbol{A}}_{t-1})}{\lambda_{\bar{\boldsymbol{a}}_{t-1}}(\bar{\boldsymbol{U}}_{t-1},\bar{\boldsymbol{M}}_{t-1};P)}\mathbb{E}_P\left[\{b_{\boldsymbol{a}_t}(\boldsymbol{U}_t;P)-b_{\boldsymbol{a}_{t-1}}(\boldsymbol{U}_{t-1};P)\}\bigg | \boldsymbol{U}_{t-1},\boldsymbol{M}_{t-1},{A}_{t-1}={a}_{t-1}\right]\right)\\ %%%% Gk, Bk incluent tout le reste, mais on garder A_k parce que pourquoi pas
& \qquad=\mathbb{E}_P\left(\dfrac{I_{\bar{\boldsymbol{a}}_{t-1}}(\bar{\boldsymbol{A}}_{t-1})}{\lambda_{\bar{\boldsymbol{a}}_{t-1}}(\bar{\boldsymbol{U}}_{t-1},\bar{\boldsymbol{M}}_{t-1};P)}\mathbb{E}_P\left[\{b_{\boldsymbol{a}_t}(\boldsymbol{U}_t;P)-b_{\boldsymbol{a}_{t-1}}(\boldsymbol{U}_{t-1};P)\}\bigg | \boldsymbol{U}_{t-1},{A}_{t-1}={a}_{t-1}\right]\right)\\ %%%% On enleve Bk car les quantites ne dependent que de Gk
& \qquad= 
\scalebox{0.9}{$
\mathbb{E}_P\left\{\dfrac{I_{\bar{\boldsymbol{a}}_{t-1}}(\bar{\boldsymbol{A}}_{t-1})}{\lambda_{\bar{\boldsymbol{a}}_{t-1}}(\bar{\boldsymbol{U}}_{t-1},\bar{\boldsymbol{M}}_{t-1};P)}(\mathbb{E}_P\{b_{\boldsymbol{a}_t}(\boldsymbol{U}_t;P)\mid \boldsymbol{U}_{t-1},\boldsymbol{a}_{t-1}\}-\mathbb{E}_P[\mathbb{E}_P\{b_{\boldsymbol{a}_t}(\boldsymbol{U}_t;P)\mid \boldsymbol{U}_{t-1},\boldsymbol{a}_{t-1}\}\mid \boldsymbol{U}_{t-1},{a}_{t-1}])\right\}$}\\ %%%% On divise l'esperance en deux. Pour la 2e on applique la dfn de b_ak et la loi des esperances totales.
& \qquad= \mathbb{E}_P\left(\dfrac{I_{\bar{\boldsymbol{a}}_{t-1}}(\bar{\boldsymbol{A}}_{t-1})}{\lambda_{\bar{\boldsymbol{a}}_{t-1}}(\bar{\boldsymbol{U}}_{t-1},\bar{\boldsymbol{M}}_{t-1};P)}[\mathbb{E}_P\{b_{\boldsymbol{a}_t}(\boldsymbol{U}_t;P)\mid \boldsymbol{U}_{t-1},\boldsymbol{a}_{t-1}\}-\mathbb{E}_P\{b_{\boldsymbol{a}_t}(\boldsymbol{U}_t;P)\mid \boldsymbol{U}_{t-1},{a}_{t-1}\}]\right)\\
& \qquad=0.
\end{align*}}
Likewise, 
\begin{align*}
     &E_P\{r_{_t}(\bar{\boldsymbol{A}}_{t-1},\bar{\boldsymbol{U}}_{t},\bar{\boldsymbol{M}}_{t-1};s_{\boldsymbol{a},t},P)r_{i}(\bar{\boldsymbol{A}}_{i-1},\bar{\boldsymbol{U}}_{i},\bar{\boldsymbol{M}}_{i-1};s_{\boldsymbol{a},t},P)\} \\ %%%% esperances totales
     & \qquad = \mathbb{E}_P[\mathbb{E}_P\{r_{_t}(\bar{\boldsymbol{A}}_{t-1},\bar{\boldsymbol{U}}_{t},\bar{\boldsymbol{M}}_{t-1};s_{\boldsymbol{a},t},P)r_{i}(\bar{\boldsymbol{A}}_{i-1},\bar{\boldsymbol{U}}_{i},\bar{\boldsymbol{M}}_{i-1};s_{\boldsymbol{a},t},P)\mid \bar{\boldsymbol{U}}_{t-1},\bar{\boldsymbol{M}}_{t-1},\bar{\boldsymbol{A}}_{t-1}\}]\\ %%%% Meme principe que la precendente
     & \qquad=\mathbb{E}_P\left(\dfrac{I_{\bar{\boldsymbol{a}}_{t-1}}(\bar{\boldsymbol{A}}_{t-1})r_{i}(\bar{\boldsymbol{A}}_{i-1},\bar{\boldsymbol{U}}_{i},\bar{\boldsymbol{M}}_{i-1};s_{\boldsymbol{a},t},P)}{\lambda_{\bar{\boldsymbol{a}}_{t-1}}(\bar{\boldsymbol{U}}_{t-1},\bar{\boldsymbol{M}}_{t-1};P)} \times \right.  \\
     & \qquad \qquad \left. \mathbb{E}_P\left[\{b_{\boldsymbol{a}_t}(\boldsymbol{U}_t, \boldsymbol{M}_t;P)-b_{\boldsymbol{a}_{t-1}}(\boldsymbol{U}_{t-1},\boldsymbol{M}_{t-1};P)\}\bigg | \bar{\boldsymbol{U}}_{t-1},\bar{\boldsymbol{M}}_{t-1},\bar{\boldsymbol{A}}_{t-1}\right]\right)=0.
\end{align*}
We conclude that 
$cov_P\left\{r_{_t}(\bar{\boldsymbol{A}}_{t-1},\bar{\boldsymbol{U}}_{t},\bar{\boldsymbol{M}}_{t-1};s_{\boldsymbol{a},t},P),r_{i}(\bar{\boldsymbol{A}}_{i-1},\bar{\boldsymbol{U}}_{i},\bar{\boldsymbol{M}}_{i-1};s_{\boldsymbol{a},i},P)\right\}=0$. 

Furthermore,
{\footnotesize
\begin{align*}
     &cov_P\left[\dfrac{I_{\bar{\boldsymbol{a}}_{T}}(\bar{\boldsymbol{A}}_{T})}{\lambda_{\bar{\boldsymbol{a}}_T}(\bar{\boldsymbol{U}}_T,\bar{\boldsymbol{M}}_T;P)}\{Y-b_{\boldsymbol{a}_T}(\boldsymbol{U}_T,\boldsymbol{M}_T;P)\},r_{_t}(\bar{\boldsymbol{A}}_{t-1},\bar{\boldsymbol{U}}_{t},\bar{\boldsymbol{M}}_{t-1};s_{\boldsymbol{a},t},P)\right] \\
     & \qquad = \mathbb{E}_P\left(\mathbb{E}_P\left[\dfrac{I_{\bar{\boldsymbol{a}}_{T}}(\bar{\boldsymbol{A}}_{T})}{\lambda_{\bar{\boldsymbol{a}}_T}(\bar{\boldsymbol{U}}_T,\bar{\boldsymbol{M}}_T;P)}\{Y-b_{\boldsymbol{a}_T}(\boldsymbol{U}_T,\boldsymbol{M}_T;P)\}r_{t}(\bar{\boldsymbol{A}}_{t-1},\bar{\boldsymbol{U}}_{t},\bar{\boldsymbol{M}}_{t-1};s_{\boldsymbol{a},t},P)\bigg | \bar{\boldsymbol{U}}_{T},\bar{\boldsymbol{M}}_{T},\bar{\boldsymbol{A}}_{T}\right]\right)\\
     & \qquad=\mathbb{E}_P\left(\dfrac{I_{\bar{\boldsymbol{a}}_{T}}(\bar{\boldsymbol{A}}_{T})r_{t}(\bar{\boldsymbol{A}}_{t-1},\bar{\boldsymbol{U}}_{t},\bar{\boldsymbol{M}}_{t-1};s_{\boldsymbol{a},t},P)}{\lambda_{\bar{\boldsymbol{a}}_{T}}(\bar{\boldsymbol{U}}_{T},\bar{\boldsymbol{M}}_{T};P)}\mathbb{E}_P\left[\{Y-b_{\boldsymbol{a}_T}(\boldsymbol{U}_T,\boldsymbol{M}_T;P)\}\bigg | \bar{\boldsymbol{U}}_{T},\bar{\boldsymbol{M}}_{T},\bar{\boldsymbol{A}}_{T}\right]\right)\\
     & \qquad=\mathbb{E}_P\left(\dfrac{I_{\bar{\boldsymbol{a}}_{T}}(\bar{\boldsymbol{A}}_{T})r_{t}(\bar{\boldsymbol{A}}_{t-1},\bar{\boldsymbol{U}}_{t},\bar{\boldsymbol{M}}_{t-1};s_{\boldsymbol{a},t},P)}{\lambda_{\bar{\boldsymbol{a}}_{T}}(\bar{\boldsymbol{U}}_{T},\bar{\boldsymbol{M}}_{T};P)}\mathbb{E}_P\left[\{Y-b_{\boldsymbol{a}_T}(\boldsymbol{U}_T,\boldsymbol{M}_T;P)\}\bigg | \boldsymbol{U}_{T},\boldsymbol{M}_{T},{A}_{T}\right]\right)\\
     & \qquad=
     \scalebox{0.9}{$
     \mathbb{E}_P\left(\dfrac{I_{\bar{\boldsymbol{a}}_{T}}(\bar{\boldsymbol{A}}_{T})r_{t}(\bar{\boldsymbol{A}}_{t-1},\bar{\boldsymbol{U}}_{t},\bar{\boldsymbol{M}}_{t-1};s_{\boldsymbol{a},t},P)}{\lambda_{\bar{\boldsymbol{a}}_{T}}(\bar{\boldsymbol{U}}_{T},\bar{\boldsymbol{M}}_{T};P)}[\mathbb{E}_P(Y\big | \boldsymbol{U}_{T},\boldsymbol{M}_{T},\boldsymbol{A}_{T})-\mathbb{E}_P\{b_{\boldsymbol{a}_T}(\boldsymbol{U}_T,\boldsymbol{M}_T;P)\big | \boldsymbol{U}_{T},\boldsymbol{M}_{T},{A}_{T}\}]\right)$}\\
     & \qquad=\mathbb{E}_P\left[\dfrac{I_{\bar{\boldsymbol{a}}_{T}}(\bar{\boldsymbol{A}}_{T})r_{t}(\bar{\boldsymbol{A}}_{t-1},\bar{\boldsymbol{U}}_{t},\bar{\boldsymbol{M}}_{t-1};s_{\boldsymbol{a},t},P)}{\lambda_{\bar{\boldsymbol{a}}_{T}}(\bar{\boldsymbol{U}}_{T},\bar{\boldsymbol{M}}_{T};P)}\{b_{\boldsymbol{a}_T}(\boldsymbol{U}_T,\boldsymbol{M}_T;P)-b_{\boldsymbol{a}_T}(\boldsymbol{U}_T,\boldsymbol{M}_T;P)\}\right]=0.
\end{align*}}
We conclude that
\begin{align*}
     &var_P\{\psi_{P,\boldsymbol{a}}(\boldsymbol{U},\boldsymbol{M};\mathcal{G})\}-var_P\{\psi_{P,a}(\boldsymbol{U};\mathcal{G})\} \\
     & \qquad= var_P\left[\dfrac{I_{\bar{\boldsymbol{a}}_{T}}(\bar{\boldsymbol{A}}_{T})\{Y-b_{\boldsymbol{a}_T}(\boldsymbol{U}_T,\boldsymbol{M}_T;P)\}}{\lambda_{\bar{\boldsymbol{a}}_T}(\bar{\boldsymbol{U}}_T,\bar{\boldsymbol{M}}_T;P)}\right]-var_P\left[\dfrac{I_{\bar{\boldsymbol{a}}_{T}}(\bar{\boldsymbol{A}}_{T})\{Y-b_{\boldsymbol{a}_T}(\boldsymbol{U}_T;P)\}}{\lambda_{\bar{\boldsymbol{a}}_T}(\bar{\boldsymbol{U}}_T;P)}\right]\\
& \qquad + \sum_{t=0}^{T}\Biggl(var_P\left[\dfrac{I_{\bar{\boldsymbol{a}}_{t-1}}(\bar{\boldsymbol{A}}_{t-1}) \{b_{\boldsymbol{a}_t}(\boldsymbol{U}_{t},\boldsymbol{M}_{t};P)-b_{\boldsymbol{a}_{t-1}}(\boldsymbol{U}_{t-1},\boldsymbol{M}_{t-1};P)\}}{\lambda_{\bar{a}_{t-1}}(\bar{\boldsymbol{U}}_{t-1},\bar{\boldsymbol{M}}_{t-1};P)}\right] \\
& \qquad \qquad -var_P\left[\dfrac{I_{\bar{\boldsymbol{a}}_{t-1}}(\bar{\boldsymbol{A}}_{t-1}) \{b_{\boldsymbol{a}_t}(\boldsymbol{U}_{t};P)-b_{\boldsymbol{a}_{t-1}}(\boldsymbol{U}_{t-1};P)\}}{\lambda_{\bar{a}_{t-1}}(\bar{\boldsymbol{U}}_{t-1};P)}\right]\Biggl).
\end{align*}
Next, we show that
\begin{equation*}
\scalebox{0.92}{$
var_P\left[\dfrac{I_{\bar{\boldsymbol{a}}_{t-1}}(\bar{\boldsymbol{A}}_{t-1}) \{b_{\boldsymbol{a}_t}(\boldsymbol{U}_{t};P)-b_{\boldsymbol{a}_{t-1}}(\boldsymbol{U}_{t-1};P)\}}{\lambda_{\bar{a}_{t-1}}(\bar{\boldsymbol{U}}_{t-1},\bar{\boldsymbol{M}}_{t-1};P)}\right]\geq var_P\left[\dfrac{I_{\bar{\boldsymbol{a}}_{t-1}}(\bar{\boldsymbol{A}}_{t-1}) \{b_{\boldsymbol{a}_t}(\boldsymbol{U}_{t};P)-b_{\boldsymbol{a}_{t-1}}(\boldsymbol{U}_{t-1};P)\}}{\lambda_{\bar{a}_{t-1}}(\bar{\boldsymbol{U}}_{t-1};P)}\right]$},    
\end{equation*}
for $t \in \{0,...,T\}$, where $\lambda_{\bar{a}_{t-1}}(\bar{\boldsymbol{U}}_{t-1};P)=\displaystyle{\prod_{k=0}^{t-1}}P(A_k=a_k\mid \boldsymbol{U}_k)$.

To do this, we first show that for any  $t \in \{0,...,T\}$ and $T\geq1$
\begin{equation}
\scalebox{0.9}{$
\mathbb{E}_P\left\{ \dfrac{1}{\lambda_{\bar{a}_{t-1}}(\boldsymbol{U}_{t-1},\boldsymbol{M}_{t-1};P)}\middle | \bar{\boldsymbol{U}}_{t-1} = \bar{\boldsymbol{u}}_{t-1},\boldsymbol{\bar{A}}_{t-1}=\boldsymbol{\bar{a}}_{t-1}\right\} = \dfrac{1}{\lambda_{\bar{a}_{t-1}}(\bar{\boldsymbol{U}}_{t-1}=\bar{\boldsymbol{u}}_{t-1},\boldsymbol{\bar{A}}_{t-2} = \boldsymbol{\bar{a}}_{t-2};P)}$}.
\end{equation}
and 
$$
    b_{\boldsymbol{a}_{t}}(\bar{\boldsymbol{U}}_{t};P)
    = b_{\boldsymbol{a}_{t}}(\boldsymbol{U}_{t};P). 
$$
On one hand for any $t \in \{0,...,T\}$ 
\begin{align*}
&  \mathbb{E}_P\left\{\dfrac{1}{\lambda_{\boldsymbol{\bar{a}}_{t-1}}(\boldsymbol{U}_{t-1},\boldsymbol{M}_{t-1};P)}\middle | \bar{\boldsymbol{U}}_{t-1} = \bar{\boldsymbol{u}}_{t-1},\bar{A}_{t-1}=\bar{a}_{t-1}\right\}\pi_{{a}_{t-1}}(\bar{\boldsymbol{U}}_{t-1} = \bar{\boldsymbol{u}}_{t-1},\boldsymbol{\bar{A}}_{t-2} =\boldsymbol{\bar{a}}_{t-2};P)\\
& = \scalebox{0.9}{$\mathbb{E}_P\left\{\dfrac{1}{\lambda_{\boldsymbol{\bar{a}}_{t-1}}(\boldsymbol{U}_{t-1},\boldsymbol{M}_{t-1};P)}\middle | \bar{\boldsymbol{U}}_{t-1} = \bar{\boldsymbol{u}}_{t-1},\boldsymbol{\bar{A}}_{t-1}=\boldsymbol{\bar{a}}_{t-1}\right\}
P(A_{t-1}=a_{t-1}\mid \bar{\boldsymbol{U}}_{t-1} = \bar{\boldsymbol{u}}_{t-1},\boldsymbol{\bar{A}}_{t-2}=\boldsymbol{\bar{a}}_{t-2})$}\\
& = \mathbb{E}_P\left\{\dfrac{1}{\lambda_{\boldsymbol{\bar{a}}_{t-2}}(\boldsymbol{U}_{t-2},\boldsymbol{M}_{t-2};P)}\times \dfrac{I_{a_{t-1}}(A_{t-1})}{\pi_{{a}_{t-1}}(\boldsymbol{U}_{t-1},\boldsymbol{M}_{t-1};P)}\middle | \bar{\boldsymbol{U}}_{t-1} = \bar{\boldsymbol{u}}_{t-1},\boldsymbol{\bar{A}}_{t-2}=\boldsymbol{\bar{a}}_{t-2}\right\}\\ %%%% Separation du lambda en deux termes; + E(XI(A = a)\mid B) = E(X\mid A=a, B)P(A = a\mid B)
& =\mathbb{E}_P\left[\dfrac{1}{\lambda_{\boldsymbol{\bar{a}}_{t-2}}(\boldsymbol{U}_{t-2},\boldsymbol{M}_{t-2};P)}\times \dfrac{E_P\{I_{a_{t-1}}(A_{t-1})\mid \boldsymbol{U}_{t-1},\boldsymbol{M}_{t-1}\}}{\pi_{{a}_{t-1}}(\boldsymbol{U}_{t-1},\boldsymbol{M}_{t-1};P)}\middle | \bar{\boldsymbol{U}}_{t-1} = \bar{\boldsymbol{u}}_{t-1},\boldsymbol{\bar{A}}_{t-2}=\boldsymbol{\bar{a}}_{t-2}\right]\\ % Esperances totales et on fait sortir tout sauf l'indicatrice, comme tout le reste est constant sachant la condition
& = \mathbb{E}_P\left\{ \dfrac{1}{\lambda_{\boldsymbol{\bar{a}}_{t-2}}(\boldsymbol{U}_{t-2},\boldsymbol{M}_{t-2};P)}\times \dfrac{P(A_{t-1}=a_{t-1}\mid \boldsymbol{U}_{t-1},\boldsymbol{M}_{t-1})}{\pi_{{a}_{t-1}}(\boldsymbol{U}_{t-1},\boldsymbol{M}_{t-1};P)}\middle | \bar{\boldsymbol{U}}_{t-1} = \bar{\boldsymbol{u}}_{t-1},\boldsymbol{\bar{A}}_{t-2}=\boldsymbol{\bar{a}}_{t-2}\right\}\\ %%%% Le reste suit par les definitions
& = \mathbb{E}_P\left\{\dfrac{1}{\lambda_{\boldsymbol{\bar{a}}_{t-1}}(\boldsymbol{U}_{t-2},\boldsymbol{M}_{t-2};P)}\times \dfrac{\pi_{{a}_{t-1}}(\boldsymbol{U}_{t-1},\boldsymbol{M}_{t-1};P)}{\pi_{{a}_{t-1}}(\boldsymbol{U}_{t-1},\boldsymbol{M}_{t-1};P)}\middle | \bar{\boldsymbol{U}}_{t-1} = \bar{\boldsymbol{u}}_{t-1},\boldsymbol{\bar{A}}_{t-2}=\boldsymbol{\bar{a}}_{t-2}\right\}\\
& = \mathbb{E}_P\left\{ \dfrac{1}{\lambda_{\boldsymbol{\bar{a}}_{t-2}}(\boldsymbol{U}_{t-2},\boldsymbol{M}_{t-2};P)}\middle | \bar{\boldsymbol{U}}_{t-1} = \bar{\boldsymbol{u}}_{t-1},\boldsymbol{\bar{A}}_{t-2}=\boldsymbol{\bar{a}}_{t-2}\right\}\\
& = \mathbb{E}_P\left\{ \dfrac{1}{\lambda_{\boldsymbol{\bar{a}}_{t-2}}(\boldsymbol{U}_{t-2},\boldsymbol{M}_{t-2};P)}\middle | \bar{\boldsymbol{U}}_{t-2} = \bar{\boldsymbol{u}}_{t-2},\boldsymbol{\bar{A}}_{t-2}=\boldsymbol{\bar{a}}_{t-2}\right\}.
\end{align*}
where the last equality follows from (\ref{eq18}).
%and remarking that $Q_t \equiv b_{a_t}(\boldsymbol{Z}_t; P)$.

On the other hand
$$
    b_{{a}_{t}}(\bar{\boldsymbol{U}}_{t};P)
    =b_{{a}_{t}}(\boldsymbol{U}_{t},\bar{\boldsymbol{U}}_{t-1};P). 
$$
If $\bar{\boldsymbol{U}}_{t-1} \subset \boldsymbol{U}_{t}$ then $b_{{a}_{t}}(\bar{\boldsymbol{U}}_{t};P)=b_{{a}_{t}}(\boldsymbol{U}_{t};P)$

If $\bar{\boldsymbol{U}}_{t-1} \subset \boldsymbol{M}_{t}$ then we obtain $b_{{a}_{t}}(\bar{\boldsymbol{U}}_{t};P)=b_{{a}_{t}}(\boldsymbol{U}_{t};P)$ by (\ref{eq:lemma2.3}) applied to  $\bar{\boldsymbol{U}}_{t-1} \subset \boldsymbol{M}_{t}$.

Then, for any $t \in \{0,...,T\}$,
	\begin{eqnarray}
		& & \mathbb{E}_P\left\{ \dfrac{1}{\lambda_{\boldsymbol{\bar{a}}_{t-1}}(\boldsymbol{U}_{t-1},\boldsymbol{M}_{t-1};P)}\middle | \bar{\boldsymbol{U}}_{t-1} = \bar{\boldsymbol{u}}_{t-1},\boldsymbol{\bar{A}}_{t-1}=\boldsymbol{\bar{a}}_{t-1}\right\} \nonumber \\
		& & 
        \scalebox{0.9}{$
        \qquad = \dfrac{1}{\pi_{{a}_t}(\bar{\boldsymbol{U}}_{t-1} = \bar{\boldsymbol{u}}_{t-1},\boldsymbol{\bar{A}}_{t-2} = \boldsymbol{\bar{a}}_{t-2};P)}
		\mathbb{E}_P\left\{\dfrac{1}{\lambda_{\bar{a}_{t-2}}(\boldsymbol{U}_{t-2},\boldsymbol{M}_{t-2};P)}\middle | \bar{\boldsymbol{U}}_{t-2} = \bar{\boldsymbol{u}}_{t-2},\boldsymbol{\bar{A}}_{t-2}=\boldsymbol{\bar{a}}_{t-2}\right\}$} \nonumber \\
		& & \qquad =
		\dfrac{1}{\pi_{{a}_{t-1}}(\bar{\boldsymbol{U}}_{t-1} = \bar{\boldsymbol{u}}_{t-1},\boldsymbol{\bar{A}}_{t-2} = \boldsymbol{\bar{a}}_{t-2};P)}\times	\dfrac{1}{\pi_{{a}_{t-1}}(\bar{\boldsymbol{U}}_{t-2} = \bar{\boldsymbol{u}}_{t-2},\boldsymbol{\bar{A}}_{t-3} = \boldsymbol{\bar{a}}_{t-3};P)} \times \nonumber \\
		&& \qquad \qquad \mathbb{E}_P\left\{\dfrac{1}{\lambda_{\boldsymbol{\bar{a}}_{t-3}}(\boldsymbol{U}_{t-3},\boldsymbol{M}_{t-3};P)}\middle | \bar{\boldsymbol{U}}_{t-3} = \bar{\boldsymbol{u}}_{t-3},\boldsymbol{\bar{A}}_{t-3}=\boldsymbol{\bar{a}}_{t-3}\right\} \nonumber \\
		& & \qquad =
		\dfrac{1}{\pi_{{a}_{t-1}}(\bar{\boldsymbol{U}}_{t-1} = \bar{\boldsymbol{u}}_{t-1},\boldsymbol{\bar{A}}_{t-2} = \boldsymbol{\bar{a}}_{t-2};P)}\times	\dfrac{1}{\pi_{{a}_{t-1}}(\bar{\boldsymbol{U}}_{t-2} = \bar{\boldsymbol{u}}_{t-2},\boldsymbol{\bar{A}}_{t-3} = \boldsymbol{\bar{a}}_{t-3};P)} \nonumber \\        
        & & \qquad \qquad \times\cdots\times \dfrac{1}{\pi_{{a}_{0}}(\boldsymbol{U}_{0} = \boldsymbol{u}_{0};P)} \nonumber \\
		& & \qquad =\dfrac{1}{\lambda_{\bar{a}_{t-1}}(\bar{\boldsymbol{U}}_{t-1} = \bar{\boldsymbol{u}}_{t-1},\boldsymbol{\bar{A}}_{t-2} = \boldsymbol{\bar{a}}_{t-2};P)}. \nonumber
	\end{eqnarray} 
	Now,
 \begin{align*}
& var_P\left[\dfrac{I_{\bar{\boldsymbol{a}}_{t-1}}(\bar{\boldsymbol{A}}_{t-1})\{b_{{a}_t}(\boldsymbol{U}_{t};P)-b_{{a}_{t-1}}(\boldsymbol{U}_{t-1};P)\}}{\lambda_{\boldsymbol{\bar{a}}_{t-1}}(\boldsymbol{U}_{t-1},\boldsymbol{M}_{t-1};P)}\right]\\
& \qquad = \mathbb{E}_P\left(var_P\left[\dfrac{I_{\bar{\boldsymbol{a}}_{t-1}}(\bar{\boldsymbol{A}}_{t-1})\{b_{{a}_t}(\boldsymbol{U}_{t};P)-b_{{a}_{t-1}}(\boldsymbol{U}_{t-1};P)\}}{\lambda_{\boldsymbol{\bar{a}}_{t-1}}(\boldsymbol{U}_{t-1},\boldsymbol{M}_{t-1};P)}\middle | \bar{\boldsymbol{U}}_{t},\boldsymbol{\bar{A}}_{t-1}=\boldsymbol{\bar{a}}_{t-1} \right]\right) \\
& \qquad \qquad + var_P\left(\mathbb{E}_P\left[\dfrac{I_{\bar{\boldsymbol{a}}_{t-1}}(\bar{\boldsymbol{A}}_{t-1})\{b_{{a}_t}(\boldsymbol{U}_{t};P)-b_{{a}_{t-1}}(\boldsymbol{U}_{t-1};P)\}}{\lambda_{\bar{a}_{t-1}}(\boldsymbol{U}_{t-1},\boldsymbol{M}_{t-1};P)}\middle | \bar{\boldsymbol{U}}_{t},\boldsymbol{\bar{A}}_{t-1}=\boldsymbol{\bar{a}}_{t-1} \right]\right) \\
& \qquad \geq
	var_P\left(\mathbb{E}_P\left[\dfrac{I_{\bar{\boldsymbol{a}}_{t-1}}(\bar{\boldsymbol{A}}_{t-1})\{b_{{a}_t}(\boldsymbol{U}_{t};P)-b_{{a}_{t-1}}(\boldsymbol{U}_{t-1};P)\}}{\lambda_{\boldsymbol{\bar{a}}_{t-1}}(\boldsymbol{U}_{t-1},\boldsymbol{M}_{t-1};P)}\middle | \bar{\boldsymbol{U}}_{t},\boldsymbol{\bar{A}}_{t-1}=\boldsymbol{\bar{a}}_{t-1} \right]\right). 
 \end{align*}
Moreover,
\begin{align*}
 &\mathbb{E}_P\left[\dfrac{I_{\bar{\boldsymbol{a}}_{t-1}}(\bar{\boldsymbol{A}}_{t-1})\{b_{{a}_t}(\boldsymbol{U}_{t};P)-b_{{a}_{t-1}}(\boldsymbol{U}_{t-1};P)\}}{\lambda_{\bar{a}_{t-1}}(\boldsymbol{U}_{t-1},\boldsymbol{M}_{t-1};P)}\middle | \bar{\boldsymbol{U}}_{t}= \bar{\boldsymbol{u}}_{t},\boldsymbol{\bar{A}}_{t-1}=\boldsymbol{\bar{a}}_{t-1} \right]\\
 & \qquad  = \mathbb{E}_P\left[\dfrac{I_{\bar{\boldsymbol{a}}_{t-1}}(\bar{\boldsymbol{A}}_{t-1})\{b_{{a}_t}(\boldsymbol{U}_{t};P)-b_{{a}_{t-1}}(\boldsymbol{U}_{t-1};P)\}}{\lambda_{\boldsymbol{\bar{a}}_{t-1}}(\boldsymbol{U}_{t-1},\boldsymbol{M}_{t-1};P)}\middle | \bar{\boldsymbol{U}}_{t}=\bar{\boldsymbol{u}}_{t} ,\boldsymbol{\bar{A}}_{t-1}=\boldsymbol{\bar{a}}_{t-1} \right]\\
 & \qquad = \mathbb{E}_P\left[\dfrac{I_{\bar{\boldsymbol{a}}_{t-1}}(\bar{\boldsymbol{A}}_{t-1})\{b_{{a}_t}(\bar{\boldsymbol{U}}_t;P)-b_{{a}_{t-1}}(\bar{\boldsymbol{U}}_{t-1};P)\}}{\lambda_{\boldsymbol{\bar{a}}_{t-1}}(\boldsymbol{U}_{t-1},\boldsymbol{M}_{t-1};P)}\middle | \bar{\boldsymbol{U}}_{t}=\bar{\boldsymbol{u}}_{t},\boldsymbol{\bar{A}}_{t-1}=\boldsymbol{\bar{a}}_{t-1} \right] \\
 & \qquad = \scalebox{0.9}{$I_{\bar{\boldsymbol{a}}_{t-1}}(\bar{\boldsymbol{A}}_{t-1})\{b_{{a}_t}(\bar{\boldsymbol{U}}_t;P)-b_{{a}_{t-1}}(\bar{\boldsymbol{U}}_{t-1};P)\}\mathbb{E}_P\left\{\dfrac{1}{\lambda_{\bar{a}_{t-1}}(\bar{\boldsymbol{U}}_{t-1},\bar{\boldsymbol{M}}_{t-1};P)}\middle | \bar{\boldsymbol{U}}_{t}=\bar{\boldsymbol{u}}_{t},\boldsymbol{\bar{A}}_{t-1}=\boldsymbol{\bar{a}}_{t-1} \right\}$}\\
 & \qquad = \scalebox{0.9}{$I_{\bar{\boldsymbol{a}}_{t-1}}(\bar{\boldsymbol{A}}_{t-1})\{b_{{a}_t}(\boldsymbol{U}_t;P)-b_{{a}_{t-1}}(\boldsymbol{U}_{t-1};P)\}\mathbb{E}_P\left\{\dfrac{1}{\lambda_{\bar{a}_{t-1}}(\boldsymbol{U}_{t-1},\boldsymbol{M}_{t-1};P)}\middle | \bar{\boldsymbol{U}}_{t} = \bar{\boldsymbol{u}}_{t},\boldsymbol{\bar{A}}_{t-1}=\boldsymbol{\bar{a}}_{t-1} \right\}$}\\
 & \qquad = \dfrac{I_{\bar{a}_{t-1}}(\bar{a}_{t-1})\{b_{{a}_t}(\boldsymbol{U}_t;P)-b_{{a}_{t-1}}(\boldsymbol{U}_{t-1};P)\}}{\lambda_{\bar{a}_{t-1}}(\bar{\boldsymbol{U}}_{t-1} = \bar{\boldsymbol{u}}_{t-1},\boldsymbol{\bar{A}}_{t-2} = \boldsymbol{\bar{a}}_{t-2};P)}.
\end{align*}
Thus, 
\begin{equation*}
\scalebox{0.9}{$
   var_P\left[\dfrac{I_{\bar{\boldsymbol{a}}_{t-1}}(\bar{\boldsymbol{A}}_{t-1})\{b_{\boldsymbol{a}_t}(\boldsymbol{U}_{t};P)-b_{\boldsymbol{a}_{t-1}}(\boldsymbol{U}_{t-1};P)\}}{\lambda_{\bar{a}_{t-1}}(\boldsymbol{U}_{t-1},\boldsymbol{M}_{t-1};P)}\right]\geq var_P\left[\dfrac{I_{\bar{\boldsymbol{a}}_{t-1}}(\bar{\boldsymbol{A}}_{t-1})\{b_{\boldsymbol{a}_t}(\boldsymbol{U}_{t};P)-b_{\boldsymbol{a}_{t-1}}(\boldsymbol{U}_{t-1};P)\}}{\lambda_{\bar{a}_{t-1}}(\bar{\boldsymbol{U}}_{t-1},\bar{A}_{t-2};P)}\right]$}
\end{equation*}
$\Longleftrightarrow$ 
\begin{equation*}
\scalebox{0.9}{$
   var_P\left[\dfrac{I_{\bar{\boldsymbol{a}}_{t-1}}(\bar{\boldsymbol{A}}_{t-1})\{b_{\boldsymbol{a}_t}(\boldsymbol{U}_{t};P)-b_{\boldsymbol{a}_{t-1}}(\boldsymbol{U}_{t-1};P)\}}{\prod_{k=0}^{t-1}P(A_k=a_k\mid \boldsymbol{U}_k,\boldsymbol{M}_k)}\right]\geq var_P\left[\dfrac{I_{\bar{\boldsymbol{a}}_{t-1}}(\bar{\boldsymbol{A}}_{t-1})\{b_{\boldsymbol{a}_t}(\boldsymbol{U}_{t};P)-b_{\boldsymbol{a}_{t-1}}(\boldsymbol{U}_{t-1};P)\}}{\prod_{k=0}^{t-1}P(A_k=a_k\mid \boldsymbol{U}_k,\bar{\boldsymbol{U}}_{k-1},\bar{A}_{k-1})}\right].(*)$}
\end{equation*}
Next, we show that 
\begin{equation*}
\scalebox{0.9}{$
    var_P\left[\dfrac{I_{\bar{\boldsymbol{a}}_{t-1}}(\bar{\boldsymbol{A}}_{t-1})\{b_{\boldsymbol{a}_t}(\boldsymbol{U}_{t};P)-b_{\boldsymbol{a}_{t-1}}(\boldsymbol{U}_{t-1};P)\}}{\prod_{k=0}^{t-1}P(A_k=a_k\mid \boldsymbol{U}_k,\boldsymbol{M}_{k})}\right]\geq var_P\left[\dfrac{I_{\bar{\boldsymbol{a}}_{t-1}}(\bar{\boldsymbol{A}}_{t-1})\{b_{\boldsymbol{a}_t}(\boldsymbol{U}_{t};P)-b_{\boldsymbol{a}_{t-1}}(\boldsymbol{U}_{t-1};P)\}}{\prod_{k=0}^{t-1}P(A_k=a_k\mid \boldsymbol{U}_k)}\right].$}
\end{equation*}

\noindent To do this, we first show that 
 \begin{align*}
 \scalebox{0.9}{$
 var_P\left[\dfrac{I_{\bar{\boldsymbol{a}}_{t-1}}(\bar{\boldsymbol{A}}_{t-1})\{b_{\boldsymbol{a}_t}(\boldsymbol{U}_{t};P)-b_{\boldsymbol{a}_{t-1}}(\boldsymbol{U}_{t-1};P)\}}{\prod_{k=0}^{t-1}P(A_k=a_k\mid \boldsymbol{U}_k,\bar{\boldsymbol{U}}_{k-1},\bar{A}_{k-1})}\right]\geq var_P\left[\dfrac{I_{\bar{\boldsymbol{a}}_{t-1}}(\bar{\boldsymbol{A}}_{t-1})\{b_{\boldsymbol{a}_t}(\boldsymbol{U}_{t};P)-b_{\boldsymbol{a}_{t-1}}(\boldsymbol{U}_{t-1};P)\}}{\prod_{k=0}^{t-1}P(A_k=a_k\mid \boldsymbol{U}_k)}\right]$}  
 \end{align*}
 effectively bounding the right-hand side of $(*)$. 

\noindent First case: if ${(\bar{\boldsymbol{U}}_{k-1},\bar{A}_{k-1})\subseteq \boldsymbol{U}_{k}}$, then $(\boldsymbol{U}_{k},\bar{\boldsymbol{U}}_{k-1},\bar{A}_{k-1}) = \boldsymbol{U}_{k}$ and we have 
\begin{equation*}
 \scalebox{0.9}{$
    var_P\left[\dfrac{I_{\bar{\boldsymbol{a}}_{t-1}}(\bar{\boldsymbol{A}}_{t-1})\{b_{\boldsymbol{a}_t}(\boldsymbol{U}_{t};P)-b_{\boldsymbol{a}_{t-1}}(\boldsymbol{U}_{t-1};P)\}}{\prod_{k=0}^{t-1}P(A_k=a_k\mid \boldsymbol{U}_k,\bar{\boldsymbol{U}}_{k-1},\bar{A}_{k-1})}\right] = var_P\left[\dfrac{I_{\bar{\boldsymbol{a}}_{t-1}}(\bar{\boldsymbol{A}}_{t-1})\{b_{\boldsymbol{a}_t}(\boldsymbol{U}_{t};P)-b_{\boldsymbol{a}_{t-1}}(\boldsymbol{U}_{t-1};P)\}}{\prod_{k=0}^{t-1}P(A_k=a_k\mid \boldsymbol{U}_k)}\right].$}
\end{equation*}
Second case: if $(\bar{\boldsymbol{U}}_{k-1},\bar{A}_{k-1})\subseteq \boldsymbol{M}_{k}$, then denoting $\boldsymbol{M}_{k}^*:=(\bar{\boldsymbol{U}}_{k-1},\bar{A}_{k-1})$ and employing arguments similar to those that led to relation $(*)$ yield: 

\scalebox{0.9}{$
\begin{aligned}
var_P\left[\dfrac{I_{\bar{\boldsymbol{a}}_{t-1}}(\bar{\boldsymbol{A}}_{t-1})\{b_{\boldsymbol{a}_t}(\boldsymbol{U}_{t};P)-b_{\boldsymbol{a}_{t-1}}(\boldsymbol{U}_{t-1};P)\}}{\prod_{k=0}^{t-1}P(A_k=a_k\mid \boldsymbol{U}_k,\bar{\boldsymbol{U}}_{k-1},\bar{A}_{k-1})}\right] &= var_P\left[\dfrac{I_{\bar{\boldsymbol{a}}_{t-1}}(\bar{\boldsymbol{A}}_{t-1})\{b_{\boldsymbol{a}_t}(\boldsymbol{U}_{t};P)-b_{\boldsymbol{a}_{t-1}}(\boldsymbol{U}_{t-1};P)\}}{\prod_{k=0}^{t-1}P(A_k=a_k\mid \boldsymbol{U}_k,\boldsymbol{M}_k^*)}\right] \\
& \geq var_P\left[\dfrac{I_{\bar{\boldsymbol{a}}_{t-1}}(\bar{\boldsymbol{A}}_{t-1})\{b_{\boldsymbol{a}_t}(\boldsymbol{U}_{t};P)-b_{\boldsymbol{a}_{t-1}}(\boldsymbol{U}_{t-1};P)\}}{\prod_{k=0}^{t-1}P(A_k=a_k\mid \boldsymbol{U}_k)}\right].
\end{aligned}
$}

Third case: if $\bar{\boldsymbol{U}}_{k-1}\subseteq \boldsymbol{U}_{k}$ and $\bar{A}_{k-1}\in \boldsymbol{M}_{k}$, then $(\boldsymbol{U}_{k},\bar{\boldsymbol{U}}_{k-1})=\boldsymbol{U}_{k}$ and similarly to the second case we get: 
\begin{eqnarray*}
\scalebox{0.9}{$
   var_P\left[\dfrac{I_{\bar{\boldsymbol{a}}_{t-1}}(\bar{\boldsymbol{A}}_{t-1})\{b_{\boldsymbol{a}_t}(\boldsymbol{U}_{t};P)-b_{\boldsymbol{a}_{t-1}}(\boldsymbol{U}_{t-1};P)\}}{\prod_{k=0}^{t-1}P(A_k=a_k\mid \boldsymbol{U}_k,\bar{\boldsymbol{U}}_{k-1},\bar{A}_{k-1})}\right]\geq var_P\left[\dfrac{I_{\bar{\boldsymbol{a}}_{t-1}}(\bar{\boldsymbol{A}}_{t-1})\{b_{\boldsymbol{a}_t}(\boldsymbol{U}_{t};P)-b_{\boldsymbol{a}_{t-1}}(\boldsymbol{U}_{t-1};P)\}}{\prod_{k=0}^{t-1}P(A_k=a_k\mid \boldsymbol{U}_k)}\right].$}
\end{eqnarray*}

Fourth case: if $\bar{A}_{k-1}\in \boldsymbol{U}_{k}$ and $\bar{\boldsymbol{U}}_{k-1}\subseteq \boldsymbol{M}_{k}$, then $(\boldsymbol{U}_{k},\bar{A}_{k-1})=\boldsymbol{U}_{k}$ and arguments similar to those employed in the second case once again yield:
\begin{equation*}
\scalebox{0.9}{$
    var_P\left[\dfrac{I_{\bar{\boldsymbol{a}}_{t-1}}(\bar{\boldsymbol{A}}_{t-1})\{b_{\boldsymbol{a}_t}(\boldsymbol{U}_{t};P)-b_{\boldsymbol{a}_{t-1}}(\boldsymbol{U}_{t-1};P)\}}{\prod_{k=0}^{t-1}P(A_k=a_k\mid \boldsymbol{U}_k,\bar{\boldsymbol{U}}_{k-1},\bar{A}_{k-1})}\right]\geq var_P\left[\dfrac{I_{\bar{\boldsymbol{a}}_{t-1}}(\bar{\boldsymbol{A}}_{t-1})\{b_{\boldsymbol{a}_t}(\boldsymbol{U}_{t};P)-b_{\boldsymbol{a}_{t-1}}(\boldsymbol{U}_{t-1};P)\}}{\prod_{k=0}^{t-1}P(A_k=a_k\mid \boldsymbol{U}_k)}\right].$}
\end{equation*}

As such, in all cases,
\begin{equation*}
\scalebox{0.8}{$
    var_P\left[\dfrac{I_{\bar{\boldsymbol{a}}_{t-1}}(\bar{\boldsymbol{A}}_{t-1})\{b_{\boldsymbol{a}_t}(\boldsymbol{U}_{t};P)-b_{\boldsymbol{a}_{t-1}}(\boldsymbol{U}_{t-1};P)\}}{\prod_{k=0}^{t-1}P(A_k=a_k\mid \boldsymbol{U}_k,\bar{\boldsymbol{U}}_{k-1},\bar{A}_{k-1})}\right]\geq var_P\left[\dfrac{I_{\bar{\boldsymbol{a}}_{t-1}}(\bar{\boldsymbol{A}}_{t-1})\{b_{\boldsymbol{a}_t}(\boldsymbol{U}_{t};P)-b_{\boldsymbol{a}_{t-1}}(\boldsymbol{U}_{t-1};P)\}}{\prod_{k=0}^{t-1}P(A_k=a_k\mid \boldsymbol{U}_k)}\right].(**)$}
\end{equation*}

Using $(*)$ and $(**)$ we conclude that, 
\begin{equation*}
    var_P\{\psi_{P,\boldsymbol{a}}(\boldsymbol{U},\boldsymbol{M};\mathcal{G})\}-var_P\{\psi_{P,a}(\boldsymbol{U};\mathcal{G})\}=\sigma^{2}_{\boldsymbol{\boldsymbol{a}},\boldsymbol{U},\boldsymbol{M}}(P)-\sigma^{2}_{\boldsymbol{\boldsymbol{a}},\boldsymbol{U}}(P)\geq0.
\end{equation*}
We derive the formula for  $\sigma^{2}_{\boldsymbol{\Delta},\boldsymbol{U},\boldsymbol{M}}(P)-\sigma^{2}_{\boldsymbol{\Delta},\boldsymbol{U}}(P)$ analogously. Specifically, we note $\boldsymbol{c}\equiv (c_a)_{a \in \mathcal{\boldsymbol{A}}}$. For all $\boldsymbol{Z}$, we define $\psi_{P}(\boldsymbol{Z};\mathcal{G})\equiv \{\psi_{P,a}(\boldsymbol{Z};\mathcal{G})\}_{a \in \mathcal{\boldsymbol{A}}}$ and  $\sum_{a \in \mathcal{\boldsymbol{A}}} c_a\psi_{P,a}(\boldsymbol{Z};\mathcal{G})=\boldsymbol{c}^{T}\psi_{P,\Delta}(\boldsymbol{Z};\mathcal{G})$
\begin{eqnarray}
&& var_P\{\boldsymbol{c}^{T}\psi_{P,\Delta}(\boldsymbol{U},\boldsymbol{M};\mathcal{G})\}-var_P\{\boldsymbol{c}^{T}\psi_{P,\Delta}(\boldsymbol{U};\mathcal{G})\} \\
&&\qquad =	var_P\left[\sum_{a \in \mathcal{\boldsymbol{A}}} c_a\dfrac{I_{\bar{\boldsymbol{a}}_{T}}(\bar{\boldsymbol{A}}_{T})}{\lambda_{\bar{a}_{T}}(\boldsymbol{U}_{T},\boldsymbol{M}_{T};P)}\{Y-b_{\boldsymbol{a}_T}(\boldsymbol{U}_{T};P)\}\right] \nonumber \\
&&\qquad \qquad - var_P\left[\sum_{a \in \mathcal{\boldsymbol{A}}} c_a\dfrac{I_{\bar{a}_{T}}(\bar{A}_{T})}{\lambda_{\bar{a}_{T}}(\bar{\boldsymbol{U}}_{T};P)}\{Y-b_{\boldsymbol{a}_T}(\boldsymbol{U}_{T};P)\}\right] \nonumber \\
&& \qquad \qquad + \sum_{t=0}^{T}var_P\left[\sum_{a \in \mathcal{\boldsymbol{A}}} c_a\dfrac{I_{\bar{\boldsymbol{a}}_{t-1}}(\bar{\boldsymbol{A}}_{t-1})}{\lambda_{\bar{a}_{t-1}}(\boldsymbol{U}_{t-1},\boldsymbol{M}_{t-1};P)}\{b_{\boldsymbol{a}_t}(\boldsymbol{U}_{t};P)-b_{\boldsymbol{a}_{t-1}}(\boldsymbol{U}_{t-1};P)\}\right] \nonumber \\
&& \qquad \qquad -\sum_{t=0}^{T}var_P\left[\sum_{a \in \mathcal{\boldsymbol{A}}} c_a\dfrac{I_{\bar{\boldsymbol{a}}_{t-1}}(\bar{\boldsymbol{A}}_{t-1})}{\lambda_{\bar{a}_{t-1}}(\bar{\boldsymbol{U}}_{t-1};P)}\{b_{\boldsymbol{a}_t}(\boldsymbol{U}_{t};P)-b_{\boldsymbol{a}_{t-1}}(\boldsymbol{U}_{t-1};P)\}\right]\geq0. \nonumber 
\end{eqnarray}

Since
\begin{equation*}
    Y\indep_{\mathcal{G}} \boldsymbol{M}_T\mid  \ \boldsymbol{U}_T,{A}_T
\end{equation*}
and
\begin{equation*}
Q_{t-1}\indep_{\mathcal{G}} \boldsymbol{M}_{t-1}\mid  \ \boldsymbol{U}_{t-1},{A}_{t-1}, \  for \  t = T,\dots,1,
\end{equation*}
it follows that for every subset \(\boldsymbol{M}'_{t} \subseteq \boldsymbol{M}_t\), one also has 
\begin{equation*}
    Y\indep_{\mathcal{G}} \boldsymbol{M}'_T\mid  \ \boldsymbol{U}_T,{A}_T
\end{equation*}
and
\begin{equation*}
Q_{t-1}\indep_{\mathcal{G}} \boldsymbol{M}'_{t-1}\mid  \ \boldsymbol{U}_{t-1},{A}_{t-1}, \  for \  t = T,\dots,1,
\end{equation*}
that is, conditional independence is stable under taking subsets.  
Applying exactly the same reasoning with \(\boldsymbol{M}'_{t}\) in place of \(\boldsymbol{M}_t\), we obtain $\sigma^2_{\boldsymbol{a},\boldsymbol{U},\boldsymbol{M}}(P) \;\;\ge\;\;\sigma^2_{\boldsymbol{a},\boldsymbol{U},\boldsymbol{M}'}(P)
\;\;\ge\;\;\sigma^2_{\boldsymbol{a},\boldsymbol{U}}(P)
$ and $\sigma^2_{\Delta,\boldsymbol{U},\boldsymbol{M}}(P) \;\;\ge\;\;\sigma^2_{\Delta,\boldsymbol{U},\boldsymbol{M}'}(P)
\;\;\ge\;\;\sigma^2_{\Delta,\boldsymbol{U}}(P).
$ 

This concludes the proof of Lemma \ref{lemma:2}.

\end{document}